\documentclass[12pt,a4paper]{article}
\usepackage{hyperref}
\usepackage[pdftex]{graphicx}
\usepackage{amssymb,amsthm,amsmath,hyperref,txfonts}

\usepackage{mathtools}

\usepackage[dvips]{color}
\usepackage{graphicx,float,xcolor}
\usepackage{epstopdf}
\usepackage{lineno}
\usepackage{tikz}
\usetikzlibrary{positioning,decorations,arrows,patterns,shapes,backgrounds}

\date{}                                           

\newtheorem{theorem}{Theorem}[section]
\newtheorem{remark}{Remark}[section]

\newtheorem{lemma}[theorem]{Lemma}

\newtheorem {corollary}[theorem]{Corollary}
\newtheorem {proposition}[theorem]{Proposition}

\newcommand{\uu}{{\mathbf u}}

\def\XXint#1#2#3{{\setbox0=\hbox{$#1{#2#3}{%
				\int}$ }
		\vcenter{\hbox{$#2#3$ }}\kern-.6\wd0}}

\def\edc{\end{document}}
\numberwithin{equation}{section}

\allowdisplaybreaks

\begin{document}

\title{ Well-posedness for a diffuse interface model of non-Newtonian two-phase flows}
\author{Fang Li\thanks{E-mail: fangli@nwu.edu.cn}\thanks{Corresponding author}
        \thanks{Partly supported by the National Natural Science Foundation of China (No. 11501445). }
        \ ; Duan Xingyu\thanks{E-mail: duanxingyu@stumail.nwu.edu.cn} \\
	\textit{\small Department of Mathematics and CNS, Northwest University, Xi'an, P. R. China}\\
     Guo Zhenhua\thanks{E-mail: zhguo@gxu.edu.cn}
     \thanks{Partly supported by the National Natural Science Foundation of China (No. 11931013). }\\
   \textit{\small School of Mathematics and Information Science, Guangxi University,}
    \\\textit{\small   Nanning, P. R. China}
}
\maketitle
\begin{abstract}
\noindent
The evolution of two partially miscible, nonhomogeneous, incompressible viscous fluids
of non-Newtonian type, can be governed by the Navier-Stokes-Cahn-Hilliard system.  In the present work, we prove the global existence of weak solutions for the case of initial density containing zero and the concentration depending viscosity with free energy potential equal to the Landau potential  in a bounded domain of three dimensions. Furthermore, we show that a strong solutions exist locally in time in the case of three dimensions periodic domain ${\mathbb T}^3.$ The proof relies on a suitable semi-Galerkin scheme and the monotonicity method.
\end{abstract}


\section{Introduction}\label{1}

\quad The dynamics of a binary fluid mixture in presence of a surfactant is effectively modeled through the socalled diffuse-interface (or phase-field) approach (see Anderson-McFadden-Wheeler \cite{Anderson and McFadden-1998}). Within this framework, we consider a model for the flow of two macroscopically immiscible, nonhomogeneous incompressible non-Newtonian fluids. In contrast to classical sharp interface models, a partial mixing of the fluids is taken into account, which leads to a so-called diffuse interface model. This has the advantage that flows beyond the occurrence of topological singularities e.g. due to droplet collision or pinch-off can be described. More precisely, we consider the following Navier-Stokes-Cahn-Hilliard (NSCH) system
\begin{eqnarray}\label{E1-1}
	\begin{cases}
		\partial_{t}\rho+div(\rho {\mathbf u})=0\quad\quad\mbox{in}\;\Omega\times(0,T),\\
		\rho\partial_t {\mathbf u}+\rho {\mathbf u}\cdot \nabla {\mathbf u}-div(\nu(\phi)(1+|\mathbb{D}{\mathbf u}|^2)^{\frac{p-2}{2}}\mathbb{D}{\mathbf u})+\nabla P=-div(\nabla\phi\otimes\nabla\phi)\;\;\mbox{in}\;\Omega\times(0,T),\\
		div{\mathbf u}=0\quad\quad\quad\quad\quad\;\;\mbox{in}\;\Omega\times(0,T),\\
		\rho \partial_t \phi+\rho {\mathbf u}\cdot\nabla\phi=\Delta\mu\;\;\;\mbox{in }\;\Omega\times(0,T),\\
		\rho\mu=-\Delta\phi+\rho\Psi^\prime(\phi)\quad\;\mbox{in }\;\Omega\times(0,T),
	\end{cases}
\end{eqnarray}
where
\begin{eqnarray*}
\mathbb{D}{\mathbf u}=\frac12(\nabla{\mathbf u}+(\nabla{\mathbf u})^t),
\end{eqnarray*}
and $(\nabla{\mathbf u})^t$ is the transpose matrix of $\nabla{\mathbf u}$. Here, $T > 0$ is a given final time of arbitrary magnitude, $\Omega\subset {\mathbb R}^3$ is a bounded domain with sufficiently smooth boundary,
$t\geqslant 0$ is time, $x$ is the spatial coordinate, $\rho=\rho(x,t)$ is the total density of the mixture, ${\mathbf u}={\mathbf u}(x,t)$ is the mass-averaged velocity of the mixture, $P=P(x,t)$ is the pressure of the mixture, $\phi=\phi(x,t)$ is the difference of fluids concentrations, $\mu=\mu(x,t)$ is the chemical potential. In the NSCH system \eqref{E1-1}, the fluid is Newtonian when $p = 2,$ the fluid us shear thinning fluids when $p\in[1,2)$, the fluid captures shear thickening effects when $p>2.$ The system \eqref{E1-1} is subjected to the initial and boundary conditions
\begin{eqnarray}\label{E1-2}
\begin{cases}
{\mathbf u}=0 &\mbox{ on }\partial\Omega\times(0,T),\\
\partial_{\bf n}\phi=\partial_{\bf n}\mu=0&\mbox{ on }\partial\Omega\times(0,T),\\
(\rho{\mathbf u})_{t=0}=\mathbf{m}_0,~~\phi|_{t=0}=\phi_0, &\mbox{ in }\Omega.
\end{cases}
\end{eqnarray}
The notation $\partial_{n} f$ denotes the normal derivative of a function $f$ on the boundary with outer unit normal $\mathbf{n}.$ The viscosity coefficient $\nu$ is assumed to be $\nu=\nu(s)\in W^{1,\infty}({\mathbb R})$ satisfying
\begin{eqnarray}\label{E1-3}
0<\nu_*\leqslant \nu(s)\leqslant \nu^* (\forall s\in {\mathbb R}).
\end{eqnarray}
For the potential $\Psi(s),$ we consider the Landau potential
\begin{eqnarray}\label{E1-4}
\Psi(s)=\frac14(s^2-1)^2 (\forall s\in {\mathbb R}).
\end{eqnarray}

In order to describe a general two-phase flow with droplet formation and coalescence of several droplet, diffuse interface models were developed, which take a mixing of the two macroscopically immiscible fluids and a small mesoscopic length-scale into account (see Anderson-McFadden-Wheeler \cite{Anderson and McFadden-1998} on this topic). A paradigm model of the Diffuse Interface theory for two-phase flows is the NSCH system, also called Model H after the seminal work \cite{Hohenberg-Halperin-1977} on dynamic critical phenomena. Over the past years there have been important developments concerning on the mathematical modeling and analysis on the NSCH systems for binary mixtures. Proceeding along the historical course, the achievements on incompressible mixtures are summarized here.

For homogeneous incompressible mixtures, the original Model H describes incompressible binary mixtures whose total density is constant. It corresponds to the system \eqref{E1-1} with density $\rho=1$ and $p=2.$ Gurtin-Polignone-Vi$\tilde{n}$als derived the system within the framework of continuum mechanics in \cite{Gurtin-Polignone-Vinals-1996}. Liu-Shen derived it through an energetic variational approach in \cite{Liu-Shen-2003} (also see \cite{Giga-ShenKirshtein-Liu-2018}). Abels in \cite{Abels-2009} and Giorgini-Miranville-Temam in \cite{Giorgini-Miranville-Temam-2019} established the existence, uniqueness and regularity of weak and strong solutions for the case of the Flory-Huggins potential
$$\Psi(s)=\frac{\theta}{2}[(1+s)\ln(1+s)+(1-s)\ln(1-s)]-\frac{\theta_0}{2}s^2~(\forall s\in[-1,1]),$$
with $0<\theta<\theta_0.$ Gui-Li proved the global well-posedness of the Cauchy problem of the two-dimension incompressible NSCH system with periodic domain by using energy estimates and the logarithmic Sobolev inequality in \cite{Gui-2018}. Giorgini-Miranville-Temam showed the uniqueness of weak solutions and the existence and uniqueness of global strong solutions to the incompressible NSCH system in a two-dimensional  bounded smooth domain, where the initial velocity satisfied some condition in \cite{Giorgini-2019}. In the case of the Landau potential \eqref{E1-4} and its polynomial generalizations, the analysis has been carried out in \cite{Boyer-1999,Cao-2012,Gal-Grasselli-2010,Giorgini-Miranville-Temam-2019}. Under the condition of allowing diffusion coefficient degradation, Boyer proved the existence and uniqueness of solutions to the incompressible NSCH system under shear through order parameter formulation for two dimensions and three dimensions in \cite{Boyer-1999}. Cao-Gal obtained the global regularity of strong solutions for the NSCH system in two dimension with mixed partial viscosity and mobility in \cite{Cao-2012}, where the authors also proved the global existence and uniqueness of classical solutions. Lam-Wu derived a class of thermodynamically consistent NSCH system for two-phase fluid flows with density contrast, based on a volume averaged velocity in \cite{Lam-Wu-2018}. Moreover, the authors proved the existence of global weak solutions to the incompressible NSCH on a smooth bounded domain in two and three dimensions, and the existence of a unique global strong solution in two dimensions under some technical assumption on the coefficients. First analytic results on the system \eqref{E1-1} for incompressible non-Newtonian fluids of power-law type were obtained by Kim-Consiglieri-Rorigues in \cite{Kim-Consiglieri-Rodrigues-2006}, and the authors proved existence of weak solutions if $p\geqslant \frac{3d+2}{d+2}$ in the space dimension $d$ ($d=2,3$). Grasselli-Pra\v{z}ak discussed the longtime behavior of solutions of \eqref{E1-1}-\eqref{E1-2} or incompressible non-Newtonian fluids of power-law type for  $p>\frac{3d+2}{d+2}$ in the space dimension $d$ ($d=2,3$), assuming periodic boundary conditions and a regular free energy density in \cite{Grasselli-Prazak-2011}. For the same range of $p$ results on existence of weak solutions with a singular free energy density and the longtime behavior were obtained by Bosia in the case of a bounded domain in ${\mathbb R}^3$ in \cite{Bosia-2013}. Abels-Diening-Terasawa proved the existence of weak solutions for a phase field model for the flow of two partly miscible incompressible, viscous fluids of non-Newtonian (power law) type in \cite{Abels-Diening-Terasawa-2014} .

The nonhomogeneous incompressible mixtures, relies on the assumptions that the density of the mixture is an independent variable of the system and the mass-averaged velocity of the mixture is divergence-free. This lead to the system \eqref{E1-1} with $p=2,$ whose derivation is based on the conservation of mass and linear momentum and the second law of thermodynamic in the form of dissipation inequality as in \cite{2008-Abels-Feireisl,Heida-2012,Gurtin-Polignone-Vinals-1996}. Recently, Giorgini-Temam in \cite{Giorgini-2020} considered the case of initial density away from zero and concentration-depending viscosity with free energy potential equal to either the Landau potential or the logarithmic potential, proved the existence of weak solutions to the nonhomogeneous incompressible NSCH system in a two-dimension/three-dimension bounded smooth domain, and showed that the existence of strong solutions are local-in-time in three dimensions and global-in-time in two dimensions with bounded and strictly positive density. Fang-Nie-Guo in \cite{Fang-Nie-Guo-2025} established a blow-up criterion of local strong solution to the initial-boundary-value problem for the case of initial density away from zero, and obtained the global existence and the decay-in-time of strong solution provided that the initial datum $\|\nabla{\mathbf u}_0\|_{{\bf L}^{2}(\Omega)}+\|\nabla \mu_0\|_{{\bf L}^{2}(\Omega)}+\|\rho_0\|_{L^{\infty}(\Omega)}$ is suitably small.

To the best of our knowledge, there is no analytic result concerning the initial-boundary-value problem \eqref{E1-1}-\eqref{E1-4} in the literature. In this paper, we present a mathematical analysis concerning weak and strong solutions for the problem \eqref{E1-1}-\eqref{E1-4}. The analysis is inspired by the classical techniques for the nonhomogeneous incompressible non-Newtonian fluid. Our specific purpose is twofold. First, we aim to prove the existence of global weak solutions in the general setting with physically consistent non-constant viscosity parameter and Landau potential. The second object is to show the existence of strong solutions by maintaining the most general assumptions on the density and velocity.

Let $X$ be a (real) Banach or Hilbert space with norm denoted by $\| \cdot\|_X$. The boldface letter ${\mathbf X}$ stands for the vectorial space $X^d$ ($d$ is the spatial dimension), which consists of vector-valued functions $\uu$ with all components belonging to $X,$ with norm
$\|\cdot\|_{{\mathbf X}}.$ Given a measurable set $J\subset{\mathbb R},$ $L^q(J;X)$ with $q\in [1,+\infty)$ denotes the space of Bochner measurable $q-$integralbe/essentially bounded functions with values in the Banach space $X.$ If $J=(a,b)$ is an interval, we write for simplicity $L^q(a,b; X).$ For  $q\in [1,+\infty],$ $W^{1,q}(J;X)$ denotes the space of functions $f$ such that $f\in L^q(J;X)$ with $\partial_tf\in L^q(J;X),$ where $\partial_tf$ means the vector-valued distributional derivative of $f.$ When $q=2,$ we set $H^1(J,X):=W^{1,2}(J;X).$ Besides, $C(J;X)$ denotes the space of functions that are strong continuous from $J$ to $X,$ and $C_w(J;X_s)$ denotes the space of functions that are or weak continuous from $J$ to $X$ with strong topology.

For $s\in {\mathbb N}$ and $q\in{\mathbb Z}^+,$ $L^q(\Omega)$ denotes the Lebesgue spaces on $\Omega$ and $W^{s,q}(\Omega)$ denotes the Sobolev space of functions in $L^q(\Omega)$ with distributional derivatives of order less than or equal to $s$ in $L^q(\Omega),$ and its norm is denoted by $\| \cdot \|_{W^{s,q}(\Omega)}$. If $q=2$ and $s\in {\mathbb Z}^+,$ we use the standard notation $H^s(\Omega)$ for $W^{s,2}(\Omega).$ The Hilbert spaces of solenoidal vector-valued functions, are introduced as follows. ${\mathcal C}_{0,\sigma}^{\infty}(\Omega)$ denotes the space
of divergence free vector fields in ${\mathcal C}_{0}^{\infty}(\Omega)$. The spaces $\bf{H}_{\sigma}$ and $\bf{V}_{\sigma}$ are defined as the closure of ${\mathcal C}_{0,\sigma}^{\infty}(\Omega)$ with respect to the ${\bf L}^2(\Omega)$ and ${\bf H}_{0}^1(\Omega)$ norms, respectively.
For simplicity, the inner product in the basic space $L^2(\Omega)$ (as well as ${\bf L}^2(\Omega)$) is denoted by $(\cdot,\cdot).$

Our result on existence of global weak solutions is stated as follows.

\begin{theorem}	\label{thm2-1}
Let $\Omega\subset {\mathbb R}^3$ be a bounded domain with $\partial\Omega\in {\mathcal C}^2,$ $T$ be a positive time and $\in(\frac52,+\infty).$ Suppose that the initial data $(\rho_0,{\mathbf u}_0, \phi_0,\mu_0)$ satisfies
\begin{align*}
0\leqslant\rho_{0}\in L^{\infty}(\Omega),~\phi_{0}\in H^{1}(\Omega),~\sqrt{\rho_0}{\mathbf u}_{0}\in {\bf L}^{2}(\Omega),
~\nabla\mu_{0}\in {\bf L}^{2}(\Omega).
\end{align*}
Then the problem \eqref{E1-1}-\eqref{E1-4} admits a global weak solution $(\rho,{\mathbf u}, \phi,\mu)$ on $[0,T]$ satisfying the following regularity properties
\begin{align*}
	&\rho\in C([0,T];L^{r}(\Omega))\cap L^{\infty}(\Omega\times(0,T))\;(r\in[1,+\infty)),\;
    \partial_{t}\rho\in L^{\infty}(0,T;H^{-1}(\Omega)),\notag\\
	&\sqrt{\rho}{\mathbf u}\in {\bf L}^{\infty}(0,T;{\bf L}^{2}(\Omega)),\;{\mathbf u}\in {\bf L}^{p}([0,T],{\bf W}_{0}^{1,p}(\Omega)),\notag\\
	&\sqrt{\rho}\phi\in L^{\infty}(0,T;L^{2}(\Omega)),\;\phi\in L^{\infty}(0,T;H^{2}(\Omega))\cap C([0,T];H^{1}(\Omega)),\notag\\
	&\sqrt{\rho}\mu\in L^{\infty}(0,T;L^{2}(\Omega)),\;\mu\in L^{\infty}(0,T;H^{1}(\Omega))\cap L^{2}(0,T;H^{2}(\Omega)),
\end{align*}
and the following equalities hold
\begin{align*}
(i)~~&\int_{0}^{T}\int_{\Omega}\rho\partial_{t}\psi dxdt+\int_{0}^{T}\int_{\Omega}\rho {\mathbf u}\cdot\nabla\psi dxdt=0\\
&\mbox{ for all $\psi\in C_{c}^{\infty}(\Omega\times(0,T))$};\\
(ii)~~&-\int_{0}^{T}\int_{\Omega}\rho {\mathbf u}\cdot\partial_{t}{\bf w} dxdt
    -\int_{0}^{T}\int_{\Omega}\rho {\mathbf u}\otimes {\mathbf u}:\nabla {\bf w} dxdt\\
    &+\int_{0}^{T}\int_{\Omega}\nu(\phi)(1+|\mathbb{D}{\mathbf u}|^2)^{\frac{p-2}{2}}\mathbb{D}{\mathbf u}:\mathbb{D} {\bf w}dxdt\\
	&=\int_{\Omega}\rho_{0}{\mathbf u}_{0}\cdot {\bf w}(0)dx+\int_{0}^{T}\int_{\Omega}\nabla\phi\otimes\nabla\phi:\nabla {\bf w}dxdt\\
&\mbox{ for all ${\bf w}\in C_{c}^{1}([0,T);{\bf W}^{1,p}_{0}(\Omega))$};\\
(iii)~~&-\int_{0}^{T}\int_{\Omega}\rho\phi\partial_{t}wdxdt-\int_{0}^{T}\int_{\Omega}\rho {\mathbf u}\phi\cdot\nabla wdxdt+\int_{0}^{T}\int_{\Omega}\nabla\mu\cdot\nabla wdxdt=\int_{\Omega}\rho_{0}\phi_{0}w(0)dx\\
&\mbox{ for all $w\in C_{c}^{1}([0,T);H^{1}(\Omega))$;}\\
(iv)~~&\rho\mu=-\Delta\phi+\rho\Psi^\prime(\phi)\;\;\;a.e.\; in\;\Omega\times(0,T).
\end{align*}
Furthermore, $\partial_{n}\phi=0$ almost everywhere on $\partial\Omega\times(0,T).$
\end{theorem}

\begin{remark}
Due to the regularity of velocity of the mixture ${\mathbf u}$ and the Sobolev inequality for the two-spatial dimension, the global existence of weak solutions to the problem \eqref{E1-1}-\eqref{E1-4} can be shown by same argument in the proof of Theorem \ref{thm2-1}.
\end{remark}

\begin{remark}
When $p=2$, the initial density away from zero and concentration-depending viscosity with free energy potential equal to the Landau potential, the global existence of weak solutions in $\mathbb{R}^{d}$(d=2,3) is given in \cite{Giorgini-2020}.
\end{remark}

For more regularity on the initial data ${\mathbf u}_0$, we can improve the regularity of weak solutions to obtain the strong solution.

\begin{theorem}\label{thm2-2}
	Let $\Omega= \mathbb{T}^3,$ $T$ be a positive time and $\frac{5}{2}\leqslant p< 3.$ Suppose the initial data $(\rho_{0},{\mathbf u}_{0},\phi_{0},\mu_0)$ satisfy
\begin{align*}
&\rho_{0}\in L^{\infty}(\Omega),~0\leqslant\rho_{0}\leqslant\rho^*,~{\mathbf u}_{0}\in {\bf W}_{0}^{1,p}(\Omega),~
\nabla\mu_0\in {\bf L}^2(\Omega),~\phi_{0}\in H^{1}(\Omega).
\end{align*}
Then there exist a $T_{0}>0$ (depending only on initial data) and a strong solution $(\rho,{\mathbf u},\phi,\mu)$ to the problem \eqref{E1-1}-\eqref{E1-4} satisfying
	\begin{align*}
		&\rho\in C([0,T_{0}];L^{ r }(\Omega))\cap L^{\infty}(\Omega\times(0,T_{0})),\\
		&{\mathbf u}\in {\bf C}([0,T_{0}];{\bf L}^{q}(\Omega))\cap {\bf L}^{\infty}(0,T_{0};{\bf W}_{0}^{1,p}(\Omega))\cap {\bf L}^{2}(0,T_{0};{\bf H}^{2}(\Omega))\cap {\bf H}^{1}(0,T_{0};{\bf H}^{1}(\Omega)),\\
		&\phi\in C([0,T_{0}];W^{1,q}(\Omega))\cap L^{\infty}(0,T_{0};H^{2}(\Omega))\cap L^{2}(0,T_{0};W^{2,\infty}(\Omega)),\\
		&\mu\in L^{\infty}(0,T_{0};H^{1}(\Omega))\cap L^{2}(0,T_{0};H^{2}(\Omega)).
	\end{align*}
	for any $ r \in[1,\infty)$ and $q\in[1,6)$.
\end{theorem}

\begin{remark}
Under the condition of Theorem \ref{thm2-2}, the existence of a pressure $P$ with $$P \in L^{2}(0,T_{0};W^{1,\frac{2p}{3p-4}}(\Omega))$$ follows immediately from the Equations $(\ref{E1-1})_2$ and $(\ref{E1-1})_3$ by a classical consideration (see Corollary \ref{cor4-8} ).
\end{remark}

\begin{remark}
The uniqueness of strong solution $(\rho,{\mathbf u},\phi,\mu)$ in Theorem \ref{thm2-2} can not be demonstrated, due to low regularity of the initial density $\rho_0.$
\end{remark}

In the subsequent sections, the symbols $C, C_i$ stand for generic positive constants that may even change within the same line. Specific dependence of these constants in terms of the data will be pointed out if necessary. The remaining part of this paper is organized as follows. In Section \ref{S-2}, we introduce the functional settings and collect several preliminaries. In Section \ref{S-3}, we propose a semi-Galerkin scheme for a suitably regularized system, demonstrate its solvability, and complete the proof of Theorem \ref{thm2-1}. Section \ref{S-4} is devoted to the proof of Theorem \ref{thm2-2}. In the Appendix, we give the proof of Lemma 2.3 in detail and show the existence of approximated solutions.

\section{Preliminaries}\label{S-2}

\quad In this section, several preliminaries are presented. We begin with the following standard elliptic theory here, which can be found in \cite{Abels-2009}.

\begin{lemma}\label{lem2-3}
Let $\Omega\subset {\mathbb R}^3$ be a bounded domain with $C^2$-boundary. For $1\leqslant  q\leqslant  \infty,$ let
	$$L^{q}_{(m)}(\Omega)=\left\{f\in L^{q}(\Omega):m(f)=m\right\}\mbox{ with } m(f)=\frac{1}{|\Omega|}\int_{\Omega}f(x) dx. $$
If $u\in  W_{(0)}^{1,2}(\Omega)$, solves $\Delta u=f$ for some $f\in L_{(0)}^{q}(\Omega),1<q<\infty$ and $\partial\Omega$ is ${\cal C}^2$, and then it follows from standard elliptic theory that $u\in W^{2,q}(\Omega)$ and $\Delta u=f$ almost everywhere in $\Omega$ and $\partial_n u|_{\partial\Omega}=0$ in the sense of trace. If $f\in W^{1,q}(\Omega)\cap L^{q}_{(0)}(\Omega)$ and $\partial\Omega\in {\cal C}^3,$ then $u\in W^{3,q}(\Omega).$ Moreover,
	$$\|u\|_{W^{k+2,q}(\Omega)}\leqslant  C_q\|f\|_{W^{k,q}(\Omega)},$$
for all any $f\in W^{k,q}(\Omega)\cap L^{q}_{(0)}(\Omega) \ ( k=0,1),$ where the positive constant $C_q$ depends only on $1<q<\infty$ and $\Omega.$
\end{lemma}

The Gronwall's inequality is stated as follows, which plays a central role in proving a priori estimates on existence of strong solutions.

\begin{lemma}\label{gronwall1}(\cite{Fang-Nie-Guo-2025})(Gronwall's inequality (differential form))

Let $\eta(\cdot)$ be a nonnegative, absolutely continuous function on $[0,T]~(T>0)$, which satisfies for a.e. t the differential inequality
\begin{equation}
	\eta^{\prime}(t)\leqslant \phi(t)\eta(t)+\psi(t)\;a.e.\;t\in[0,T],
\end{equation}
where $\phi(t)$ and $\psi(t)$ are nonnegative and summable functions on $[0,T]$. Then
\begin{equation}
	\eta(t)\leqslant e^{\int_{0}^{t}\phi(s)ds}\left[\eta(0)+\int_{0}^{t}\psi(s)ds\right]
\end{equation}
for all $t\in[0,T].$ In particular, if $\eta^{\prime}(t)\leqslant \phi(t)\eta(t)$ on $[0,T]$ and $\eta(0)=0$, then $\eta\equiv0$ on $[0,T]$.
\end{lemma}

The following auxiliary result is a key to estimate the velocity of mixture, whose proof is similar to Lemma 3.2 in \cite{Feireisl-2004} and given in Appendix for completeness (see Lemma A.1. in Appendix).

\begin{lemma}\label{L-2-3}
Let $\Omega\subset {\mathbb R}^3$ be a bounded domain with $\partial\Omega\in {\mathcal C}^2,$ $1<q<+\infty\;and \;{\bf v}\in {\bf W}^{1,q}(\Omega)$. If $\rho$ is a non-negative function satisfying
\begin{eqnarray*}0<M<\int_{\Omega}\rho dx~~\mbox{ and }~~\int_{\Omega}\rho^{\gamma}dx\leqslant E_{0}\end{eqnarray*}
for some $\gamma>1$, then there exists a constant $C(M,E_{0})$ such that
	\begin{equation}
		\|{\bf v}\|^{q}_{{\bf L}^{q}(\Omega)}\leqslant C(M,E_{0})\left(\|\nabla {\bf v}\|^{q}_{{\bf L}^{q}(\Omega)}+(\int_{\Omega}\rho|{\bf v}|dx)^{q}\right).
	\end{equation}
\end{lemma}

The Gagliardo-Nirenberg inequality will be frequently used later, which is given here for the convenience of readers.

\begin{lemma}	\label{gn1}(\cite{Kim-Consiglieri-Rodrigues-2006})(Gagliardo-Nirenberg inequality).

Let $\Omega\subset\mathbb{R}^{d}$ be a bounded and sufficiently regular domain. For $1\leqslant q_{1},q_{2}\leqslant \infty\;and\;0\leqslant r\leqslant l$, when $ \frac{d}{q}-r=(1-\theta)\frac{d}{q_{1}}+\theta(\frac{n}{q_{2}}-l)\;(\frac{r}{l}\leqslant\theta\leqslant 1)$, the multiplicative inequality
	\begin{equation*}
		\|\nabla^{r}{\bf v}\|_{{\bf L}^q(\Omega)}
		\leqslant C\|{\bf v}\|_{{\bf L}^{q_{1}}(\Omega)}^{1-\theta}\|\nabla^{l}{\bf v}\|_{{\bf L}^{q_{2}}(\Omega)}^{\theta}
	\end{equation*}
	holds with the following exceptions:\\
	(i) if $r=0,l\leqslant \frac{d}{q_{2}}$, $q_{1}=\infty$ and $\Omega$ unbounded, we assume in addition that or $|{\bf v}|\rightarrow0$ as $|x|\rightarrow\infty$ or 	${\bf v}\in L^{k}(\Omega)$ for some $k>0$;\\
	(ii) if $1<q_{1}<\infty$ and $l-r-\frac{d}{q_{2}}$ is a non-negative integer, then does not hold for $\theta=1$.
\end{lemma}

Finally, we give a key tool to deal with convergence of the nonlinear term by using the monotonicity method.
Let ${\bf A}_{\varepsilon}(x,\xi):\ \Omega\times {\mathbb R}^3\rightarrow {\mathbb R}^3$ be Carath\'{e}odory vector functions.
These vector functions are assumed to satisfy the minimal monotonicity and convergence conditions
\begin{eqnarray*}
	&({\bf A}_{\varepsilon}(x,\xi)-{\bf A}_{\varepsilon}(x,\eta))\cdot(\xi-\eta)\geqslant 0,& {\bf A}_{\varepsilon}(x,0)\equiv 0,\\
	&|{\bf A}_{\varepsilon}(x,\xi)|\leqslant c_0(|\xi|)<\infty, &\lim\limits_{\varepsilon\rightarrow 0}{\bf A}_{\varepsilon}(x,\xi)=A(x,\xi)
\end{eqnarray*}
for a.e. $x\in \Omega$ and any $\xi,\eta\in {\mathbb R}^3.$ In fact, the minimal monotonicity of the Carath\'{e}odory vector function is the foundation for the monotonicity method.

\begin{lemma}\label{L2-6}(\cite{Zhikov-2010-1}) Let $\Omega\subset {\mathbb R}^3$ be a bounded domain with $\partial\Omega\in {\mathcal C}^2,$ ${\bf v}_\varepsilon \rightarrow {\bf v}, {\bf A}_{\varepsilon}(x,{\bf v}_\varepsilon)\rightarrow {\bf z}\ \mbox{in } {\bf L}^1(\Omega).$
Let $K\subset \Omega$ be a measurable set such that ${\bf z}\cdot {\bf v}\in {\bf L}^1(K).$ Then
	$$\lim\inf\limits_{\varepsilon\rightarrow 0}\int_{K}{\bf A}_{\varepsilon}(x,{\bf v}_\varepsilon)\cdot {\bf v}_\varepsilon dx \geqslant \int_{K}{\bf z}\cdot {\bf v}dx.$$
	In particular,
	$${\bf z}|_{K}={\bf A}|_{K},\ \ {\bf A}={\bf A}(x,{\bf v})$$
	when $\lim\inf\limits_{\varepsilon\rightarrow 0}\int_{K}{\bf A}_{\varepsilon}(x,{\bf v}_\varepsilon)\cdot {\bf v}_\varepsilon dx =\int_{K}{\bf z}\cdot {\bf v}dx.$
\end{lemma}

\section{Proof of Theorem \ref{thm2-1}}\label{S-3}

\subsection{The semi-Galerkin scheme}

\quad Given $\rho_{0}\in L^{\infty}(\Omega)$ with $0\leqslant \rho_0\leqslant\rho^*$ for almost every $x\in\Omega,$ there exist a sequence $\{\rho_{0\delta}\}$  such that
\begin{equation}\label{E-3-1}
	\rho_{0\delta}\in C^{\infty}(\overline{\Omega}),\;\;0<\delta\leqslant\rho_{0\delta}(x)\leqslant\rho^{*}+1\;\;(\forall x\in\overline{\Omega})
\end{equation}
and
\begin{equation}\label{E-3-2}
\rho_{0\delta}\rightarrow\rho_{0}\;strongly\;in\;L^{ r }(\Omega)~(\forall r \in[1,\infty)),\;\; \;
\rho_{0\delta}\rightarrow\rho_{0}\;weak-star\;in\;L^{\infty}(\Omega)
\end{equation}
as $\delta\rightarrow 0^+.$ Consider the family of eigenfunctions $\{w_{j}\}_{j=1}^{\infty}$ and eigenvalues $\{\lambda_{j}\}_{j=1}^{\infty}$ of the Laplace operator $A=-\Delta+I$ with homogeneous Neumann boundary condition and family of eigenfunctions $\{{\mathbf w}_{j}\}_{j=1}^{\infty}$ and eigenvalues $\{\lambda^{S}_{j}\}_{j=1}^{\infty}$ of the Stokes operator $A$. For any integer $m\geqslant 1$, the $m$-dimensional orthogonal projections $\Pi_{m}$ and $P_{m}$ are defined on $V_{m}$ and ${\mathbf V}_{m}$ with respect to the inner product in $L^{2}(\Omega)$ and in ${\mathbf H}_{\sigma},$ respectively. Let
\begin{equation}\label{E-3-3}
{\mathbf u}_{0m}=P_{m}{\mathbf u}_{0}~~~~~~~~\mbox{        and         }~~~~~~~~~~~~~~\phi_{0m}=\Pi_{m}\phi_{0}.
\end{equation}
then
\begin{equation*}
	{\mathbf u}_{0m}\rightarrow{\mathbf u}_{0}\;strongly\;in\;{\bf L}^{2}(\Omega)~~
\mbox{        and         }~~\phi_{0m}\rightarrow\phi_{0}\;strongly\;in\;H^{1}(\Omega)
\end{equation*}
respectively as $m\to+\infty.$ Then, we are looking for $(\rho_{m},{\mathbf u}_{m},\phi_{m},\mu_{m})$ to solve the following system
\begin{eqnarray}\label{E3-3}
\begin{cases}
\partial_{t}\rho_{m}+{\mathbf u}_{m}\cdot\nabla\rho_{m}=0&\mbox{in}\;\Omega\times(0,T),\\
(\rho_{m}\partial_{t}{\mathbf u}_{m},{\mathbf w})+(\rho_{m}{\mathbf u}_{m}\cdot\nabla{\mathbf u}_{m},{\mathbf w})
    +(\nu(\phi_{m})(1+|{\mathbb D}{{\mathbf u}}_{m}|^{2})^{\frac{p-2}{2}}\mathbb{D}{\mathbf u}_{m},\nabla {\mathbf w})\\
~~~~~~~~~~~~~~~~=(\rho_{m}\mu_{m}\nabla\phi_{m},{\mathbf w})-(\rho_{m}\nabla\Psi(\phi_{m}),{\mathbf w})~~~~({\bf w}\in {\bf V}_m)&\mbox{in}\;(0,T),\\
(\rho_{m}\partial_{t}\phi_{m},w)+(\rho_{m} {\mathbf u}_{m}\cdot\nabla\phi_{m},w)+(\nabla\mu_{m},\nabla w)=0~~~~(w\in V_m)&\mbox{in}\;(0,T)\\
(\rho_{m}\mu_{m},w)=(\nabla\phi_{m},\nabla w)+(\rho_m\Psi^\prime(\phi_{m}),w)~~~~(w\in V_m)&\mbox{in}\;(0,T),\\
{\mathbf u}_{m}=0,\partial_{n}\mu_{m}=\partial_{n}\phi_{m}=0 &\mbox{ on }\partial\Omega\times (0,T),\\
\rho_{m}(\cdot,0)=\rho_{0\delta},{\mathbf u}_{m}(\cdot,0)={\mathbf u}_{0m},\phi_{m}(\cdot,0)=\phi_{0m}&\mbox{ in }\Omega.
\end{cases}
\end{eqnarray}
Due to Proposition A.1.in Appendix,  the semi-Galerkin scheme \eqref{E3-3} admits a unique solution $(\rho_m,{\mathbf u}_m,\phi_m,\mu_m)$ on $[0,T]$ satisfying
\begin{eqnarray*}
\rho_m\in C([0,T]\times \overline{\Omega}),~~{\bf u}_m\in{\bf H}^1(0,T;{\bf V}_m),~~
\phi_m\in H^1(0,T;V_m),~~\mu_m\in C([0,T];V_m).
\end{eqnarray*}

\subsection{Uniform estimates}\label{subs-3-2}

\quad To establish existence results for weak solutions, we need to prove uniform estimates for approximate solutions $(\rho_m,{\bf u}_m,\phi_m,\mu_m)$ in suitable norms.

It is easy to get from $\eqref{E3-3}_1$ that there exists a positive constant $C,$ independent of $m$ and $\delta,$ such that
\begin{equation}\label{E3-5}
	\|\rho_{m}\|_{L^{\infty}(0,T,L^{\infty}(\Omega))}\leqslant C.
\end{equation}
Obviously, one finds that
\begin{eqnarray}\label{E-3-6}
   &&\frac{d}{dt}\int_{\Omega}\left(\frac{1}{2}\rho_{m}|{\mathbf u}_{m}|^{2}+\frac{1}{2}|\nabla\phi_{m}|^{2}+\rho_{m}\Psi^{\prime}(\phi_{m})\right)dx\nonumber\\
   &&+\int_{\Omega}\left(\nu({\phi}_{m})(1+|{\mathbb D}{\mathbf u}_m|^{2})^{\frac{p-2}{2}}|\mathbb{D}{\mathbf u}_{m}|^{2}+|\nabla\mu_{m}|^{2}\right)dx=0
\end{eqnarray}
holds for $t\in [0,T].$ Due to \eqref{E-3-1}-\eqref{E-3-2}, one finds that
\begin{eqnarray*}
&&\int_{\Omega}\left(\frac{1}{2}\rho_{0\delta}|{\mathbf u}_{0m}|^{2}+\frac{1}{2}|\nabla\phi_{0m}|^{2}+\rho_{0\delta}\Psi(\phi_{0m})\right)dx\rightarrow \int_{\Omega}\left(\frac{1}{2}\rho_{0}|{\mathbf u}_{0}|^{2}+\frac{1}{2}|\nabla\phi_{0}|^{2}+\rho_{0}\Psi(\phi_{0})\right)dx
\end{eqnarray*}
as $m\rightarrow+\infty\;and\;\delta\rightarrow 0^+$. Thus, there exists a positive constant $C,$ independent of $m$ and $\delta,$ such that
\begin{eqnarray}\label{E-3-7}
&&\int_{\Omega}\left(\frac{1}{2}\rho_{0\delta}|{\mathbf u}_{0m}|^{2}+\frac{1}{2}|\nabla\phi_{0m}|^{2}+\rho_{0\delta}\Psi(\phi_{0m})\right)dx\leqslant C.
\end{eqnarray}
Moreover, it is  deduced from \eqref{E-3-6} and \eqref{E-3-7} that
\begin{equation*}
	\int_{\Omega}\rho_{m}\Psi^{\prime}(\phi_{m})dx=\frac{1}{4}	\int_{\Omega}\rho_{m}(\phi_{m}^{2}-1)^{2}dx=\frac{1}{4}	\int_{\Omega}(\rho_{m}\phi_{m}^{4}-2\rho_{m}\phi_{m}^{2}+\rho_{m})dx\leqslant C.
\end{equation*}
So,
\begin{equation*}
\frac{1}{4}	\int_{\Omega}\rho_{m}\phi_{m}^{4}dx\leqslant C+\frac{1}{2}	\int_{\Omega}\rho_{m}\phi_{m}^{2}dx\leqslant\frac{1}{8}	\int_{\Omega}\rho_{m}\phi_{m}^{4}dx+C+C\int_{\Omega}\rho_{m}dx.
\end{equation*}
Hence,
\begin{equation*}
	\int_{\Omega}\rho_{m}\phi_{m}^{4}dx\leqslant C.
\end{equation*}

Summary, one can the following lemma with the help of Lemma \ref{L-2-3}.

\begin{lemma}\label{lemma3-2}
	Under the condition of Theorem \ref{thm2-1}, there exists a positive constant $C$, independent of $m$ and $\delta$, such that
	\begin{eqnarray}\label{E-3-8}
		\begin{cases}
			\|\rho_{m}\|_{L^{\infty}(0,T;L^{\infty}(\Omega))}\leqslant C,\\
			\|\sqrt{\rho_{m}}{\mathbf u}_{m}\|_{{\bf L}^{\infty}(0,T;{\bf L}^{2}(\Omega))}\leqslant C,~~
			\|(1+|\mathbb{D}{\mathbf u}_{m}|^{2})^{\frac{p-2}{4}}\mathbb{D}{\mathbf u}_{m}\|_{{\bf L}^{2}(0,T;{\bf L}^{2}(\Omega))}\leqslant C,\\
			\|\sqrt[4]{\rho_{m}}\phi_{m}\|_{L^{\infty}(0,T;L^{4}(\Omega))}\leqslant C,~~\|\phi_{m}\|_{L^{\infty}(0,T;H^{1}(\Omega))}\leqslant C,\\
			\|\nabla\mu_{m}\|_{{\bf L}^{2}(0,T;{\bf L}^{2}(\Omega))}\leqslant C.
		\end{cases}
	\end{eqnarray}
\end{lemma}

\begin{lemma}\label{lemma3-1}
	Under the condition of Theorem \ref{thm2-1}, there exists a positive constant $C$, independent of $m$ and $\delta$, such that
	\begin{equation}\label{E3-9}
		\|\partial_{t}\rho_{m}\|_{L^{\infty}(0,T;H^{-1}(\Omega))}\leqslant C,\;\;\;\;\;\;\|{\mathbf u}_{m}\|_{{\bf L}^{p}(0,T;{\bf W}_{0}^{1,p}(\Omega))}\leqslant C.
	\end{equation}
\end{lemma}

\begin{proof}
For any $ \eta\in H^{1}_{0}(\Omega)$ such that $\|\eta\|_{H^{1}_{0}(\Omega)}\leqslant1,$ it is deduced from $\eqref{E3-3}_1$ that
	\begin{align*}
		|\int_{\Omega}\partial_{t}\rho_{m} \eta dx|=|\int_{\Omega}\rho_{m}{\mathbf u}_{m}\cdot\nabla \eta dx|
		\leqslant\|\sqrt{\rho_{m}}{\mathbf u}_{m}\|_{{\bf L}^{2}(\Omega)}\|\sqrt{\rho_{m}}\|_{L^{\infty}(\Omega)}\|\nabla \eta\|_{{\bf L}^{2}(\Omega)}\leqslant C.
	\end{align*}
So, there exists a positive constant $C$, independent of $m$ and $\delta$, such that
	\begin{equation*}
		\|\partial_{t}\rho_{m}\|_{L^{\infty}(0,T;H^{-1}(\Omega))}\leqslant C.
	\end{equation*}

Obviously, it follows from \eqref{E-3-8} that
	\begin{equation*}
		\|\mathbb{D}{\mathbf u}_{m}\|_{{\bf L}^{p}(0,T;{\bf L}^{p}(\Omega))}\leqslant C\;\mbox{and}\;\|\mathbb{D}{\mathbf u}_{m}\|_{{\bf L}^{2}(0,T;{\bf L}^{2}(\Omega))}\leqslant C.
	\end{equation*}
Thus, one deduces from Lemma \ref{L-2-3} and \eqref{E-3-8} that
	\begin{equation*}
		\|{\mathbf u}_{m}\|_{{\bf L}^{p}(0,T;{\bf W}_{0}^{1,p}(\Omega))}\leqslant C.
	\end{equation*}
\end{proof}

\begin{lemma}\label{lemma3.4}
	Under the condition of Theorem \ref{thm2-1}, there exists a positive constant $C$, independent of $m$ and $\delta$, such that
	\begin{equation}\label{E-3-10}
		\|\mu_{m}\|_{L^{2}(0,T;H^{1}(\Omega))}\leqslant C,,\;\;\;\;\;\;\|\sqrt{\rho_{m}}\mu_{m}\|_{{ L}^{4}(0,T;{L}^{2}(\Omega))}\leqslant C.
	\end{equation}
\end{lemma}

\begin{proof}
First, it follows from $\eqref{E3-3}_4$ that
		\begin{align*}
			\int_{\Omega}\rho_{m}\mu^{2}_{m}dx&=\int_{\Omega}\nabla\phi_{m}\cdot\nabla\mu_{m} dx
			+\int_{\Omega}\rho_{m}\mu_{m}(\phi^{3}_{m}-\phi_{m})dx\\
			&\leqslant\|\nabla\phi_{m}\|_{{\bf L}^{2}(\Omega)}\|\nabla\mu_{m}\|_{{\bf L}^{2}(\Omega)}
			+\|\rho_{m}\|^{\frac{1}{2}}_{L^{\infty}(\Omega)}\|\sqrt{\rho_{m}}\mu_{m}\|_{L^{2}(\Omega)}(\|\phi^{3}_{m}\|_{L^{2}(\Omega)}+\|\phi_{m}\|_{L^{2}(\Omega)})\\
			&\leqslant\frac{1}{2}\|\sqrt{\rho_{m}}\mu_{m}\|^{2}_{L^{2}(\Omega)}+\|\nabla\phi_{m}\|^{2}_{{\bf L}^{2}(\Omega)}+\|\nabla\mu_{m}\|^{2}_{\mathbf{L}^{2}(\Omega)} +C(\|\phi_{m}\|^{6}_{H^{1}(\Omega)}+\|\phi_{m}\|^{2}_{L^{2}(\Omega)}),
		\end{align*}
which ensures that
		\begin{equation}\label{E-3-11}
			\|\sqrt{\rho_{m}}\mu_{m}\|_{L^{2}(0,T;L^{2}(\Omega))}\leqslant C.
		\end{equation}

Next, it follows from Lemma \ref{L-2-3} that
\begin{equation*}
\int_{0}^{T}\|\mu_{m}(\tau)\|^{2}_{L^{2}(\Omega)}d\tau\leqslant C\int_{0}^{T}[\|\nabla\mu_{m}(\tau)\|^{2}_{{\bf L}^{2}(\Omega)}+(\int_{\Omega}\rho_{m}|\mu_{m}|dx)^2] d\tau.
\end{equation*}
Since
	\begin{equation*}
	\int_{0}^{T}\int_{\Omega}\rho_{m}|\mu_{m}|dxd\tau
     \leqslant\int_{0}^{T}\|\sqrt{\rho_{m}}\|_{L^{2}(\Omega)}\|\sqrt{\rho_{m}}\mu_{m}\|_{L^{2}(\Omega)}d\tau\leqslant C
		\end{equation*}
due to \eqref{E-3-11}, one deduces from \eqref{E-3-8} that there exists a positive constant $C$, independent of $m$ and $\delta$, such that again, one can get that
\begin{equation*}
\|\mu_{m}\|_{L^{2}(0,T;H^{1}(\Omega))}\leqslant C.
\end{equation*}

Furthermore, one deduce from \eqref{E-3-8} that
		\begin{align}\label{E3-12}
			\int_{\Omega}\rho_{m}|\mu_{m}|^2 dx
			&=\int_{\Omega}\nabla\phi_{m}\cdot\nabla\mu_{m}dx+\int_{\Omega}\rho_{m}\Psi^\prime(\phi_{m})\mu_{m}dx\notag\\
			&\leqslant\|\nabla\phi_{m}\|_{{\bf L}^{2}(\Omega)}\|\nabla\mu_{m}\|_{{\bf L}^{2}(\Omega)}
			+\|\sqrt{\rho_{m}}\|_{L^{\infty}(\Omega)}\|\sqrt{\rho_{m}}\mu_{m}\|_{L^{2}(\Omega)}\|\phi^{3}_{m}-\phi_{m}\|_{L^{2}(\Omega)}\notag\\
			&\leqslant\|\nabla\phi_{m}\|_{{\bf L}^{2}(\Omega)}\|\nabla\mu_{m}\|_{{\bf L}^{2}(\Omega)}
			+C\|\sqrt{\rho_{m}}\mu_{m}\|_{L^{2}(\Omega)}(\|\phi_{m}\|^{3}_{L^{6}(\Omega)}+\phi_{m}\|_{L^{2}(\Omega)})\notag\\
			&\leqslant C\|\nabla\mu_{m}\|_{{\bf L}^{2}(\Omega)}+C\|\sqrt{\rho_{m}}\mu_{m}\|_{L^{2}(\Omega)}.
		\end{align}
One gets that
\begin{equation*}
\int_{0}^{T}\|\sqrt{\rho_{m}}\mu_{m}\|^{4}_{L^{2}(\Omega)}dt\leqslant C\int_{0}^{T}\|\nabla\mu_{m}\|^{2}_{{\bf L}^{2}(\Omega)}dt
+C\int_{0}^{T}\|\sqrt{\rho_{m}}\mu_{m}\|^{2}_{L^{2}(\Omega)}dt\leqslant C.
\end{equation*}
Thus, it follows from \eqref{E-3-8} and \eqref{E-3-11} that
		\begin{equation*}
			\|\sqrt{\rho_{m}}\mu_{m}\|_{L^{4}(0,T;L^{2}(\Omega))}\leqslant C.
		\end{equation*}
	\end{proof}

\begin{lemma}
	Under the condition of Theorem \ref{thm2-1}, there exists a positive constant $C$, independent of $m$ and $\delta$, such that
	\begin{equation}\label{E-3-13}
		\|\phi_{m}\|_{L^{4}(0,T;H^{2}(\Omega))}\leqslant C.
	\end{equation}
\end{lemma}

\begin{proof}
One takes $w=-\Delta\phi_{m}$ in $\eqref{E3-3}_4$ and using \eqref{E-3-8}, one arrives at
\begin{align*}
	&\int_{\Omega}|\Delta\phi_{m}|^{2}dx=\int_{\Omega}\left[-\rho_{m}\mu_{m}\Delta\phi_{m}+\rho_{m}\Psi^\prime(\phi_{m})\Delta\phi_{m}\right]dx\\
    &\leqslant\|\Delta\phi_{m}\|_{L^{2}(\Omega)}\|\rho_{m}\|_{L^{\infty}(\Omega)}(\|\phi_{m}\|^{3}_{L^{2}(\Omega)}
    +\|\phi_{m}\|_{L^{2}(\Omega)})+\|\Delta\phi_{m}\|_{L^{2}(\Omega)}\|\sqrt{\rho_{m}}\|_{L^{\infty}(\Omega)}\|\sqrt{\rho_{m}}\mu_{m}\|_{L^{2}(\Omega)}\\
	&\leqslant C\|\Delta\phi_{m}\|_{L^{2}(\Omega)}(\|\phi_{m}\|^{3}_{L^{2}(\Omega)}+\|\phi_{m}\|_{L^{2}(\Omega)})
    +C\|\Delta\phi_{m}\|_{L^{2}(\Omega)}\|\sqrt{\rho}\mu_{m}\|_{L^{2}(\Omega)}\\
	&\leqslant\frac12\|\Delta\phi_{m}\|^{2}_{L^{2}(\Omega)}+C(\|\phi_{m}\|^{6}_{H^{1}(\Omega)}+\|\phi_{m}\|^{2}_{L^{2}(\Omega)}+\|\sqrt{\rho_{m}}\mu_{m}\|^{2}_{L^{2}(\Omega)})\\
	&\leqslant\frac12\|\Delta\phi_{m}\|^{2}_{L^{2}(\Omega)}+C(1+\|\sqrt{\rho_{m}}\mu_{m}\|^{2}_{L^{2}(\Omega)}),
\end{align*}
which implies that
\begin{equation}\label{E3-14}
	\|\Delta\phi_{m}\|^{2}_{L^{2}(\Omega)}\leqslant C(1+\|\sqrt{\rho_{m}}\mu_{m}\|^{2}_{L^{2}(\Omega)}).
\end{equation}
So, it follows \eqref{E-3-10} that
\begin{equation*}
	\|\nabla^{2}\phi_{m}\|_{{\bf L}^{4}(0,T;{\bf L}^{2}(\Omega))}\leqslant C.
\end{equation*}
Furthermore, it is deduced from Lemma \ref{lem2-3} that
\begin{equation*}
	\|\phi_{m}\|_{L^{4}(0,T;H^{2}(\Omega))}\leqslant C.
\end{equation*}
\end{proof}

\begin{lemma}\label{lemmma3-7}
	Under the condition of Theorem \ref{thm2-1}, there exists a positive constant $C$, independent of $m$ and $\delta$, such that
	\begin{equation}\label{E3-15}
\|\sqrt{\rho_{m}}\mu_{m}\|_{L^{\infty}(0,T;L^{2}(\Omega))}\leqslant C,\;\|\sqrt{\rho_{m}}\partial_t\phi_{m}\|_{L^{2}(0,T;L^{2}(\Omega))}\leqslant C,\;	\|\phi_{m}\|_{L^{\infty}(0,T;H^{2}(\Omega))}\leqslant C.
	\end{equation}
\end{lemma}

\begin{proof}	
 First, one selects $w=\partial_t\mu_{m}$ in $\eqref{E3-3}_4$ to get that
\begin{equation*}
	\frac{1}{2}\int_{\Omega}\rho_{m}(\partial_{t}\mu_{m})^{2}_{m}dx=\int_{\Omega}\left[\nabla\phi_{m}\cdot\nabla\partial_{t}\mu_{m}+\rho_{m}\Psi^\prime(\phi_{m}\partial_{t}\mu_{m})\right]dx.
\end{equation*}
Since
\begin{equation*}
		\frac{1}{2}\int_{\Omega}\partial_{t}\rho_{m}\mu^{2}_{m}dx=	\frac{1}{2}\int_{\Omega}-div(\rho_{m}{\mathbf u}_{m})\mu^{2}_{m}dx=\int_{\Omega}\rho_{m}\mu_{m}{\mathbf u}_{m}\cdot\nabla\mu_{m}dx,
\end{equation*}
one gets that
\begin{equation}\label{E3-19}
\frac{1}{2}\frac{d}{dt}\int_{\Omega}\rho_{m}\mu^{2}_{m}dx
=\int_{\Omega}\left(\nabla\phi_{m}\cdot\nabla\partial_{t}\mu_{m}+\rho_{m}\Psi^\prime(\phi_{m})\partial_{t}\mu_{m}
+\rho_{m}\mu_{m}\nabla\mu_{m}\cdot{\mathbf u}_{m}\right)dx.
\end{equation}
Selecting $w=\mu_{m}$ in $\eqref{E3-3}_4$ and taking the time derivative of the resulting equation, one finds that
\begin{align}\label{E3-20}
\frac{d}{dt}\int_{\Omega}\rho_{m}\mu^{2}_{m}dx&=\int_{\Omega}\nabla\partial_{t}\phi_{m}\cdot\nabla\mu_{m}dx+\int_{\Omega}\nabla\phi_{m}\cdot\nabla\partial_{t}\mu_{m}dx\\
&+\int_{\Omega}\partial_{t}\rho_{m}\mu_{m}\Psi^\prime(\phi_{m})dx+\int_{\Omega}\rho_{m}\partial_{t}\mu_{m}\Psi^\prime(\phi_{m})dx+\int_{\Omega}\rho_{m}\mu_{m}\Psi^{''}(\phi_{m})\partial_{t}\phi_{m}dx\notag\\
&=\int_{\Omega}\nabla\partial_{t}\phi_{m}\cdot\nabla\mu_{m}dx+\int_{\Omega}\nabla\phi_{m}\cdot\nabla\partial_{t}\mu_{m}dx+\int_{\Omega}\rho_{m}{\mathbf u}_{m}\cdot\nabla\mu_{m}\Psi^\prime(\phi_{m})dx\notag\\
&+\int_{\Omega}\rho_{m}\mu_{m}{\mathbf u}_{m}\cdot\nabla\Psi^\prime(\phi_{m})dx+\int_{\Omega}\rho_{m}\partial_{t}\mu_{m}\Psi^\prime(\phi_{m})dx+\int_{\Omega}\rho_{m}\mu_{m}\Psi^{''}(\phi_{m})\partial_{t}\phi_{m}dx.\notag
\end{align}
Taking $w=\partial_t\phi_{m}$ in $\eqref{E3-3}_3$, one obtains that
\begin{equation}\label{E3-21}
	\int_{\Omega}\left(\rho_{m}\partial_t|\phi_{m}|^{2}+\rho_{m}\partial_{t}\phi_{m}{\mathbf u}_{m}\cdot\nabla\phi_{m}+\nabla\mu_{m}\cdot\nabla\partial_{t}\phi_{m}\right)dx=0.
\end{equation}
Combining \eqref{E3-19}-\eqref{E3-21}, one can arrive at
\begin{align*}
	&\frac{1}{2}\frac{d}{dt}\int_{\Omega}\rho_{m}\mu^{2}_{m}dx+\int_{\Omega}\rho_{m}(\partial_t\phi_{m})^{2}dx\\
   &=\underset{I_{1}}{\underbrace{\int_{\Omega}-\rho_{m}\mu_{m}{\bf u}_{m}\cdot\nabla\mu_{m} dx}}
   +\underset{I_{2}}{\underbrace{\int_{\Omega}-\rho_{m} {\mathbf u}_{m}\cdot\nabla\phi_{m}\partial_t\phi_{m}dx}}\\
	&+\underset{I_{3}}{\underbrace{\int_{\Omega}3\rho_{m}\mu_{m}\phi_{m}^{2}{\mathbf u}_{m}\cdot\nabla\phi_{m} dx}}
    -\underset{I_{4}}{\underbrace{\int_{\Omega}\rho_{m}\mu_{m}{\mathbf u}_{m}\cdot\nabla\phi_{m} dx}}
    +\underset{I_{5}}{\underbrace{\int_{\Omega}\rho_{m}\phi_{m}^{3}{\mathbf u}_{m}\cdot\nabla\mu_{m}dx}}\\
	&-\underset{I_{6}}{\underbrace{\int_{\Omega}\rho_{m}\phi_{m} {\mathbf u}_{m}\cdot\nabla\mu_{m} dx}}
      +\underset{I_{7}}{\underbrace{3\int_{\Omega}\rho_{m}\phi_{m}^{2}\mu_{m}\partial_t\phi_{m}dx}}
      -\underset{I_{8}}{\underbrace{\int_{\Omega}\rho_{m}\mu_{m}\partial_t\phi_{m}dx}}.
\end{align*}
Now, one estimates each term $I_j$ for the case of $\frac52<p< 3$ as follows
\begin{align*}
&|I_{1}|=|\int_{\Omega}-\rho_{m}\mu_{m}{\bf u}_{m}\cdot\nabla\mu_{m} dx|\\
&\;\;\;\;\leqslant\int_{\Omega}|\sqrt{\rho_{m}}\mu_{m}|^{\theta}|\rho_{m}|^{1-\frac{\theta}{2}}
          |\mu_{m}|^{1-\theta}|{\mathbf u}_{m}||\nabla\mu_{m}| dx\\
&\;\;\;\;\leqslant C\|{\mathbf u}_{m}\|_{{\bf L}^{\frac{3p}{3-p}}(\Omega)}\|\sqrt{\rho_{m}}\mu_{m}\|^{\frac{2p-3}{p}}_{L^{2}(\Omega)}
	\|\mu_{m}\|^{\frac{3-p}{p}}_{L^{6}(\Omega)}\|\nabla\mu_{m}\|_{{\bf L}^{2}(\Omega)}\\
&\;\;\;\;\leqslant C\|{\mathbf u}_{m}\|_{{\bf W}_{0}^{1,p}(\Omega)}\|\sqrt{\rho_{m}}\mu_{m}\|^{\frac{3p}{3-p}}_{L^{2}(\Omega)}
     \|\mu_{m}\|^{\frac{3}{p}}_{H^{1}(\Omega)}\\
&\;\;\;\;\leqslant C\|{\mathbf u}_{m}\|^{p}_{{\bf W}_{0}^{1,p}(\Omega)}
    +C\|\sqrt{\rho_{m}}\mu_{m}\|^{\frac{2p-3}{p-1}}_{L^{2}(\Omega)}\|\mu_{m}\|^{\frac{3}{p-1}}_{H^{1}(\Omega)}\\
&\;\;\;\;\leqslant C\|{\mathbf u}_{m}\|^{p}_{{\bf W}_{0}^{1,p}(\Omega)}
	+C\|\sqrt{\rho}\mu_{m}\|^{\frac{2p-3}{p-1}}_{L^{2}(\Omega)}\|\mu_{m}\|^{\frac{2p-3}{p-1}}_{H^{1}(\Omega)}
   \|\mu_{m}\|^{\frac{2(3-p)}{p-1}}_{H^{1}(\Omega)}\\
	&\;\;\;\;\leqslant C\|{\mathbf u}_{m}\|^{p}_{{\bf W}_{0}^{1,p}(\Omega)}
    +C\|\sqrt{\rho_{m}}\mu_{m}\|^{2}_{L^{2}(\Omega)}\|\mu_{m}\|^{2}_{H^{1}(\Omega)}+C\|\mu_{m}\|^{2}_{H^{1}(\Omega)},\\
&|I_{2}|=|\int_{\Omega}-\rho_{m} {\mathbf u}_{m}\cdot\nabla\phi_{m}\partial_t\phi_{m}dx|\\
	&\;\;\;\;\leqslant\int_{\Omega}(\sqrt{\rho_{m}}|{\mathbf u}_{m}|)^{\frac{4p-6}{5p-6}}\rho_{m}^{\frac{p}{2(5p-6)}}|\nabla\phi_{m}|
	|\sqrt{\rho_{m}}\partial_t\phi_{m}||{\mathbf u}_{m}|^{\frac{p}{5p-6}}dx\\
	&\;\;\;\;\leqslant\|\sqrt{\rho_{m}}{\mathbf u}_{m}\|^{\theta}_{{\bf L}^{2}(\Omega)}
    \|{\mathbf u}_{m}\|^{1-\theta}_{{\bf L}^{\frac{3p}{3-p}}(\Omega)}
	\|\sqrt{\rho_{m}}\partial_t\phi_{m}\|_{L^{2}(\Omega)}\|\nabla\phi_{m}\|_{{\bf L}^{6}(\Omega)}\\
	&\;\;\;\;\leqslant\|\sqrt{\rho_{m}}{\mathbf u}_{m}\|^{\frac{4p-6}{5p-6}}_{{\bf L}^{2}(\Omega)}
      \|{\mathbf u}_{m}\|^{\frac{p}{5p-6}}_{{\bf W}_{0}^{1,p}(\Omega)}
	\|\sqrt{\rho_{m}}\partial_t\phi_{m}\|_{L^{2}(\Omega)}\|\phi_{m}\|_{H^{2}(\Omega)}\\
	&\;\;\;\;\leqslant\frac{1}{16}\|\sqrt{\rho_{m}}\partial_t\phi_{m}\|^{2}_{L^{2}(\Omega)}
    +C\left(\|\sqrt{\rho_{m}}{\mathbf u}_{m}\|^{\frac{8(2p-3)}{5p-10}}_{{\bf L}^{2}(\Omega)}
    +\|\phi_{m}\|^{4}_{H^{2}(\Omega)}+\|{\mathbf u}_{m}\|^{p}_{{\bf W}_{0}^{1,p}(\Omega)}\right),\\
	&|I_{3}|=|\int_{\Omega}3\rho_{m}\mu_{m}\phi_{m}^{2}{\mathbf u}_{m}\cdot\nabla\phi_{m} dx|\\
	&\;\;\;\;\leqslant\int_{\Omega}(\sqrt{\rho_{m}}|{\mathbf u}_{m}|)^{\frac{4p-6}{5p-6}}\rho_{m}^{\frac{3p-3}{5p-6}}
     |{\mathbf u}_{m}|^{\frac{p}{5p-6}}| |\mu_{m}||\phi_{m}^{2}||\nabla\phi_{m}|dx\\
	&\;\;\;\;\leqslant C \|\sqrt{\rho_{m}}{\mathbf u}_{m}\|^{\frac{4p-6}{5p-6}}_{{\bf L}^{2}(\Omega)}
	\|{\mathbf u}_{m}\|^{\frac{p}{5p-6}}_{{\bf L}^{\frac{3p}{3-p}}(\Omega)}\|\mu_{m}\|_{L^{6}(\Omega)}
     \|\phi^{2}_{m}\|_{L^{3}(\Omega)}\|\nabla\phi_{m}\|_{{\bf L}^{6}(\Omega)}\\
	&\;\;\;\;\leqslant C \|\sqrt{\rho_{m}}{\mathbf u}_{m}\|^{\frac{4p-6}{5p-6}}_{{\bf L}^{2}(\Omega)}
	\|{\mathbf u}_{m}\|^{\frac{p}{5p-6}}_{{\bf L}^{\frac{3p}{3-p}}(\Omega)}\|\mu_{m}\|_{H^{1}(\Omega)}
     \|\phi_{m}\|^{2}_{H^{1}(\Omega)}\|\phi_{m}\|_{H^{2}(\Omega)}\\
	&\;\;\;\;\leqslant C\left(\|\sqrt{\rho_{m}}{\mathbf u}_{m}\|^{\frac{8(2p-3)}{5p-10}}_{{\bf L}^{2}(\Omega)}
    +\|{\mathbf u}_{m}\|^{p}_{{\bf W}_{0}^{1,p}(\Omega)}
     +\|\mu_{m}\|^{2}_{H^{1}(\Omega)}+\|\phi_{m}\|^{4}_{H^{2}(\Omega)}\right),\\
&|I_{4}|=|\int_{\Omega}\rho_{m}\mu_{m}{\mathbf u}_{m}\cdot\nabla\phi_{m} dx|\\
	&\;\;\;\;\leqslant\|\nabla\phi_{m}\|_{{\bf L}^{6}(\Omega)}\|{\mathbf u}_{m}
      \|_{{\bf L}^{3}(\Omega)}\|\sqrt{\rho_{m}}\mu_{m}\|_{L^{2}(\Omega)}\|\sqrt{\rho_{m}}\|_{L^{\infty}(\Omega)}\\
	&\;\;\;\;\leqslant C\left(\|\sqrt{\rho_{m}}\mu_{m}\|^{2}_{L^{2}(\Omega)}\|\phi_{m}\|^{2}_{H^{2}(\Omega)}
     +\|{\mathbf u}_{m}\|^{p}_{{\bf W}_{0}^{1,p}(\Omega)}\right),\\
&|I_{5}|=|\int_{\Omega}\rho_{m}\phi_{m}^{3}{\mathbf u}_{m}\cdot\nabla\mu_{m}dx|\\
	&\;\;\;\;\leqslant\int_{\Omega}(\sqrt{\rho_{m}}
     |{\mathbf u}_{m}|)^{\frac{10p^{2}-37p+30}{5(5p-6)(p-2)}}|\rho_{m}|^{\frac{40p^{2}+123p+90}{10(5p-6)(p-2)}}
     |\phi^{3}_{m}||\nabla\mu_{m}||{\mathbf u}_{m}|^{\frac{15p^{2}-43p+30}{5(5p-6)(p-2)}}dx\\
	&\;\;\;\;\leqslant C\|\nabla\mu_{m}\|_{{\bf L}^{2}(\Omega)}\|\phi_{m}^{3}\|_{L^{\frac{15p(p-2)}{30p^{2}-86p+60}}(\Omega)}
     \|\sqrt{\rho_{m}}{\mathbf u}_{m}\|^{\frac{10p^{2}-37p+30}{5(5p-6)(p-2)}}_{{\bf L}^{2}(\Omega)}
	\|{\mathbf u}_{m}\|^{\frac{15p^{2}-43p+30}{5(5p-6)(p-2)}}_{{\bf L}^{\frac{3p}{3-p}}(\Omega)}\\
	&\;\;\;\;\leqslant C\|\nabla\mu_{m}\|_{{\bf L}^{2}(\Omega)}\|\phi_{m}\|^{3}_{L^{\frac{45p(p-2)}{30p^{2}-86p+60}}(\Omega)}
    \|\sqrt{\rho_{m}}{\mathbf u}_{m}\|^{\frac{10p^{2}-37p+30}{5(5p-6)(p-2)}}_{{\bf L}^{2}(\Omega)}
	\|{\mathbf u}_{m}\|^{\frac{15p^{2}-43p+30}{5(5p-6)(p-2)}}_{{\bf W}_{0}^{1,p}(\Omega)}\\
	&\;\;\;\;\leqslant C\left(\|\mu_{m}\|^{2}_{H^{1}(\Omega)}+\|\phi_{m}\|^{\frac{60p(p-2)}{13p-30}}_{L^{\frac{45p(p-2)}{30p^{2}-86p+60}}(\Omega)}
	+\|\sqrt{\rho_{m}}{\mathbf u}_{m}\|^{\frac{4p}{5p-10}}_{{\bf L}^{2}(\Omega)}+\|{\mathbf u}_{m}\|^{p}_{{\bf W}_{0}^{1,p}(\Omega)}\right)\\
	&\;\;\;\;\leqslant C\left(\|\phi_{m}\|^{4}_{H^{2}(\Omega)}+\|\mu_{m}\|^{2}_{H^{1}(\Omega)}
    +\|{\mathbf u}_{m}\|^{p}_{{\bf W}_{0}^{1,p}(\Omega)}+\|\sqrt{\rho_{m}}{\mathbf u}_{m}\|^{\frac{4p}{5p-10}}_{{\bf L}^{2}(\Omega)}\right),\\
&|I_{6}|=|\int_{\Omega}\rho_{m}\phi_{m} {\mathbf u}_{m}\cdot\nabla\mu_{m} dx|\\
	&\;\;\;\;\leqslant\|\nabla\mu_{m}\|_{{\bf L}^{2}(\Omega)}\|\phi_{m}\|_{L^{6}(\Omega)}
    \|{\mathbf u}_{m}\|_{{\bf L}^{3}(\Omega)}\|\rho_{m}\|_{L^{\infty}(\Omega)}\\
	&\;\;\;\;\leqslant C\|\mu_{m}\|_{H^{1 }(\Omega)}\|{\mathbf u}_{m}\|_{_{{\bf W}_{0}^{1,p}}}\\
	&\;\;\;\;\leqslant C\left(\|\mu_{m}\|^{2}_{H^{1}(\Omega)}+\|{\mathbf u}_{m}\|^{2}_{{\bf W}_{0}^{1,p}}\right),\\
&|I_{7}|=|3\int_{\Omega}\rho_{m}\phi_{m}^{2}\mu_{m}\partial_t\phi_{m}dx|\\
	&\;\;\;\;\leqslant\|\sqrt{\rho_{m}}\partial_t\phi_{m}\|_{L^{2}(\Omega)}\|\mu_{m}\|_{L^{6}(\Omega)}\|\phi^{2}_{m}\|_{L^{3}(\Omega)}\\
	&\;\;\;\;\leqslant C\|\sqrt{\rho_{m}}\partial_t\phi_{m}\|_{L^{2}(\Omega)}\|\mu_{m}\|_{H^{1}(\Omega)}\\
	&\;\;\;\;\leqslant\frac{1}{16}\|\sqrt{\rho_{m}}\partial_t\phi_{m}\|^{2}_{L^{2}(\Omega)}+C\|\mu_{m}\|^{2}_{H^{1}(\Omega)},\\
&|I_{8}|=|\int_{\Omega}\rho_{m}\mu_{m}\partial_t\phi_{m}dx|\\
	&\;\;\;\;\leqslant\|\sqrt{\rho_{m}}\partial_t\phi_{m}\|_{L^{2}(\Omega)}\|\sqrt{\rho_{m}}\|_{L^{\infty}(\Omega)}
   \|\mu_{m}\|_{L^{6}(\Omega)}\|\phi_{m}\|_{L^{3}(\Omega)}\\
	&\;\;\;\;\leqslant C\|\sqrt{\rho_{m}}\partial_t\phi_{m}\|_{L^{2}(\Omega)}\|\mu_{m}\|_{H^{1}(\Omega)}\\
	&\;\;\;\;\leqslant\frac{1}{16}\|\sqrt{\rho_{m}}\partial_t\phi_{m}\|^{2}_{L^{2}(\Omega)}+C\|\mu_{m}\|^{2}_{H^{1}(\Omega)}.
\end{align*}
For the case of $p> 3$, one estimates $I_{1}-I_{3}$ and $I_{5}$ different as follows
\begin{align*}
&|I_{1}|=|\int_{\Omega}-\rho_{m}\mu_{m}{\bf u}_{m}\cdot\nabla\mu_{m} dx|\\
	&\;\;\;\;\leqslant \|\sqrt{\rho_{m}}\|_{L^{\infty}(\Omega)}
    \|{\mathbf u}_{m}\|_{{\bf L}^{\infty}(\Omega)}\|\sqrt{\rho_{m}}\mu_{m}\|_{L^{2}(\Omega)}
	\|\nabla\mu_{m}\|_{{\bf L}^{2}(\Omega)}\\
    &\;\;\;\;\leqslant C\|{\mathbf u}_{m}\|_{{\bf W}_{0}^{1,p}(\Omega)}\|\sqrt{\rho_{m}}\mu_{m}\|_{L^{2}(\Omega)}
    \|\mu_{m}\|_{H^{1}(\Omega)}\\
	&\;\;\;\;\leqslant C\|{\mathbf u}_{m}\|^{p}_{{\bf W}_{0}^{1,p}(\Omega)}
    +C\|\sqrt{\rho_{m}}\mu_{m}\|^{2}_{L^{2}(\Omega)}\|\mu_{m}\|^{2}_{H^{1}(\Omega)},\\
&|I_{2}|=|\int_{\Omega}-\rho_{m} {\mathbf u}_{m}\cdot\nabla\phi_{m}\partial_t\phi_{m}dx|\\
	&\;\;\;\;\leqslant\|\sqrt{\rho_{m}}\|_{L^{\infty}(\Omega)}
     \|{\mathbf u}_{m}\|_{{\bf L}^{\infty}(\Omega)}\|\sqrt{\rho_{m}}\partial_t\phi_{m}\|_{L^{2}(\Omega)}\|\nabla\phi_{m}\|_{{\bf L}^{2}(\Omega)}\\
	&\;\;\;\;\leqslant\|{\mathbf u}_{m}\|_{{\bf W}_{0}^{1,p}(\Omega)}
	\|\sqrt{\rho_{m}}\partial_t\phi_{m}\|_{L^{2}(\Omega)}\|\phi_{m}\|_{H^{1}(\Omega)}\\
	&\;\;\;\;\leqslant\frac{1}{16}\|\sqrt{\rho_{m}}\partial_t\phi_{m}\|^{2}_{L^{2}(\Omega)}
    +C\|{\mathbf u}_{m}\|^{p}_{{\bf W}_{0}^{1,p}(\Omega)}+C\|\phi_{m}\|^{\frac{2p}{p-2}}_{H^{1}(\Omega)},\\
	&|I_{3}|=|\int_{\Omega}3\rho_{m}\mu_{m}\phi_{m}^{2}{\mathbf u}_{m}\cdot\nabla\phi_{m} dx|\\
	&\;\;\;\;\leqslant C \|\rho_{m}\|_{ L^{\infty}(\Omega)}
	\|{\mathbf u}_{m}\|_{{\bf L}^{\infty}(\Omega)}\|\mu_{m}\|_{L^{3}(\Omega)}\|\phi_{m}\|^{2}_{L^{6}(\Omega)}
    \|\nabla\phi_{m}\|_{{\bf L}^{2}(\Omega)}\\
	&\;\;\;\;\leqslant C\|{\mathbf u}_{m}\|_{{\bf W}_{0}^{1,p}(\Omega)}\|\mu_{m}\|_{H^{1}(\Omega)}\|\phi_{m}\|^{3}_{H^{1}(\Omega)}\\
	&\;\;\;\;\leqslant C\|{\mathbf u}_{m}\|^{p}_{{\bf W}_{0}^{1,p}(\Omega)}+C\|\mu_{m}\|^{2}_{H^{1}(\Omega)}\|\phi_{m}\|^{6}_{H^{1}(\Omega)},\\
&|I_{5}|=|\int_{\Omega}\rho_{m}\phi_{m}^{3}{\mathbf u}_{m}\cdot\nabla\mu_{m}dx|\\
	&\;\;\;\;\leqslant \|\sqrt{\rho_{m}}\|_{L^{\infty}(\Omega)}\|\nabla\mu_{m}\|_{{\bf L}^{2}(\Omega)}
    \|\phi_{m}^{3}\|_{L^{2}(\Omega)}\|{\mathbf u}_{m}\|_{{\bf L}^{\infty}(\Omega)}\\
	&\;\;\;\;\leqslant C\|\mu_{m}\|_{H^{1}(\Omega)}\|\phi_{m}\|^{3}_{H^{1}(\Omega)}\|{\mathbf u}_{m}\|_{{\bf W}_{0}^{1,p}(\Omega)}\\
	&\;\;\;\;\leqslant C\|{\mathbf u}_{m}\|^{p}_{{\bf W}_{0}^{1,p}(\Omega)}+C\|\mu_{m}\|^{2}_{H^{1}(\Omega)}\|\phi_{m}\|^{6}_{H^{1}(\Omega)}.
\end{align*}
Similarly, one estimates $I_{1}-I_{3}$ and $I_{5}$ different for the case of $p=3$ as follows. For  for any $r\geqslant9,$
\begin{align*}
&|I_{1}|=|\int_{\Omega}-\rho_{m}\mu_{m}{\bf u}_{m}\cdot\nabla\mu_{m} dx|\\
	&\;\;\;\;\leqslant C\|{\mathbf u}_{m}\|_{{\bf L}^{r}(\Omega)}\|\sqrt{\rho_{m}}\mu_{m}\|^{\frac{r-3}{r}}_{L^{2}(\Omega)}
    \|\mu_{m}\|^{\frac{3}{p}}_{L^{6}(\Omega)}\|\nabla\mu_{m}\|_{{\bf L}^{2}(\Omega)}\\
	&\;\;\;\;\leqslant C\|{\mathbf u}_{m}\|_{{\bf W}_{0}^{1,p}(\Omega)}\|\sqrt{\rho_{m}}\mu_{m}\|^{\frac{r-3}{r}}_{L^{2}(\Omega)}
	\|\mu_{m}\|^{\frac{r+3}{r}}_{H^{1}(\Omega)}\\
	&\;\;\;\;\leqslant C\|{\mathbf u}_{m}\|_{{\bf W}_{0}^{1,3}(\Omega)}\|\sqrt{\rho_{m}}\mu_{m}\|^{\frac{r-3}{r}}_{L^{2}(\Omega)}
    \|\mu_{m}\|^{\frac{r-3}{r}}_{H^{1}(\Omega)}
	\|\mu_{m}\|^{\frac{6}{r}}_{H^{1}(\Omega)}\\
	&\;\;\;\;\leqslant C\|\sqrt{\rho_{m}}\mu_{m}\|^{2}_{L^{2}(\Omega)}\|\mu_{m}\|^{2}_{H^{1}(\Omega)}
    +C\left(\|{\mathbf u}_{m}\|^{p}_{{\bf W}_{0}^{1,p}(\Omega)}+\|\mu_{m}\|^{2}_{H^{1}(\Omega)}\right),\\
&|I_{2}|=|\int_{\Omega}-\rho_{m} {\mathbf u}_{m}\cdot\nabla\phi_{m}\partial_t\phi_{m}dx|\\
	&\;\;\;\;\leqslant C\|\sqrt{\rho_{m}}{\mathbf u}_{m}\|^{\frac{2r-6}{3r-6}}_{{\bf{L}}^{2}(\Omega)}
     \|{\mathbf u}_{m}\|^{\frac{r}{3r-6}}_{{\bf L}^{r}(\Omega)}\|\sqrt{\rho_{m}}\partial_t\phi_{m}\|_{L^{2}(\Omega)}
     \|\nabla\phi_{m}\|_{{\bf L}^{6}(\Omega)}\\
	&\;\;\;\;\leqslant C\|\sqrt{\rho_{m}}{\mathbf u}_{m}\|^{\frac{2r-6}{3r-6}}_{{\bf{L}}^{2}(\Omega)}
     \|{\mathbf u}_{m}\|^{\frac{r}{3r-6}}_{{\bf W}_{0}^{1,p}(\Omega)}\|\sqrt{\rho_{m}}\partial_t\phi_{m}\|_{L^{2}(\Omega)}
     \|\phi_{m}\|_{H^{2}(\Omega)}\\
	&\;\;\;\;\leqslant\frac{1}{16}\|\sqrt{\rho_{m}}\partial_t\phi_{m}\|^{2}_{L^{2}(\Omega)}
     +C\left(\|{\mathbf u}_{m}\|^{p}_{{\bf W}_{0}^{1,p}(\Omega)}
   +\|\phi_{m}\|^{4}_{H^{2}(\Omega)}+\|\sqrt{\rho_{m}}{\mathbf u}_{m}\|^{\frac{24r-72}{5r-18}}_{{\bf{L}}^{2}(\Omega)}\right),\\
&|I_{3}|=|\int_{\Omega}3\rho_{m}\mu_{m}\phi_{m}^{2}{\mathbf u}_{m}\cdot\nabla\phi_{m} dx|\\
		&\;\;\;\;\leqslant C\|\sqrt{\rho_{m}}{\mathbf u}_{m}\|^{\frac{2r-6}{3r-6}}_{{\bf{L}}^{2}(\Omega)}
      \|{\mathbf u}_{m}\|^{\frac{r}{3r-6}}_{{\bf L}^{r}(\Omega)}\|\nabla\phi_{m}\|_{{\bf L}^{6}(\Omega)}\|\phi_{m}^{2}\|_{L^{3}(\Omega)}
      \|\mu_{m}\|_{L^{6}(\Omega)}\\
	&\;\;\;\;\leqslant C\|\sqrt{\rho_{m}}{\mathbf u}_{m}\|^{\frac{2r-6}{3r-6}}_{{\bf{L}}^{2}(\Omega)}
      \|{\mathbf u}_{m}\|^{\frac{r}{3r-6}}_{{\bf W}_{0}^{1,p}(\Omega)}\|\mu_{m}\|_{H^{1}(\Omega)}\|\phi_{m}\|^{2}_{H^{1}(\Omega)}
       \|\phi_{m}\|_{H^{2}(\Omega)}\\
	&\;\;\;\;\leqslant C\left(\|{\mathbf u}_{m}\|^{p}_{{\bf W}_{0}^{1,p}(\Omega)}+\|\mu_{m}\|^{2}_{H^{1}(\Omega)}
    +\|\phi_{m}\|^{4}_{H^{2}(\Omega)}+\|\sqrt{\rho_{m}}{\mathbf u}_{m}\|^{\frac{24r-72}{5r-18}}_{{\bf{L}}^{2}(\Omega)}\right),\\
&|I_{5}|=|\int_{\Omega}\rho_{m}\phi_{m}^{3}{\mathbf u}_{m}\cdot\nabla\mu_{m}dx|\\
	&\;\;\;\;\leqslant C\|\nabla\mu_{m}\|_{{\bf L}^{2}(\Omega)}\|\phi_{m}^{3}\|_{L^{\frac{2r}{r-2}}(\Omega)}
   \|{\mathbf u}_{m}\|_{{\bf L}^{r}(\Omega)}\\
	&\;\;\;\;\leqslant C\|\mu_{m}\|_{H^{1}(\Omega)}\|\phi_{m}\|^{\frac{6}{r}}_{H^{2}(\Omega)}\|{\mathbf u}_{m}\|_{{\bf W}_{0}^{1,p}(\Omega)}\\
	&\;\;\;\;\leqslant C\left(|{\mathbf u}_{m}\|^{p}_{{\bf W}_{0}^{1,p}(\Omega)}+\|\mu_{m}\|^{2}_{H^{1}(\Omega)}
    +\|\phi_{m}\|^{\frac{36}{r}}_{H^{2}(\Omega)}\right).
\end{align*}
Collecting the estimates for $I_{1}-I_{8}$, one can obtain that
\begin{align*}
	&\frac{d}{dt}\int_{\Omega}\rho_{m}\mu^{2}_{m}dx+\int_{\Omega}\rho(\partial_t\phi_{m})^{2}dx\\
	&\leqslant C\|\sqrt{\rho_{m}}\mu_{m}\|^{2}_{L^{2}(\Omega)}\left(\|\mu_{m}\|^{2}_{H^{1}(\Omega)}+\|\phi_{m}\|^{2}_{H^{2}(\Omega)}\right)\\
	&+ C\left(\|\mu_{m}\|^{2}_{H^{1}(\Omega)}\|\phi_{m}\|^{6}_{H^{1}(\Omega)}+\|\mu_{m}\|^{2}_{H^{1}(\Omega)}+\|\phi_{m}\|^{4}_{H^{2}(\Omega)}+\|{\mathbf u}_{m}\|^{p}_{{\bf W}_{0}^{1,p}(\Omega)}+1\right).
\end{align*}
So, one deduces from Gronwall's inequality and \eqref{E-3-8}-\eqref{E-3-13} that
\begin{equation}\label{E3-19}
	\|\sqrt{\rho_{m}}\mu_{m}\|_{L^{\infty}(0,T;L^{2}(\Omega))}\leqslant C,\;\;\;\|\sqrt{\rho_{m}}\partial_t\phi_{m}\|_{L^{2}(0,T;L^{2}(\Omega))}\leqslant C.
\end{equation}
where $C$ is independent of $m$ and $\delta.$ Moreover, it is deduced from $\eqref{E3-14}$ that
\begin{equation}\label{E3-20}
	\|\Delta\phi_{m}\|_{L^{\infty}(0,T;L^{2}(\Omega))}\leqslant C.
\end{equation}
Furthermore, one gets from Lemma \ref{lem2-3} that
\begin{equation*}
	\|\phi_{m}\|_{L^{\infty}(0,T;H^{2}(\Omega))}\leqslant C,
\end{equation*}
where $C$ is independent of $m$ and $\delta$.
	\end{proof}

\begin{lemma}\label{lemmma3-8}
	Under the condition of Theorem \ref{thm2-1}, there exists a positive constant $C$, independent of $m$ and $\delta$, such that
	 \begin{equation}\label{E3-21}
\|\mu_{m}\|_{L^{2}(0,T;H^{2}(\Omega))}\leqslant C,\;\;\;\|\phi_{m}\|_{ L^{2}(0,T;W^{2,\infty}(\Omega))}\leqslant C.
	\end{equation}
\end{lemma}

\begin{proof}
 Taking $w=\Delta\mu_{m}$ in $\eqref{E3-3}_{3},$ one finds that
\begin{align*}
	&\int_{\Omega}|\Delta\mu_{m}|^{2}dx=\int_{\Omega}\left(\rho_{m}\Delta\mu_{m}{\mathbf u}_{m}\cdot\nabla\phi_{m}+\rho_{m}\partial_{t}\phi_{m}\Delta\mu_{m}\right) dx\\
	&\leqslant\|\sqrt{\rho_{m}}\|_{L^{\infty}(\Omega)}\|\sqrt{\rho_{m}}\partial_t\phi_{m}\|_{L^{2}(\Omega)}\|\Delta\mu_{m}\|_{L^{2}(\Omega)}+\|\rho_{m}\|_{L^{\infty}(\Omega)}\|{\mathbf u}_{m}\|_{{\bf L}^{3}(\Omega)}\|\nabla\phi_{m}\|_{{\bf L}^{6}(\Omega)}\|\Delta\mu_{m}\|_{L^{2}(\Omega)}\\
	&\leqslant\frac{1}{2}\|\Delta\mu_{m}\|^{2}_{L^{2}(\Omega)}+C\left(\|\sqrt{\rho_{m}}\partial_t\phi_{m}\|^{2}_{L^{2}(\Omega)}+\|{\mathbf u}_{m}\|^{p}_{{\mathbf{W}}_{0}^{1,p}(\Omega)}\|\phi_{m}\|^{2}_{H^{2}(\Omega)}\right).
\end{align*}
So, it follows Lemma \ref{lem2-3} and \eqref{E-3-8}-\eqref{E3-15} that there exists a positive constant $C$, independent of $m$ and $\delta,$ such that
\begin{equation}\label{E3-22}
	\|\mu_{m}\|_{L^{2}(0,T;H^{2}(\Omega))}\leqslant C.
\end{equation}
The equation $\rho_{m}\mu_{m}=-\Delta\phi_{m}+\rho_{m}\Psi^\prime(\phi_{m})$, gives that
 \begin{align*}
 	\|\Delta\phi_{m}\|_{L^{\infty}(\Omega)}\leqslant\|\rho_{m}\mu_{m}\|_{L^{\infty}(\Omega)}+\|\rho_{m}\Psi^\prime(\phi_{m})\|_{L^{\infty}(\Omega)}\leqslant C\left(\|\mu_{m}\|_{H^{2}(\Omega)}+\|\phi_{m}\|_{H^{2}(\Omega)}\right).
 \end{align*}
Using Lemma \ref{lem2-3} again, one deduces from  \eqref{E3-15} and \eqref{E3-22} that
 \begin{equation*}
	\|\phi_{m}\|_{ L^{2}(0,T;W^{2,\infty}(\Omega))}\leqslant C,
\end{equation*}
where $C$ is independent of $m$ and $\delta$.
\end{proof}

\begin{lemma}\label{lemma3-9}
	Under the condition of Theorem \ref{thm2-1}, there exists a positive constant $C$, independent of $m$ and $\delta$, such that
	\begin{equation}\label{E3-23}
	\|\mu_{m}\|_{L^{\infty}(0,T;H^{1}(\Omega))}\leqslant C,\;\;\;\|\partial_t\phi_{m}\|_{L^{2}(0,T;H^{1}(\Omega))}\leqslant C.
	\end{equation}
\end{lemma}

\begin{proof}
First, one obtains from $\eqref{E3-3}_4$ that
\begin{align}
   &\int_{\Omega}|\nabla\partial_t\phi_{m}|^{2}+3\rho_{m}\phi_{m}^{2}(\partial_t\phi_{m})^{2}dx\notag\\
   &=\int_{\Omega}\rho_{m} {\mathbf u}_{m}\cdot\nabla\mu_{m}\partial_t\phi_{m}dx
   +\int_{\Omega}\rho_{m}\mu_{m} {\mathbf u}_{m}\cdot\nabla\partial_t\phi_{m}dx+\int_{\Omega}\rho_{m}\partial_t\mu_m\partial_t\phi_{m}dx\notag\\
   &-3\int_{\Omega}\rho_{m}\phi^{2}_{m}{\mathbf u}_{m}\cdot\nabla\phi_{m}\partial_t\phi_{m}dx
    +\int_{\Omega}\rho_{m}\partial_t\phi_{m}{\mathbf u}_{m}\cdot\nabla\phi_{m} dx-\int_{\Omega}\rho_{m}\phi^{3}_{m}{\mathbf u}_{m}\cdot\nabla\partial_t\phi_{m}dx\notag\\
	&+\int_{\Omega}\rho_{m}\phi_{m} {\mathbf u}_{m}\cdot\nabla\partial_t\phi_{m}dx+\int_{\Omega}\rho_{m}(\partial_t\phi_{m})^{2}dx\label{(3a)}.
\end{align}
Second, one selects $w=\partial_t\mu_m$ in $\eqref{E3-3}_4$ and $w=\partial_t\phi_{m}$ in $\eqref{E3-3}_4$ to get that
\begin{equation}\label{(3b)}
	\frac{d}{dt}\int_{\Omega}\frac{1}{2}|\nabla\mu_{m}|^{2}dx
 =-\int_{\Omega}\rho_{m}\partial_t\mu_m\partial_t\phi_{m}dx-\int_{\Omega}\rho_{m}{\mathbf u}_{m}\cdot\nabla\phi_{m}\partial_t\mu_mdx
\end{equation}
and
\begin{equation}\label{(3c)}
	\int_{\Omega}\rho_{m}(\partial_t\phi_{m})^{2}dx+\int_{\Omega}\rho_{m}{\mathbf u}_{m}\cdot\nabla\phi_{m}\partial_t\phi_{m}dx=\int_{\Omega}\nabla\mu_{m}\cdot\nabla\partial_t\phi_{m}dx.
\end{equation}
Third, one takes the time derivative of $\eqref{E3-3}_4$ and testes the resultant by $w={\mathbf u}_{m}\cdot\nabla\phi_{m}$ to obtain that
\begin{align}\label{(3d)}
  &\int_{\Omega}\rho_{m}\partial_t\mu_m{\mathbf u}_{m}\cdot\nabla\phi_{m} dx\notag\\
  &=-\int_{\Omega}\rho_{m} {\mathbf u}_{m}\cdot\nabla\mu_{m}({\mathbf u}_{m}\cdot\nabla\phi_{m})dx
  -\int_{\Omega}\rho_{m}\mu_{m} {\mathbf u}_{m}\cdot(\nabla {\mathbf u}_{m}\cdot\nabla\phi_{m})dx
  +\int_{\Omega}\nabla\partial_t\phi_{m}\cdot(\nabla {\mathbf u}_{m}\cdot\nabla\phi_{m})dx\notag\\
	&-\int_{\Omega}\rho_{m}\mu_{m} {\mathbf u}_{m}\cdot({\mathbf u}_{m}\cdot\nabla^{2}\phi_{m})dx
    +\int_{\Omega}{\mathbf u}_{m}\cdot(\nabla^{2} \phi_{m}\cdot\nabla\partial_t\phi_{m})dx
    +3\int_{\Omega}\rho_{m}\phi^{2}_{m}{\mathbf u}_{m}\cdot\nabla\phi_{m}({\mathbf u}_{m}\cdot\nabla\phi_{m})dx\notag\\
	&+\int_{\Omega}\rho_{m}\phi^{3}_{m}{\mathbf u}_{m}\cdot\nabla({\mathbf u}_{m}\cdot\nabla\phi_{m})dx
   -\int_{\Omega}\rho_{m} {\mathbf u}_{m}\cdot\nabla\phi_{m}({\mathbf u}_{m}\cdot\nabla\phi_{m})dx
   -\int_{\Omega}\rho_{m}\phi_{m} {\mathbf u}_{m}\cdot\nabla({\mathbf u}_{m}\cdot\nabla\phi_{m})dx\notag\\
	&+3\int_{\Omega}\rho_{m}\phi^{2}_{m}\partial_t\phi_{m}({\mathbf u}_{m}\cdot\nabla\phi_{m})dx
	-\int_{\Omega}\rho_{m}\partial_t\phi_{m}({\mathbf u}_{m}\cdot\nabla\phi_{m})dx.
\end{align}
Now, one sums \eqref{(3a)}-\eqref{(3d)} to arrive at
\begin{align}\label{(E3-35)}
&\frac{d}{dt}\int_{\Omega}\left(\frac{1}{2}|\nabla\mu_{m}|^{2}dx+\int_{\Omega}|\nabla\partial_t\phi_{m}|^{2}dx+\int_{\Omega}3\rho_{m}\phi_{m}^{2}(\partial_t\phi_{m})^{2}\right)dx\notag\\
	&+3\int_{\Omega}\rho_{m}\phi_{m}^{2}|{\mathbf u}_{m}|^{2}|\nabla\phi_{m}|^{2}dx+\int_{\Omega}\rho_{m}(\partial_t\phi_{m})^{2}dx\notag\\
	&=\underset{I_{1}}{\underbrace{\int_{\Omega}\rho_{m} {\mathbf u}_{m}\cdot\nabla\mu_{m}\partial_t\phi_{m}dx}}+\underset{I_{2}}{\underbrace{\int_{\Omega}\rho_{m}\mu_{m}{\mathbf u}_{m}\cdot\nabla\partial_t\phi_{m}dx}}+\underset{I_{3}}{\underbrace{\int_{\Omega}\rho_{m}\phi_{m} {\mathbf u}_{m}\cdot\nabla\partial_t\phi_{m}dx}}\notag\\
	&-\underset{I_{4}}{\underbrace{\int_{\Omega}\rho_{m}\phi^{3}_{m}{\mathbf u}_{m}\cdot\nabla\partial_t\phi_{m}dx}}-\underset{I_{5}}{\underbrace{6\int_{\Omega}\rho_{m}\phi^{2}_{m}{\mathbf u}_{m}\cdot\nabla\phi_{m}\partial_t\phi_{m}dx}}-\underset{I_{6}}{\underbrace{\int_{\Omega}\nabla\partial_t\phi_{m}\cdot(\nabla {\mathbf u}_{m}\cdot\nabla\phi_{m})dx}}\notag\\
	&-\underset{I_{7}}{\underbrace{\int_{\Omega}{\mathbf u}_{m}\cdot(\nabla^{2} \phi_{m}\cdot\nabla\partial_t\phi_{m})dx}}+\underset{I_{8}}{\underbrace{\int_{\Omega}\rho_{m} {\mathbf u}_{m}\cdot\nabla\mu_{m}({\mathbf u}_{m}\cdot\nabla\phi_{m})dx}}+\underset{I_{9}}{\underbrace{\int_{\Omega}\rho_{m}\mu_{m} {\mathbf u}_{m}\cdot(\nabla {\mathbf u}_{m}\cdot\nabla\phi_{m})dx}}\notag\\
	&+\underset{I_{10}}{\underbrace{\int_{\Omega}\rho_{m}\mu_{m} {\mathbf u}_{m}\cdot({\mathbf u}_{m}\cdot\nabla^{2}\phi_{m})dx}}-\underset{I_{11}}{\underbrace{\int_{\Omega}\rho_{m}\phi^{3}_{m}{\mathbf u}_{m}\cdot({\mathbf u}_{m}\cdot\nabla^{2}\phi_{m})dx}}-\underset{I_{12}}{\underbrace{\int_{\Omega}\rho_{m}\phi^{3}_{m}{\mathbf u}_{m}\cdot(\nabla {\mathbf u}_{m}\cdot\nabla\phi_{m})dx}}\notag\\
	&+\underset{I_{13}}{\underbrace{\int_{\Omega}\rho_{m}|{\mathbf u}_{m}|^{2}|\nabla\phi_{m}|^{2}dx}}+\underset{I_{14}}{\underbrace{\int_{\Omega}\rho_{m}\phi_{m} {\mathbf u}_{m}\cdot(\nabla {\mathbf u}_{m}\cdot\nabla\phi_{m})dx}}+\underset{I_{15}}{\underbrace{\int_{\Omega}\rho_{m}\phi_{m} {\mathbf u}_{m}\cdot( {\mathbf u}_{m}\cdot\nabla^{2}\phi_{m})dx}}.
\end{align}
For the case of $\frac{5}{2}<p<3$, the estimation of $I_{1}-I_{15}$ are stated as follows
\begin{align*}
	&|I_{1}|=|\int_{\Omega}\rho_{m} {\mathbf u}_{m}\cdot\nabla\mu_{m}\partial_t\phi_{m}dx|\\
	&\;\;\;\;\leqslant\int_{\Omega}\rho_{m} |{\mathbf u}_{m}||\nabla\mu_{m}|^{\frac{3-p}{p}}|\nabla\mu_{m}|^{\frac{2p-3}{p}}|\partial_t\phi_{m}|dx\\
	&\;\;\;\;\leqslant C\|\sqrt{\rho_{m}}\partial_t\phi_{m}\|_{L^{2}(\Omega)}\|\nabla\mu_{m}\|^{\frac{3-p}{p}}_{{\bf L}^{6}(\Omega))}\|\nabla\mu_{m}\|^{\frac{2p-3}{p}}_{{\bf L}^{2}(\Omega)}\|{\mathbf u}_{m}\|_{{\bf L}^{\frac{3p}{3-p}}(\Omega))}\\
	&\;\;\;\;\leqslant C\|\sqrt{\rho_{m}}\partial_t\phi_{m}\|_{L^{2}(\Omega)}\|\mu_{m}\|^{\frac{3-p}{p}}_{H^{2}(\Omega)}\|\nabla\mu_{m}\|^{\frac{3p-3}{p}}_{{\bf L}^{2}(\Omega)}\|{\mathbf u}_{m}\|_{{\bf W}_{0}^{1,p}(\Omega)}\\
	&\;\;\;\;\leqslant C\|\sqrt{\rho_{m}}\partial_t\phi_{m}\|^{2}_{L^{2}(\Omega)}+C\left(\|\mu_{m}\|^{2}_{H^{2}(\Omega)}+\|\nabla\mu_{m}\|^{2}_{{\bf L}^{2}(\Omega)}\|{\mathbf u}_{m}\|^{p}_{{\bf W}_{0}^{1,p}(\Omega)}\right),\\
	&|I_{2}|=|\int_{\Omega}\rho_{m}\mu_{m}{\mathbf u}_{m}\cdot\nabla\partial_t\phi_{m}dx|\\
	&\;\;\;\;\leqslant\|\nabla\partial_t\phi_{m}\|_{{\bf L}^{2}(\Omega)}\|{\mathbf u}_{m}\|_{{\bf L}^{\frac{3p}{3-p}}(\Omega))}\|\mu_{m}\|_{L^{\frac{6p}{5p-6}}(\Omega)}\|\rho_{m}\|_{L^{\infty}(\Omega)}\\
	&\;\;\;\;\leqslant C\|\nabla\partial_t\phi_{m}\|_{{\bf L}^{2}(\Omega)}\|{\mathbf u}_{m}\|_{{\bf W}_{0}^{1,p}(\Omega))}\|\nabla\mu_{m}\|_{{\bf L}^{2}(\Omega)}\\
	&\;\;\;\;\leqslant\epsilon\|\nabla\partial_t\phi_{m}\|^{2}_{{\bf L}^{2}(\Omega)}+C\|{\mathbf u}_{m}\|^{p}_{{\bf W}_{0}^{1,p}(\Omega)}\|\nabla\mu_{m}\|^{2}_{{\bf L}^{2}(\Omega)},\\
	&|I_{3}|=|\int_{\Omega}\rho_{m}\phi_{m} {\mathbf u}_{m}\cdot\nabla\partial_t\phi_{m}dx|\\
	&\;\;\;\;\leqslant \|\rho_{m}\|_{L^{\infty}(\Omega)}\|\phi_{m}\|_{L^{\infty}(\Omega)}\|\nabla\partial_t\phi_{m}\|_{{\bf L}^{2}(\Omega)}\|{\mathbf u}_{m}\|_{{\bf L}^{2}(\Omega)}\\
	&\;\;\;\;\leqslant C\|\phi_{m}\|_{H^{2}(\Omega)}\|\nabla\partial_t\phi_{m}\|_{{\bf L}^{2}(\Omega)}\|{\mathbf u}_{m}\|_{{\bf L}^{2}(\Omega)}\\
	&\;\;\;\;\leqslant\epsilon\|\nabla\partial_t\phi_{m}\|^{2}_{{\bf L}^{2}(\Omega)}+C\|{\mathbf u}_{m}\|^{p}_{{\bf W}_{0}^{1,p}(\Omega)}\|\phi_{m}\|^{2}_{H^{2}(\Omega)},\\
	&|I_{4}|=|\int_{\Omega}\rho_{m}\phi^{3}_{m}{\mathbf u}_{m}\cdot\nabla\partial_t\phi_{m}dx|\\
	&\;\;\;\;\leqslant\|\phi^{3}_{m}\|_{L^{\infty}(\Omega)}\|\rho_{m}\|_{L^{\infty}(\Omega)}\|{\mathbf u}_{m}\|_{{\bf L}^{2}(\Omega)}\|\nabla\partial_t\phi_{m}\|_{{\bf L}^{2}(\Omega)}\\
	&\;\;\;\;\leqslant C\|\phi_{m}\|^{3}_{H^{2}(\Omega)}\|{\mathbf u}_{m}\|_{{\bf W}_{0}^{1,p}(\Omega)}\|\nabla\partial_t\phi_{m}\|_{{\bf L}^{2}(\Omega)}\\
	&\;\;\;\;\leqslant\epsilon\|\nabla\partial_t\phi_{m}\|^{2}_{{\bf L}^{2}(\Omega)}+C\|\phi_{m}\|^{6}_{H^{2}(\Omega)}\|{\mathbf u}_{m}\|^{p}_{{\bf W}_{0}^{1,p}(\Omega)},\\
	&|I_{5}|=|6\int_{\Omega}\rho_{m}\phi^{2}_{m}{\mathbf u}_{m}\cdot\nabla\phi_{m}\partial_t\phi_{m}dx|\\
	&\;\;\;\;\leqslant\|\sqrt{\rho_{m}}\phi_{m} |{\mathbf u}_{m}||\nabla\phi_{m}|\|_{{\bf L}^{2}(\Omega)}\|\phi_{m}\|_{L^{\infty}(\Omega)}\|\sqrt{\rho_{m}}\partial_t\phi_{m}\|_{L^{2}(\Omega)}\\
	&\;\;\;\;\leqslant C\|\sqrt{\rho_{m}}\partial_t\phi_{m}\|^{2}_{L^{2}(\Omega)}+\frac{1}{2}\|\sqrt{\rho_{m}}\phi_{m}|{\mathbf u}_{m}||\nabla\phi_{m}|\|^{2}_{{\bf L}^{2}(\Omega)},\\
	&|I_{6}|=|\int_{\Omega}\nabla\partial_t\phi_{m}\cdot(\nabla {\mathbf u}_{m}\cdot\nabla\phi_{m})dx|\\
	&\;\;\;\;\leqslant\|\nabla\partial_t\phi_{m}\|_{{\bf L}^{2}(\Omega)}\|\nabla {\mathbf u}_{m}\|_{{\bf L}^{p}(\Omega)}\|\nabla\phi_{m}\|_{{\bf L}^{\frac{2p}{p-2}}(\Omega)}\\
	&\;\;\;\;\leqslant\epsilon\|\nabla\partial_t\phi_{m}\|^{2}_{{\bf L}^{2}(\Omega)}+C\left(\| {\mathbf u}_{m}\|^{p}_{{\bf W}_{0}^{1,p}(\Omega)}+\|\phi_{m}\|^{2}_{W^{2,\infty}(\Omega)}\right),\\
	&|I_{7}|=|\int_{\Omega}{\mathbf u}_{m}\cdot(\nabla^{2} \phi_{m}\cdot\nabla\partial_t\phi_{m})dx|\\
	&\;\;\;\;\leqslant\|\nabla\partial_t\phi_{m}\|_{{\bf L}^{2}(\Omega)}\| {\mathbf u}_{m}\|_{{\bf L}^{\frac{3p}{3-p}}(\Omega)}\|\nabla^{2}\phi_{m}\|_{{\bf L}^{\frac{6p}{5p-6}}(\Omega)}\\	&\;\;\;\;\leqslant\epsilon\|\nabla\partial_t\phi_{m}\|^{2}_{{\bf L}^{2}(\Omega)}+C\left(\| {\mathbf u}_{m}\|^{p}_{{\bf W}_{0}^{1,p}(\Omega)}+\|\phi_{m}\|^{2}_{W^{2,\infty}(\Omega)}\right),\\
	&|I_{8}|=|\int_{\Omega}\rho_{m} {\mathbf u}_{m}\cdot\nabla\mu_{m}({\mathbf u}_{m}\cdot\nabla\phi_{m})dx|\\
&\;\;\;\;\leqslant \|\sqrt{\rho_{m}}\|_{{\bf L}^{\infty}(\Omega)}\|\sqrt{\rho_{m}}{\mathbf u}_{m}\|_{{\bf L}^{2}(\Omega)}\|{\mathbf u}_{m}\|_{{\bf L}^{\frac{3p}{3-p}}(\Omega)}\|\nabla\mu_{m}\|_{{\bf L}^{6}(\Omega)}\|\nabla\phi_{m}\|_{{\bf L}^{\frac{3p}{2p-3}}(\Omega)}\\
&\;\;\;\;\leqslant C\|{\mathbf u}_{m}\|_{{\bf W}_{0}^{1,p}(\Omega)}\|\mu_{m}\|_{H^{2}(\Omega)}\|\phi_{m}\|_{H^{2}(\Omega)}\\
	&\;\;\;\;\leqslant C\|{\mathbf u}_{m}\|^{p}_{{\bf W}_{0}^{1,p}(\Omega)}+C\|\mu_{m}\|^{2}_{H^{2}(\Omega)}\|\phi_{m}\|^{2}_{H^{2}(\Omega)},\\
	&|I_{9}|=|\int_{\Omega}\rho_{m}\mu_{m} {\mathbf u}_{m}\cdot(\nabla {\mathbf u}_{m}\cdot\nabla\phi_{m})dx|\\
	&\;\;\;\;\leqslant \|\sqrt{\rho_{m}}\|_{{\bf L}^{\infty}(\Omega)}\|\nabla {\mathbf u}_{m}\|_{{\bf L}^{p}(\Omega)}\|\sqrt{\rho_{m}}{\mathbf u}_{m}\|_{{\bf L}^{2}(\Omega)}\|\mu_{m}\|_{L^{\infty}(\Omega)}\|\nabla\phi_{m}\|_{{\bf L}^{\frac{2p}{p-2}}(\Omega)}\\
	&\;\;\;\;\leqslant C|{\mathbf u}_{m}\|_{{\bf W}_{0}^{1,p}(\Omega)}\|\mu_{m}\|_{H^{2}(\Omega)}\|\nabla\phi_{m}\|_{{\bf L}^{\frac{2p}{p-2}}(\Omega)}\\
	&\;\;\;\;\leqslant C\|{\mathbf u}_{m}\|^{p}_{{\bf W}_{0}^{1,p}(\Omega)}+C\left(\|\mu_{m}\|^{2}_{H^{2}(\Omega)}+\|\phi_{m}\|^{2}_{W^{2,\infty}(\Omega)}\right),\\
	&|I_{10}|=|\int_{\Omega}\rho_{m}\mu_{m} {\mathbf u}_{m}\cdot({\mathbf u}_{m}\cdot\nabla^{2}\phi_{m})dx|\\
	&\;\;\;\;\leqslant \|\sqrt{\rho_{m}}\|_{{\bf L}^{\infty}(\Omega)}\|\sqrt{\rho_{m}}{\mathbf u}_{m}\|_{{\bf L}^{2}(\Omega)}\| {\mathbf u}_{m}\|_{{\bf L}^{\frac{3p}{3-p}}(\Omega)}\|\mu_{m}\|_{L^{\infty}(\Omega)}\|\nabla^{2}\phi_{m}\|_{{\bf L}^{\frac{6p}{5p-6}}(\Omega)}\\
	&\;\;\;\;\leqslant C\| {\mathbf u}_{m}\|_{{\bf W}_{0}^{1,p}(\Omega)}\|\mu_{m}\|_{H^{2}(\Omega)}\|\nabla^{2}\phi_{m}\|_{{\bf L}^{\frac{6p}{5p-6}}(\Omega)}\\
&\;\;\;\;\leqslant C\|{\mathbf u}_{m}\|^{p}_{{\bf W}_{0}^{1,p}(\Omega)}+C\left(\|\mu_{m}\|^{2}_{H^{2}(\Omega)}+\|\phi_{m}\|^{2}_{W^{2,\infty}(\Omega)}\right),\\
	&|I_{11}|=|\int_{\Omega}\rho_{m}\phi^{3}_{m}{\mathbf u}_{m}\cdot({\mathbf u}_{m}\cdot\nabla^{2}\phi_{m})dx|\\
	&\;\;\;\;\leqslant \|\rho_{m}\|_{L^{\infty}(\Omega)}\|\phi_{m}^{3}\|_{ L^{\infty}(\Omega)}\|{\mathbf u}_{m}\|^{2}_{{\bf L}^{4}(\Omega)}\|\nabla^{2}\phi_{m}\|_{{\bf L}^{2}(\Omega)}\\
	&\;\;\;\;\leqslant C\|\phi_{m}\|^{3}_{ H^{2}(\Omega)}\|{\mathbf u}_{m}\|^{2}_{{\bf W}_{0}^{1,p}(\Omega)}\|\phi_{m}\|_{H^{2}(\Omega)}\\
	&\;\;\;\;\leqslant C\|{\mathbf u}_{m}\|^{p}_{{\bf W}_{0}^{1,p}(\Omega)}\|\phi_m\|^{4}_{H^{2}(\Omega)},\\
	&|I_{12}|=|\int_{\Omega}\rho_{m}\phi^{3}_{m}{\mathbf u}_{m}\cdot(\nabla {\mathbf u}_{m}\cdot\nabla\phi_{m})dx|\\
	&\;\;\;\;\leqslant \|\rho_{m}\|_{L^{\infty}(\Omega)}\|\phi_{m}^{3}\|_{ L^{\infty}(\Omega)}\|{\mathbf u}_{m}\|_{{\bf L}^{3}(\Omega)}\|\nabla{\mathbf u}_{m}\|_{{\bf L}^{2}(\Omega)}\|\nabla\phi_{m}\|_{{\bf L}^{6}(\Omega)}\\
	&\;\;\;\;\leqslant C\|\phi_{m}\|^{3}_{ H^{2}(\Omega)}\|{\mathbf u}_{m}\|^{2}_{{\bf W}_{0}^{1,p}(\Omega)}\|\phi_{m}\|_{H^{2}(\Omega)}\\
	&\;\;\;\;\leqslant C\|{\mathbf u}_{m}\|^{p}_{{\bf W}_{0}^{1,p}(\Omega)}\|\phi_m\|^{4}_{H^{2}(\Omega)},\\
	&|I_{13}|=|\int_{\Omega}\rho_{m}{\mathbf u}_{m}^{2}\nabla\phi_{m}^{2}dx|\\
	&\;\;\;\;\leqslant \|\rho_{m}\|_{L^{\infty}(\Omega)}\|{\mathbf u}_{m}\|^{2}_{{\bf L}^{4}(\Omega)}\|\nabla\phi_{m}\|^{2}_{{\bf L}^{4}(\Omega)}\\
	&\;\;\;\;\leqslant C\|{\mathbf u}_{m}\|^{2}_{{\bf W}_{0}^{1,p}(\Omega)}\|\phi_{m}\|^{2}_{H^{2}(\Omega)}\\
	&\;\;\;\;\leqslant C\|{\mathbf u}_{m}\|^{p}_{{\bf W}_{0}^{1,p}(\Omega)}\|\phi_m\|^{4}_{H^{2}(\Omega)},\\
	&|I_{14}|=|\int_{\Omega}\rho_{m}\phi_{m} {\mathbf u}_{m}\cdot(\nabla {\mathbf u}_{m}\cdot\nabla\phi_{m})dx|\\
	&\;\;\;\;\leqslant \|\rho_{m}\|_{L^{\infty}(\Omega)}\|\phi_m\|_{L^{\infty}(\Omega)}\|{\mathbf u}_{m}\|_{{\bf L}^{3}(\Omega)}\|\nabla{\mathbf u}_{m}\|_{{\bf L}^{2}(\Omega)}\|\nabla\phi_{m}\|_{{\bf L}^{6}(\Omega)}\\
	&\;\;\;\;\leqslant C\|\phi_m\|_{H^{2}(\Omega)}\|{\mathbf u}_{m}\|^{2}_{{\bf W}_{0}^{1,p}(\Omega)}\|\phi_{m}\|_{H^{2}(\Omega)}\\
	&\;\;\;\;\leqslant C\|{\mathbf u}_{m}\|^{p}_{{\bf W}_{0}^{1,p}(\Omega)}\|\phi_m\|^{2}_{H^{2}(\Omega)},\\
	&|I_{15}|=|\int_{\Omega}\rho_{m}\phi_{m} {\mathbf u}_{m}\cdot( {\mathbf u}_{m}\cdot\nabla^{2}\phi_{m})dx|\\
	&\;\;\;\;\leqslant \|\rho_{m}\|_{L^{\infty}(\Omega)}\|\phi_{m}\|_{L^{\infty}(\Omega)}\|{\mathbf u}_{m}\|^{2}_{L^{4}(\Omega)}\|\nabla^{2}\phi_{m}\|_{{\bf L}^{2}(\Omega)}\\
	&\;\;\;\;\leqslant C\|{\mathbf u}_{m}\|^{2}_{{\bf W}_{0}^{1,p}(\Omega)}\|\phi_{m}\|^{2}_{H^{2}(\Omega)}\\
	&\;\;\;\;\leqslant C\|{\mathbf u}_{m}\|^{p}_{{\bf W}_{0}^{1,p}(\Omega)}\|\phi_m\|^{2}_{H^{2}(\Omega)},
\end{align*}
where $\epsilon\in(0,1)$ is any fixed, and $C$ is independent of $m$ and $\delta.$ For the case of $p> 3$, one just to estimate $I_{j}\;(j=1,2,7,8,10)$  different in the following
\begin{align*}
	&|I_{1}|=|\int_{\Omega}\rho_{m} {\mathbf u}_{m}\cdot\nabla\mu_{m}\partial_t\phi_{m}dx|\\
	&\;\;\;\;\leqslant \|\sqrt{\rho_{m}}\|_{L^{\infty}(\Omega)}\|\sqrt{\rho_{m}}\partial_t\phi_{m}\|_{L^{2}(\Omega)}\|\nabla\mu_{m}\|_{{\bf L}^{2}(\Omega))}\|{\mathbf u}_{m}\|_{{\bf L}^{\infty}(\Omega))}\\
&\;\;\;\;\leqslant C\|\sqrt{\rho_{m}}\partial_t\phi_{m}\|_{L^{2}(\Omega)}\|\mu_{m}\|_{H^{1}(\Omega))}\|{\mathbf u}_{m}\|_{{\bf W}_{0}^{1,p}(\Omega)}\\
	&\;\;\;\;\leqslant\epsilon\|\sqrt{\rho_{m}}\partial_t\phi_{m}\|^{2}_{L^{2}(\Omega)}+C\|\mu_{m}\|^{2}_{H^{1}(\Omega)}\|{\mathbf u}_{m}\|^{p}_{{\bf W}_{0}^{1,p}(\Omega)},\\
	&|I_{2}|=|\int_{\Omega}\rho_{m}\mu_{m}{\mathbf u}_{m}\cdot\nabla\partial_t\phi_{m}dx|\\
	&\;\;\;\;\leqslant\|\rho_{m}\|_{L^{\infty}(\Omega)}\|\nabla\partial_t\phi_{m}\|_{{\bf L}^{2}(\Omega)}\|{\mathbf u}_{m}\|_{{\bf L}^{\infty}(\Omega)}\|\mu_{m}\|_{L^{2}(\Omega)}\\
	&\;\;\;\;\leqslant C\|\nabla\partial_t\phi_{m}\|_{{\bf L}^{2}(\Omega)}\|{\mathbf u}_{m}\|_{{\bf W}_{0}^{1,p}(\Omega)}\|\mu_{m}\|_{L^{2}(\Omega)}\\
	&\;\;\;\;\leqslant\epsilon\|\nabla\partial_t\phi_{m}\|^{2}_{{\bf L}^{2}(\Omega)}+C\|{\mathbf u}_{m}\|^{p}_{{\bf W}_{0}^{1,p}(\Omega)}\|\mu_{m}\|^{2}_{L^{2}(\Omega)},\\
	&|I_{7}|=|\int_{\Omega}{\mathbf u}_{m}\cdot(\nabla^{2} \phi_{m}\cdot\nabla\partial_t\phi_{m})dx|\\
	&\;\;\;\;\leqslant\|\nabla\partial_t\phi_{m}\|_{{\bf L}^{2}(\Omega)}\| {\mathbf u}_{m}\|_{{\bf L}^{\infty}(\Omega)}\|\nabla^{2}\phi_{m}\|_{{\bf L}^{2}(\Omega)}\\
	&\;\;\;\;\leqslant\|\nabla\partial_t\phi_{m}\|_{{\bf L}^{2}(\Omega)}\| {\mathbf u}_{m}\|_{{\bf W}_{0}^{1,p}(\Omega)}\|\phi_{m}\|_{H^{2}(\Omega)}\\
	&\;\;\;\;\leqslant\epsilon\|\nabla\partial_t\phi_{m}\|^{2}_{{\bf L}^{2}(\Omega)}+C\| {\mathbf u}_{m}\|^{p}_{{\bf W}_{0}^{1,p}(\Omega)}\|\phi_{m}\|^{2}_{H^{2}(\Omega)},\\
	&|I_{8}|=|\int_{\Omega}\rho_{m} {\mathbf u}_{m}\cdot\nabla\mu_{m}({\mathbf u}_{m}\cdot\nabla\phi_{m})dx|\\
	&\;\;\;\;\leqslant\|\sqrt{\rho_{m}}\|_{{\bf L}^{\infty}(\Omega)} \|\sqrt{\rho_{m}}{\mathbf u}_{m}\|_{{\bf L}^{2}(\Omega)}\|{\mathbf u}_{m}\|_{{\bf L}^{\infty}(\Omega)}\|\nabla\mu_{m}\|_{{\bf L}^{3}(\Omega)}\|\nabla\phi_{m}\|_{{\bf L}^{6}(\Omega)}\\
	&\;\;\;\;\leqslant C\|\sqrt{\rho_{m}}{\mathbf u}_{m}\|_{{\bf L}^{2}(\Omega)}\|{\mathbf u}_{m}\|_{{\bf W}_{0}^{1,p}(\Omega)}\|\mu_{m}\|_{H^{2}(\Omega)}\|\phi_{m}\|_{H^{2}(\Omega)}\\
	&\;\;\;\;\leqslant C\|{\mathbf u}_{m}\|^{p}_{{\bf W}_{0}^{1,p}(\Omega)}+C\|\mu_{m}\|^{2}_{H^{2}(\Omega)}\|\phi_{m}\|^{2}_{H^{2}(\Omega)},\\
	&|I_{10}|=|\int_{\Omega}\rho_{m}\mu_{m} {\mathbf u}_{m}\cdot({\mathbf u}_{m}\cdot\nabla^{2}\phi_{m})dx|\\
	&\;\;\;\;\leqslant\|\sqrt{\rho_{m}}\|_{{\bf L}^{\infty}(\Omega)} \|\sqrt{\rho_{m}}{\mathbf u}_{m}\|_{{\bf L}^{2}(\Omega)}\|{\mathbf u}_{m}\|_{{\bf L}^{\infty}(\Omega)}\|\mu_{m}\|_{L^{\infty}(\Omega)}\|\nabla^{2}\phi_{m}\|_{{\bf L}^{2}(\Omega)}\\
	&\;\;\;\;\leqslant C\|\sqrt{\rho_{m}}{\mathbf u}_{m}\|_{{\bf L}^{2}(\Omega)}\|{\mathbf u}_{m}\|_{{\bf W}_{0}^{1,p}(\Omega)}\|\mu_{m}\|_{H^{2}(\Omega)}\|\phi_{m}\|_{H^{2}(\Omega)}\\
	&\;\;\;\;\leqslant C\|{\mathbf u}_{m}\|^{p}_{{\bf W}_{0}^{1,p}(\Omega)}+C\|\phi_{m}\|^{2}_{H^{2}(\Omega)}\|\mu_{m}\|^{2}_{H^{2}(\Omega)},
\end{align*}
where $\epsilon\in(0,1)$ is any fixed, and $C$ is independent of $m$ and $\delta.$
Taking similar argument, one estimates $I_{j}\;(j=1,2,7,8,10)$ for the case $p= 3$ as follows. For any $r\geqslant 9,$
\begin{align*}
	&|I_{1}|=|\int_{\Omega}\rho_{m} {\mathbf u}_{m}\cdot\nabla\mu_{m}\partial_t\phi_{m}dx|\\
	&\;\;\;\;\leqslant C\|\sqrt{\rho_{m}}\partial_t\phi_{m}\|_{L^{2}(\Omega)}\|\nabla\mu_{m}\|_{{\bf L}^{\frac{2r}{r-2}}(\Omega))}\|{\mathbf u}_{m}\|_{{\bf L}^{r}(\Omega))}\\
	&\;\;\;\;\leqslant C\|\sqrt{\rho_{m}}\partial_t\phi_{m}\|_{L^{2}(\Omega)}\|\mu_{m}\|^{\frac{3}{r}}_{H^{2}(\Omega))}\|{\mathbf u}_{m}\|_{{\bf W}_{0}^{1,p}(\Omega)}\\
	&\;\;\;\;\leqslant\epsilon\|\sqrt{\rho_{m}}\partial_t\phi_{m}\|^{2}_{L^{2}(\Omega)}
+C\left(\|\mu_{m}\|^{\frac{18}{r}}_{H^{2}(\Omega)}+\|{\mathbf u}_{m}\|^{p}_{{\bf W}_{0}^{1,p}(\Omega)}\right),\\
	&|I_{2}|=|\int_{\Omega}\rho_{m}\mu_{m}{\mathbf u}_{m}\cdot\nabla\partial_t\phi_{m}dx|\\
	&\;\;\;\;\leqslant C\|\nabla\partial_t\phi_{m}\|_{{\bf L}^{2}(\Omega)}\|{\mathbf u}_{m}\|_{{\bf L}^{r}(\Omega)}\|\mu_{m}\|_{L^{\frac{2r}{r-2}}(\Omega)}\\
	&\;\;\;\;\leqslant C\|\nabla\partial_t\phi_{m}\|_{{\bf L}^{2}(\Omega)}\|{\mathbf u}_{m}\|_{{\bf W}_{0}^{1,p}(\Omega)}\|\mu_{m}\|_{H^{1}(\Omega)}\\
	&\;\;\;\;\leqslant\epsilon\|\nabla\partial_t\phi_{m}\|^{2}_{{\bf L}^{2}(\Omega)}+C\|{\mathbf u}_{m}\|^{p}_{{\bf W}_{0}^{1,p}(\Omega)}\|\mu_{m}\|^{2}_{H^{1}(\Omega)},\\
	&|I_{7}|=|\int_{\Omega}{\mathbf u}_{m}\cdot(\nabla^{2} \phi_{m}\cdot\nabla\partial_t\phi_{m})dx|\\
	&\;\;\;\;\leqslant\|\nabla\partial_t\phi_{m}\|_{{\bf L}^{2}(\Omega)}\| {\mathbf u}_{m}\|_{{\bf L}^{r}(\Omega)}\|\nabla^{2}\phi_{m}\|_{{\bf L}^{\frac{2r}{r-2}}(\Omega)}\\
	&\;\;\;\;\leqslant\|\nabla\partial_t\phi_{m}\|_{{\bf L}^{2}(\Omega)}\| {\mathbf u}_{m}\|_{{\bf W}_{0}^{1,p}(\Omega)}\|\phi_{m}\|^{\frac{2}{r}}_{W^{2,\infty}(\Omega)}\\
	&\;\;\;\;\leqslant\epsilon\|\nabla\partial_t\phi_{m}\|^{2}_{{\bf L}^{2}(\Omega)}+C\left(\| {\mathbf u}_{m}\|^{p}_{{\bf W}_{0}^{1,p}(\Omega)}+\|\phi_{m}\|^{\frac{12}{r}}_{W^{2,\infty}(\Omega)}\right),\\
	&|I_{8}|=|\int_{\Omega}\rho_{m} {\mathbf u}_{m}\cdot\nabla\mu_{m}({\mathbf u}_{m}\cdot\nabla\phi_{m})dx|\\
	&\;\;\;\;\leqslant C \|\sqrt{\rho_{m}}{\mathbf u}_{m}\|_{{\bf L}^{2}(\Omega)}\|{\mathbf u}_{m}\|_{{\bf L}^{R}(\Omega)}\|\nabla\mu_{m}\|_{{\bf L}^{\frac{3r}{r-3}}(\Omega)}\|\nabla\phi_{m}\|_{{\bf L}^{6}(\Omega)}\\
	&\;\;\;\;\leqslant C\|\sqrt{\rho_{m}}{\mathbf u}_{m}\|_{{\bf L}^{2}(\Omega)}\|{\mathbf u}_{m}\|_{{\bf W}_{0}^{1,p}(\Omega)}\|\mu_{m}\|_{H^{2}(\Omega)}\|\phi_{m}\|_{H^{2}(\Omega)}\\
	&\;\;\;\;\leqslant C\|{\mathbf u}_{m}\|^{p}_{{\bf W}_{0}^{1,p}(\Omega)}+C\|\mu_{m}\|^{2}_{H^{2}(\Omega)}\|\phi_{m}\|^{2}_{H^{2}(\Omega)},\\
	&|I_{10}|=|\int_{\Omega}\rho_{m}\mu_{m} {\mathbf u}_{m}\cdot({\mathbf u}_{m}\cdot\nabla^{2}\phi_{m})dx|\\
	&\;\;\;\;\leqslant C \|\sqrt{\rho_{m}}{\mathbf u}_{m}\|_{{\bf L}^{2}(\Omega)}\|{\mathbf u}_{m}\|_{{\bf L}^{r}(\Omega)}\|\mu_{m}\|_{L^{\infty}(\Omega)}\|\nabla^{2}\phi_{m}\|_{{\bf L}^{\frac{2r}{r-2}}(\Omega)}\\
	&\;\;\;\;\leqslant C\|{\mathbf u}_{m}\|_{{\bf W}_{0}^{1,p}(\Omega)}\|\mu_{m}\|_{H^{2}(\Omega)}\|\phi_{m}\|^{\frac{2}{r}}_{W^{2,\infty}(\Omega)}\\
	&\;\;\;\;\leqslant C\|{\mathbf u}_{m}\|^{p}_{{\bf W}_{0}^{1,p}(\Omega)}+C\left(\|\phi_{m}\|^{\frac{12}{r}}_{W^{2,\infty}(\Omega)}+\|\mu_{m}\|^{2}_{H^{2}(\Omega)}\right).
\end{align*}
Collecting $I_{1}-I_{15}$  together, one can deduce from \eqref{E-3-8}-\eqref{E3-21} that
\begin{align*}
&\frac{d}{dt}\int_{\Omega}|\nabla\mu_{m}|^{2}dx
+\int_{\Omega}\left(|\nabla\partial_t\phi_{m}|^{2}+\rho_{m}\phi_{m}^{2}(\partial_t\phi_{m})^{2}+\rho_{m}\phi_{m}^{2}|{\mathbf u}_{m}|^{2}|\nabla\phi_{m}|^{2}+\rho_{m}(\partial_t\phi_{m})^{2}\right)dx\notag\\
	&\leqslant C_{1}\|\nabla\mu_{m}\|^{2}_{{\bf L}^{2}(\Omega)}\|{\mathbf u}_{m}\|^{p}_{{\bf W}_{0}^{1,p}(\Omega)}\\
	&+C_{2} \left(\|{\mathbf u}_{m}\|^{2}_{{\bf W}_{0}^{1,p}(\Omega)}\|\phi_m\|^{2}_{H^{2}(\Omega)}+\|{\mathbf u}_{m}\|^{p}_{{\bf W}_{0}^{1,p}(\Omega)}+\|\mu_{m}\|^{2}_{H^{2}(\Omega)}+\|\sqrt{\rho_{m}}\partial_t\phi_{m}\|^{2}_{L^{2}(\Omega)}+\|\nabla^{2}\phi_{m}\|^{2}_{{\bf W}^{2,\infty}(\Omega)}\right)
\end{align*}
where positive constants $C_{1}$ and $C_{2}$ are independent of $m$ and $\delta$. According to Lemma \ref*{gronwall1} and \eqref{E-3-8}-\eqref{E3-21}, there exists a positive constant $C,$ independent $m$ and $\delta$, such that
\begin{equation*}
	\|\mu_{m}\|_{L^{\infty}(0,T;H^{1}(\Omega))}\leqslant C.
\end{equation*}
Further,
\begin{equation}\label{E3-29}
	\|\nabla\partial_t\phi_{m}\|_{{\bf L}^{2}(0,T;{\bf L}^{2}(\Omega))}\leqslant C.
\end{equation}
Combining \eqref{E3-15} and \eqref{E3-29}, one finds that
\begin{equation*}
	\|\partial_t\phi_{m}\|_{L^{2}(0,T;H^{1}(\Omega))}\leqslant C.
\end{equation*}
\end{proof}

\subsection{Passage to the limit}

\quad Now we are ready to pass the limit for $m\rightarrow\infty$ in the sequence $(\rho_m,{\bf u}_m,\phi_m,\mu_m).$

\begin{proposition}\label{prop3-2} Let $\frac{5}{2}<p<+\infty$ and fix $\delta>0.$ Suppose that \eqref{E1-3}-\eqref{E1-4} are satisfied, $T\in(0,+\infty),$ and the initial data $(\rho_0,{\mathbf u}_0,\phi_0,\mu_0)$ are given as in Theorem \ref{thm2-1}. The following problem
	\begin{eqnarray}\label{EA-1}
		\begin{cases}
			\partial_{t}\rho+div(\rho {\mathbf u})=0~~\quad\mbox{ in }\Omega\times (0,T),\\
			\rho\partial_t {\mathbf u}+\rho {\mathbf u}\cdot \nabla {\mathbf u}-div(\nu(\phi)(1+|\mathbb{D}{\mathbf u}|^2)^{\frac{p-2}{2}}\mathbb{D}{\mathbf u})+\nabla P=-div(\nabla\phi\otimes\nabla\phi)~~\quad\mbox{ in }\Omega\times (0,T),\\
			div{\mathbf u}=0~~\quad\mbox{ in }\Omega\times (0,T),\\
			\rho \partial_t \phi+\rho {\mathbf u}\cdot\nabla\phi=\Delta\mu~~\quad\mbox{ in }\Omega\times (0,T),\\
			\rho\mu=-\Delta\phi+\rho\Psi^\prime(\phi)~~\quad\mbox{ in }\Omega\times (0,T),\\
			{\mathbf u}=0 \qquad\qquad\qquad\qquad\qquad\qquad\;\mbox{ on }\;\partial\Omega\times(0,T),\\
			\partial_{\bf n}\phi=\partial_{\bf n}\mu=0\qquad\qquad\qquad\qquad\;\;\;\mbox{ on }\;\partial\Omega\times(0,T),\\
			\rho_{t=0}=\rho_{0\delta},~~{\mathbf u}_{t=0}={\mathbf u}_0,~~\phi|_{t=0}=\phi_0 \;\;\quad\mbox{ in }\;\Omega,
		\end{cases}
	\end{eqnarray}
	admits a weak solution $(\rho^\delta,{\bf u}^\delta,\phi^\delta,\mu^\delta)$ on $[0,T]$ such that
	\begin{align*}
		&\rho^\delta\in C([0,T];L^{r}(\Omega))\cap L^{\infty}((0,T)\times\Omega)\;(\forall r\in[1,+\infty)),\;\partial_{t}\rho^\delta\in L^{\infty}(0,T;H^{-1}(\Omega)),\\
		&\sqrt{\rho^\delta}{\mathbf u}^\delta\in{\bf L}^{\infty}(0,T;{\bf L}^{2}(\Omega)),\;
		{\mathbf u}^\delta\in{\bf L}^{p}([0,T],{\bf W}_{0}^{1,p}(\Omega)),\\
		&\sqrt{\rho^\delta}\phi^\delta\in L^{\infty}(0,T;L^{2}(\Omega)),\;\phi^\delta\in L^{\infty}(0,T;H^{2}(\Omega))\cap C([0,T];H^{1}(\Omega)),\\
		&\sqrt{\rho^\delta}\mu^\delta\in L^{\infty}(0,T;L^{2}(\Omega)),\;\mu^\delta\in L^{\infty}(0,T;H^{1}(\Omega))\cap L^{2}(0,T;H^{2}(\Omega))
	\end{align*}
	and the following system
\begin{eqnarray}\label{E-3-35}
\begin{cases}
	(\partial_{t}\rho^\delta,\eta)+(\rho^\delta{\bf u}^\delta,\nabla\eta)=0~~~~a.e.\ in\ (0,T),\\
	(\rho^\delta\partial_{t}{\mathbf u}^\delta,{\bf \zeta})+(\rho^\delta\left({\bf u}^\delta\cdot\nabla\right){\mathbf u}^\delta,{\mathbf \zeta})
			+(\nu(\psi^\delta)|1+|{\mathbb D}{{\mathbf u}}^\delta|^{2}|^{\frac{p-2}{2}}\mathbb{D}{\mathbf u}^\delta,{\mathbb D} {\mathbf \zeta})\\
	~~~~~~~~~~~~~~=(\rho^\delta\mu^\delta\nabla\phi^\delta,{\mathbf \zeta})-(\rho^\delta\nabla\Psi(\phi^\delta),{\mathbf \zeta})~~~~a.e.\ in\ (0,T),\\
			(\rho^\delta\partial_{t}\phi^\delta,\xi)+(\rho^\delta {\mathbf u}^\delta\cdot\nabla\phi^\delta,\xi)
		+(\nabla\mu^\delta,\nabla \xi)=0~~~~a.e.\ in\ (0,T),\\
			\rho^\delta\mu^\delta=-\Delta\phi^\delta+\rho^\delta\Psi^\prime(\phi^\delta)~~~~a.e.\ \mbox{ in} \ \Omega\times (0,T),\\
		div{\mathbf u}^\delta=0~~~~a.e.\ \mbox{ in }\ \Omega\times (0,T),\\
			{\mathbf u}^\delta=0,\partial_{n}\mu^\delta=\partial_{n}\phi^\delta=0 ~~~\mbox{ on }\partial\Omega\times (0,T),\\
		\rho^\delta(\cdot,0)=\rho_{0\delta},~{\mathbf u}^\delta(\cdot,0)={\mathbf u}_{0},~\phi^\delta(\cdot,0)=\phi_{0}^\delta ~~~\mbox{ in }\Omega
		\end{cases}
	\end{eqnarray}
holds for all test functions $\eta\in H^1(\Omega), {\bf \zeta}\in {\bf C}_{0,\sigma}^\infty(\Omega)$ and $\xi\in H^1(\Omega)$.
\end{proposition}

\begin{proof}
The uniform estimates \eqref{E-3-8}--\eqref{E3-23} allow us to apply the weak compactness theorems together with the J. Simon \cite{Simon-1987} to extract a suitable subsequence that converges to a limit in the corresponding topologies as $m \to +\infty$. Except for $\eqref{E-3-35}_2$, the remaining convergence arguments follow standard procedures (see, e.g., \cite{Giorgini-2020, Feireisl-2004}) and are therefore omitted for brevity. In what follows, we only verify $\eqref{E-3-35}_2$.So, there exists a subsequence of $(\rho_m, {\bf u}_m, \phi_m, \mu_m)$ and $(\rho^\delta, {\bf u}^\delta, \phi^\delta, \mu^\delta)$ such that	
	\begin{eqnarray}\label{E3-32}
		\begin{cases}
			\rho_{m}\rightarrow\rho^\delta\; weak-star\;in\;L^{\infty}(\Omega\times(0, T)),\\
			{\mathbf u}_{m}\rightarrow{\mathbf u}^\delta\;weakly\;in\;{\bf L}^{p}(0,T;{\bf W}_{0}^{1,p}(\Omega)),\\
			\phi_{m}\rightarrow\phi^\delta\;weak-star\;in\;L^{\infty}(0,T;H^{1}(\Omega)),\\
			\phi_{m}\rightarrow\phi^\delta\;weakly\;in\;L^{4}(0,T;H^{2}(\Omega)),\\
			\mu_{m}\rightarrow\mu^{\delta}\;weakly\;in\;L^{2}(0,T;H^{1}(\Omega)).
		\end{cases}
	\end{eqnarray}
The convergence of the density is obtained similarly from [\cite{PL-1996}, Lemma 2.4], that is,
	\begin{equation}\label{E3-33}
		\rho_{m}\rightarrow\rho^\delta\;strongly\; in\;C([0,T];L^{ r }(\Omega)),(\forall  r \in[1,\infty)).
	\end{equation}
In light of \eqref{E3-32} and \eqref{E3-33}, one deduces that
\begin{equation}
\rho_{m}\mu_{m}\rightarrow\rho^{\delta}\mu^{\delta}\;weakly\;in\;L^{2}(0,T;L^{\frac{6r}{r+6}}(\Omega))
\end{equation}
for any fixed $r\in[1,\infty).$ By using compact embedding, one finds that
\begin{equation}
	\mathbf{u}_{m}\rightarrow\mathbf{u}^{\delta}\;in\;{\bf L}_{w}^{p}(0,T;{\bf L}_{s}^{\frac{3p}{3-p}}(\Omega))
\end{equation}
for $\frac{5}{2}<p<3$ and
\begin{equation}\label{(convergence3-35)}
	\mathbf{u}_{m}\rightarrow\mathbf{u}^{\delta}\;in\;{\bf L}_{w}^{p}(0,T;{\bf L}_{s}^{q}(\Omega))\;(\forall 1<q<\infty)
\end{equation}
for $p\geqslant3$.  Now, it is the position to deduce the convergence of $\phi_{m}$ and $\mu_{m}.$ Combine with \eqref{E3-15} and \eqref{E3-23}, one has that
\begin{equation*}
	\phi_{m}\rightarrow\phi^\delta\;in\;C_s([0,T];H_{w}^{2}(\Omega)).
\end{equation*}
Further,
\begin{align}\label{(convergence4)}
	&\nabla\phi_{m}\rightarrow\nabla\phi^\delta\;\;strongly\; in\; {\bf C}([0,T];{\bf L}^{r}(\Omega))\;\;(\forall\;r\in(1,6)),\\
	&\phi_{m}\rightarrow\phi^\delta\;\;strongly\; in\; C([0,T];L^{r}(\Omega))\;\;(\forall\;r\in(1,\infty)).
\end{align}
So, it follows from \eqref{E3-33} that
\begin{equation}\label{(convergence4-1)}
\begin{cases}
	\rho_{m}\phi_{m}\rightarrow\rho^\delta\phi^\delta\;\;strongly\; in\;C([0,T];L^{2}(\Omega)),\\
\rho_{m}\Psi^\prime(\phi_{m})\rightarrow\rho^\delta\Psi^\prime(\phi^\delta)\;\;strongly\; in\;C([0,T];L^{2}(\Omega)),\\
\rho_{m}\Psi^\prime(\phi_{m})\nabla\phi_{m}\rightarrow\rho^\delta\Psi^\prime(\phi^\delta)\nabla\phi^{\delta}\;strongly\; in\;{\bf C}([0,T];{\bf L}^{2}(\Omega)).
\end{cases}
\end{equation}
Above all, one obtains that
\begin{eqnarray}\label{convergence-2}
	\begin{cases}
	\rho_{m}{\mathbf u}_{m}\rightarrow\rho^\delta {\mathbf u}^\delta\;weakly\;in\;{\bf L}^{p}(0,T;{\bf L}^{\frac{6p}{p+6}}(\Omega)),\\
	\rho_{m}{\mathbf u}_{m}\otimes {\mathbf u}_{m}\rightarrow \rho^\delta {\mathbf u}^\delta\otimes {\mathbf u}^\delta\;weakly\;in\;{\bf L}^{p}(0,T;{\bf L}^{\frac{2p}{p+2}}(\Omega)),\\
	\rho_{m}\phi_{m}\rightarrow\rho^\delta\phi^\delta\;strongly\;in\;C([0,T];L^{2}(\Omega)),\\
	\rho_{m}\Psi^\prime(\phi_{m})\rightarrow\rho^\delta\Psi^\prime(\phi^\delta)\;strongly\;in\;C([0,T];L^{2}(\Omega)),\\
	\rho_{m}\Psi^\prime(\phi_{m})\nabla\phi_{m}\rightarrow\rho^\delta\Psi^\prime(\phi^\delta)\nabla\phi^\delta\;strongly\;in\;{\bf C}([0,T];{\bf L}^{2}(\Omega)),\\
	\rho_{m}{\mathbf u}_{m}\phi_{m}\rightarrow \rho^\delta  {\mathbf u}^\delta\phi^\delta\;weakly\;in\;{\bf L}^{2}(0,T;{\bf L}^{\frac{3}{2}}(\Omega)),\\
	\rho_{m}\mu_{m}\rightarrow\rho^{\delta}\mu^{\delta}\;weakly\;in\;L^{2}(0,T;L^{4}(\Omega)),\\
	\rho_{m}\mu_{m}\nabla\phi_{m}\rightarrow\rho^{\delta}\mu^{\delta}\nabla\phi^{\delta}\;weakly\;in\;{\bf L}^{2}(0,T;{\bf L}^{\frac{3}{2}}(\Omega)).
	\end{cases}
	\end{eqnarray}

Now, we focus on the limit of term $\nu(\phi_{m})(1+|{\mathbb D}{{\mathbf u}}_{m}|^{2})^{\frac{p-2}{2}}\mathbb{D}{\mathbf u}_{m}.$ For the sake of discussion, we denote
$(1+|{\mathbb D}{{\mathbf u}}_{m}|^{2})^{\frac{p-2}{2}}\mathbb{D}{\mathbf u}_{m}$ and $(1+|{\mathbb D}{{\mathbf u}}^\delta|^{2})^{\frac{p-2}{2}}\mathbb{D}{\mathbf u}^\delta$ by
${\cal T}({\bf u}_m)$ and ${\cal T}({\bf u}^\delta),$ respectively. From  \eqref{E3-32} and \eqref{(convergence4)}, it is easy to get that
\begin{equation*}
	\nu(\phi_{m}){\cal T}({\mathbf u}_{m})\rightarrow\nu(\phi)\overline{{\cal T}({\mathbf u}^\delta)}\;weakly\;in\;{\bf L}^{\frac{p}{p-1}}(0,T;{\bf L}^{\frac{p}{p-1}}(\Omega)),
\end{equation*}
where $\overline{{\cal T}({\mathbf u}^\delta)}$ denotes the weak limit of ${\cal T}({\mathbf u}_{m}).$

First, one deduces from $\eqref{E3-3}_{3},$ Lemma \ref{lemma3-1}, Lemma \ref{lemmma3-7} and Lemma \ref{lemma3-9} that
\begin{equation*}
	\|\partial_{t}(\rho_{m}\mu_{m})\|_{L^{2}(0,T;H^{-1}(\Omega))}\leqslant C,
\end{equation*}
where $C$ is independent of $m$ and $\delta.$ Combining $$\|\rho_{m}\mu_{m}\|_{ L^{\infty}(0,T;L^{4}(\Omega))}\leqslant C$$ with $C$ being independent of $m$ and $\delta,$ one gets that
\begin{equation*}
	\rho_{m}\mu_{m}\rightarrow\rho^{\delta}\mu^{\delta}in\; C_s([0,T]; L_{w}^{4}(\Omega)).
\end{equation*}
Using \eqref{(convergence4)}, one obtains that
\begin{equation}\label{(convergence4-2)}
	\rho_{m}\mu_{m}\nabla\phi_{m}\rightarrow\rho^{\delta}\mu^{\delta}\nabla\phi^{\delta}\;\; in\;{\bf C}_s([0,T];{\bf L}_{w}^{2}(\Omega)).
\end{equation}
So, one can deduce from \eqref{(convergence3-35)} and \eqref{(convergence4-2)} that
\begin{align}
&|\int_{0}^{t}\int_{\Omega}\rho_{m}\mu_{m}\nabla\phi_m\cdot{\mathbf u}_{m}dxd\tau-\int_{0}^{t}\int_{\Omega}\rho^{\delta}\mu^{\delta}\nabla\phi^{\delta}\cdot{\mathbf u}^{\delta}dxd\tau|\notag\\
&\leqslant|\int_{0}^{t}\int_{\Omega}\rho_{m}\mu_{m}\nabla\phi_m\cdot({\mathbf u}_{m}-{\mathbf u}^{\delta})dxd\tau|
+|\int_{0}^{t}\int_{\Omega}(\rho_{m}\mu_{m}\nabla\phi_m-\rho^{\delta}\mu^{\delta}\nabla\phi^{\delta})\cdot{\mathbf u}^{\delta}dxd\tau|\notag\\
&\leqslant\int_{0}^{t}\|\rho_{m}\mu_{m}\nabla\phi_m\|_{{\bf L}^{2}(\Omega)}\|{\mathbf u}_{m}-{\mathbf u}^{\delta}\|_{{\bf L}^{2}(\Omega)}d\tau+|\int_{0}^{t}\int_{\Omega}(\rho_{m}\mu_{m}\nabla\phi_m-\rho^{\delta}\mu^{\delta}\nabla\phi^{\delta})\cdot{\mathbf u}^{\delta}dxd\tau|\notag\\
&\leqslant\int_{0}^{t}(\|\rho_{m}\mu_{m}\nabla\phi_m\|_{{\bf L}^{2}(\Omega)}
-\mathop{\underline{\lim}}\limits_{m\rightarrow\infty}\|\rho_{m}\mu_{m}\nabla\phi_m\|_{{\bf L}^{2}(\Omega)})\|{\mathbf u}_{m}-{\mathbf u}^{\delta}\|_{{\bf L}^{2}(\Omega)}d\tau\notag\\
&+|\int_{0}^{t}\int_{\Omega}(\rho_{m}\mu_{m}\nabla\phi_m-\rho^{\delta}\mu^{\delta}\nabla\phi^{\delta})\cdot{\mathbf u}^{\delta}dxd\tau|
+\int_{0}^{t}\mathop{\underline{\lim}}\limits_{m\rightarrow\infty}\|\rho_{m}\mu_{m}\nabla\phi_m\|_{{\bf L}^{2}(\Omega)}\|{\mathbf u}_{m}-{\mathbf u}^{\delta}\|_{{\bf L}^{2}(\Omega)}d\tau\notag\\
&\leqslant C\int_{0}^{t}(\|\rho_{m}\mu_{m}\nabla\phi_m\|_{{\bf L}^{2}(\Omega)}
-\mathop{\underline{\lim}}\limits_{m\rightarrow\infty}\|\rho_{m}\mu_{m}\nabla\phi_m\|_{{\bf L}^{2}(\Omega)})d\tau\notag\\
&+|\int_{0}^{t}\int_{\Omega}(\rho_{m}\mu_{m}\nabla\phi_m-\rho^{\delta}\mu^{\delta}\nabla\phi^{\delta})\cdot{\mathbf u}^{\delta}dxd\tau|
+\int_{0}^{t}\mathop{\underline{\lim}}\limits_{m\rightarrow\infty}\|\rho_{m}\mu_{m}\nabla\phi_m\|_{{\bf L}^{2}(\Omega)}\|{\mathbf u}_{m}-{\mathbf u}^{\delta}\|_{{\bf L}^{2}(\Omega)}d\tau\notag\\
&\rightarrow 0
\end{align}
as $m\rightarrow\infty.$

Second, one takes ${\bf w}={\mathbf u}_{m}$ and ${\bf w}=\varphi$ in $\eqref{E3-3}_{2}$ respectively, to arrive at
\begin{align}\label{(1d)}
	&[\frac{1}{2}\int_{\Omega}\rho_{m}|{\mathbf u}_{m}|^{2}dx]_{0}^{\tau}-	[\frac{1}{2}\int_{\Omega}\rho_{m}{\mathbf u}_{m}\cdot\varphi dx]_{0}^{\tau}\notag\\
	&+\int_{0}^{\tau}\int_{\Omega}(\rho_{m} {\mathbf u}_{m}\cdot\partial_{t}\varphi+\rho_{m}{\mathbf u}_{m}\otimes {\mathbf u}_{m}:\nabla\varphi+\rho_{m} ({\mathbf u}_{m}-\varphi)\cdot\nabla\Psi(\phi_{m})\notag\\
	&-\rho_{m}\mu_{m}\nabla\phi_m\cdot ({\mathbf u}_{m}-\varphi)
	+\nu(\phi_{m})(1+|{\mathbb D}{\mathbf u}_{m}|^{2})^{\frac{p-2}{2}}{\mathbb D}{\mathbf u}_{m}:{\mathbb D}({\mathbf u}_{m}-\varphi))dxd\tau=0
\end{align}
holds for all $\varphi\in {\bf C}_{c}^{\infty}(\Omega\times[0,T], {\mathbb R}^{3}).$  Letting $m\rightarrow\infty$ in \eqref{(1d)}, one finds that
\begin{align}\label{(1e)}
	&[\frac{1}{2}\int_{\Omega}\rho^{\delta}|{\mathbf u}^{\delta}|^{2}dx]_{0}^{\tau}
  -[\frac{1}{2}\int_{\Omega}\rho^{\delta}{\mathbf u}^{\delta}\cdot\varphi dx]_{0}^{\tau}\notag\\
	&+\int_{0}^{\tau}\int_{\Omega}(\rho^{\delta} {\mathbf u}^{\delta}\cdot\partial_{t}\varphi+\rho^{\delta}{\mathbf u}^{\delta}\otimes {\mathbf u}^{\delta}:\nabla\varphi+\rho^{\delta} ({\mathbf u}^{\delta}-\varphi)\cdot\nabla\Psi(\phi^{\delta})
	-\rho^{\delta}\mu^{\delta}\nabla\phi^{\delta}\cdot ({\mathbf u}^{\delta}-\varphi)\notag\\
	&+\nu(\phi^{\delta})(\overline{(1+|{\mathbb D}{\mathbf u}^{\delta}|^{2})^{\frac{p-2}{2}}{\mathbb D}{\mathbf u}^{\delta}:{\mathbb D}{\mathbf u}^{\delta}}-\overline{(1+|{\mathbb D}{\mathbf u}^{\delta}|^{2})^{\frac{p-2}{2}}{\mathbb D}{\mathbf u}^{\delta}}:{\mathbb D}\varphi))dxd\tau\leqslant0
\end{align}
holds for all $\varphi\in {\bf C}_{c}^{\infty}(\Omega\times[0,T],{\mathbb R}^{3}).$ Consider a family of regularizing kernels
\begin{equation*}
	\eta_{h}(t):=\frac{1}{h}1_{[-h,0]}(t)\;\; \mathrm{and}\;\; \eta_{-h}(t):=\frac{1}{h}1_{[0,h]}(t)\; (h>0)
\end{equation*}
along with the cut-off functions
\begin{equation*}
	\zeta_{\vartheta}\in C_{c}^{\infty}(0,\tau), 0\leqslant\zeta\leqslant1, \zeta_{\vartheta}(t)=1\;\; \mathrm{whenever}\;\; t\in[\vartheta,\tau-\vartheta] (\vartheta>0).
\end{equation*}
Notice that $\eta_{h}*{\mathbf u}^{\delta}=\frac{1}{h}\int_{t}^{t+h}{\mathbf u}^{\delta}ds\in {\bf W}^{1,p}(0,T;{\bf W}_{0}^{1,p}(\Omega)).$ One can use the quantities
\begin{equation}\label{E3-40}
	\varphi_{h,\vartheta}=\zeta_{\vartheta}\eta_{-h}*\eta_{h}*(\zeta_{\vartheta}{\mathbf u}^{\delta})\; (\vartheta>0, h>0)
\end{equation}
as a test function in \eqref{(1e)}. Obviously, it leads to $\varphi_{h,\vartheta}\to{\mathbf u}^{\delta}$ in $\mathbf{L} ^{p}(0,T,\mathbf{W}_{0}^{1,p}(\Omega))$ as $\theta\to0^{+}$ and $h\to 0^{+}$. So,
\begin{align*}
&[\int_{\Omega}\rho^{\delta}{\mathbf u}^{\delta}\cdot\varphi_{h,\vartheta}dx]|_{0}^{\tau}=0 \;\;\mathrm{for}\;\mathrm{all}\;\; \vartheta,h>0, \\
&\lim\limits_{\vartheta\to0}\lim\limits_{h\to0}\int_{0}^{\tau}\int_{\Omega}\nu(\phi^{\delta})
[\overline{(1+|{\mathbb D}{\mathbf u}^{\delta}|^{2})^{\frac{p-2}{2}}{\mathbb D}{\mathbf u}^{\delta}:{\mathbb D}{\mathbf u}^{\delta}}
-\overline{(1+|{\mathbb D}{\mathbf u}^{\delta}|^{2})^{\frac{p-2}{2}}{\mathbb D}{\mathbf u}^{\delta}}:{\mathbb D}\varphi_{h,\vartheta})]dxdt
\geqslant0.
\end{align*}
Moreover,
\begin{eqnarray}
&&\int_{0}^{\tau}\int_{\Omega}\rho^{\delta}{\mathbf u}^{\delta}\cdot\partial_{t}\varphi_{h,\theta}dxdt\\
&&=\int_{0}^{\tau}\int_{\Omega}\rho^{\delta}{\mathbf u}^{\delta}\cdot\partial_{t}\zeta_{\theta}\eta_{-h}*\eta_{h}*(\zeta_{\theta}{\mathbf u}^{\delta})dxdt +\int_{{\mathbb R}^{1}}\int_{\Omega}[\eta_{h}*\rho^{\delta}\zeta_{\theta}{\mathbf u}^{\delta}]\cdot\partial_{t}[\eta_{h}*(\zeta_{\theta}{\mathbf u}^{\delta})]dxdt,\notag\\
&&\int_{0}^{\tau}\int_{\Omega}\rho^{\delta}{\mathbf u}^{\delta}\cdot\partial_{t}\zeta_{\vartheta}\eta_{-h}*\eta_{h}*(\zeta_{\vartheta}{\mathbf u}^{\delta})dxdt\notag\\
 &&=\lim_{\vartheta\to0}\int_{0}^{\tau}\left(\frac{1}{2}\int_{\Omega}\rho^{\delta}|{\mathbf u}^{\delta}|^{2}dx\right)\partial_{t}|\zeta_{\vartheta}|^{2}dt =\left.\left[\frac{1}{2}\int_{\Omega}\rho^{\delta}|{\mathbf u}^{\delta}|^{2}dx\right]\right|_{0}^{\tau}\notag
\end{eqnarray}
and
\begin{align*}
	&\int_{{\mathbb R}^{1}}\int_{\Omega}[\eta_{h}*\rho^{\delta}\zeta_{\vartheta}{\mathbf u}^{\delta}]\cdot\partial_{t}[\eta_{h}*(\zeta_{\vartheta}{\mathbf u}^{\delta})]dxdt\\
	& =\int_{{\mathbb R}^{1}}\int_{\Omega}[\eta_{h}*\rho^{\delta}\zeta_{\vartheta}{\mathbf u}^{\delta}]\cdot\partial_{t}[\eta_{h}*(\zeta_{\vartheta}{\mathbf u}^{\delta})]dxdt \\
	& =-\int_{{\mathbb R}^{1}}\int_{\Omega}\frac{\rho^{\delta}(t+h)-\rho^{\delta}(t)}{h}(\zeta_{\theta}{\mathbf u}^{\delta})(t+h)[\eta_{h}*(\zeta_{\vartheta}{\mathbf u}^{\delta})]dxdt\\
	&-\frac{1}{2}\int_{{\mathbb R}^{1}}\int_{\Omega}\rho^{\delta}\partial_{t}[\eta_{h}*(\zeta_{\vartheta}{\mathbf u}^{\delta})]^{2}dxdt.
\end{align*}
Thus,
\begin{align}\label{(I3)}
	& & & -\int_{{\mathbb R}^{1}}\int_{\Omega}\frac{\rho^{\delta}(t+h)-\rho^{\delta}(t)}{h}(\zeta_{\theta}{\mathbf u}^{\delta})(t+h)\cdot[\eta_{h}*(\zeta_{\theta}{\mathbf u}^{\delta})]dxdt \notag\\
	& & & =-\int_{{\mathbb R}^{1}} \int_{\Omega}\partial_{t}(\eta_{h}*\rho^{\delta})(\zeta_{\theta}{\mathbf u}^{\delta})(t+h)\cdot[\eta_{h}*(\zeta_{\theta}{\mathbf u}^{\delta})]dxdt\notag \\
	& & & =\int_{{\mathbb R}^{1}}\int_{\Omega}\mathrm{div}[\eta_{h}*(\rho^{\delta}{\mathbf u}^{\delta})](\zeta_{\theta}{\mathbf u}^{\delta})(t+h)\cdot[\eta_{h}*(\zeta_{\theta}{\mathbf u}^{\delta})]dxdt\notag \\
	& & & =\int_{{\mathbb R}^{1}} \int_{\Omega}\mathrm{div}[\eta_{h}*(\rho^{\delta}{\mathbf u}^{\delta})(\zeta_{\theta}{\mathbf u}^{\delta})(t+h)]\cdot[\eta_{h}*(\zeta_{\theta}{\mathbf u}^{\delta})]dxdt\notag \\
	& & &-\int_{{\mathbb R}^{1}}\int_{\Omega}[\eta_{h}*(\rho^{\delta}{\mathbf u}^{\delta})]\cdot\nabla(\zeta_{\theta}{\mathbf u}^{\delta})(t+h)\cdot[\eta_{h}*(\zeta_{\theta}{\mathbf u}^{\delta})]dxdt\notag \\
	& & & =-\int_{{\mathbb R}^{1}}\int_{\Omega}[\eta_{h}*(\rho^{\delta}{\mathbf u}^{\delta})]\cdot\nabla(\zeta_{\theta}{\mathbf u}^{\delta})(t+h)\cdot[\eta_{h}*(\zeta_{\theta}{\mathbf u}^{\delta})]dxdt
\end{align}	
and
\begin{align}\label{(I4)}
	& & & \;\;\;\;- \frac{1}{2}\int_{{\mathbb R}^{1}}\int_{\Omega}\rho^{\delta}\partial_{t}[\eta_{h}*(\zeta_{\theta}{\mathbf u}^{\delta})]^{2}dxdt \notag\\
	& & & =-\frac{1}{2}\int_{{\mathbb R}^{1}}\int_{\Omega}\partial_{t}(\rho^{\delta}(\eta_{h}*(\zeta_{\theta}{\mathbf u}^{\delta}))^{2})dxdt+\frac{1}{2}\int_{{\mathbb R}^{1}}\int_{\Omega}\partial_{t}\rho^{\delta}[\eta_{h}*(\zeta_{\theta}{\mathbf u}^{\delta})]^{2}dxdt \notag\\
	& & & =-\frac{1}{2}\int_{{\mathbb R}^{1}}\int_{\Omega}\mathrm{div}(\rho^{\delta}{\mathbf u}^{\delta})[\eta_{h}*(\zeta_{\theta}{\mathbf u}^{\delta})]^{2}dxdt \notag\\
	& & & =-\frac{1}{2}\int_{{\mathbb R}^{1}}\int_{\Omega}\mathrm{div}(\rho^{\delta}{\mathbf u}^{\delta}[\eta_{h}*(\zeta_{\theta}{\mathbf u}^{\delta})]^{2})dxdt+\frac{1}{2}\int_{{\mathbb R}^{1}}\int_{\Omega}\rho^{\delta}{\mathbf u}^{\delta}\cdot\nabla[\eta_{h}*(\zeta_{\theta}{\mathbf u}^{\delta})]^{2})dx\notag \\
	& & & =\frac{1}{2}\int_{{\mathbb R}^{1}}\int_{\Omega}\rho^{\delta}{\mathbf u}^{\delta}\cdot\nabla[\eta_{h}*(\zeta_{\theta}{\mathbf u}^{\delta})]^{2})dxdt.
\end{align}
Putting together \eqref{(I3)}-\eqref{(I4)}, one obtains that
\begin{align}
	& & & \int_{{\mathbb R}^{1}}\int_{\Omega}  [\eta_{h}*\rho^{\delta}\zeta_{\vartheta}{\mathbf u}^{\delta}]\cdot\partial_{t}[\eta_{h}*(\zeta_{\vartheta}{\mathbf u}^{\delta})]dxdt \notag\\
	& & & =-\int_{{\mathbb R}^{1}}\int_{\Omega}[\eta_{h}*(\rho^{\delta}{\mathbf u}^{\delta})]\cdot\nabla(\zeta_{\theta}{\mathbf u}^{\delta})(t+h)\cdot[\eta_{h}*(\zeta_{\theta}{\mathbf u}^{\delta})]dxdt\notag\\
	& & &+\frac{1}{2}\int_{{\mathbb R}^{1}}\int_{\Omega}\rho^{\delta}{\mathbf u}^{\delta}\cdot\nabla[\eta_{h}*(\zeta_{\theta}{\mathbf u}^{\delta})]^{2}dxdt.
\end{align}
Further, one finds that
\begin{align*}
	\lim_{\vartheta\to0}\lim_{h\to0}(\int_{{\mathbb R}^{1}}\int_{\Omega}&(\rho^{\delta}{\mathbf u}^{\delta}\otimes{\mathbf u}^{\delta}:\nabla\varphi_{h,\vartheta}-[\eta_{h}*(\rho^{\delta}{\mathbf u}^{\delta})]\cdot\nabla(\zeta_{\theta}{\mathbf u}^{\delta})(t+h)\cdot[\eta_{h}*(\zeta_{\theta}{\mathbf u}^{\delta})])dxds\\
	&+\frac{1}{2}\int_{{\mathbb R}^{1}}\int_{\Omega}\rho^{\delta}{\mathbf u}^{\delta}\cdot\nabla[\eta_{h}*(\zeta_{\theta}{\mathbf u}^{\delta})]^{2}dxdt)=0.
\end{align*}

In conclusion, one can choose $\varphi=\varphi_{h,\theta}$ in \eqref{(1e)} and let $h\to0,\vartheta\to0$ to attain the goal that
\begin{equation*}
	\int_{0}^{\tau}\int_{\Omega}\nu(\phi^{\delta})
	\overline{(1+|{\mathbb D}{\mathbf u}^{\delta}|^{2})^{\frac{p-2}{2}}{\mathbb D}{\mathbf u}^{\delta}:{\mathbb D}{\mathbf u}^{\delta}}-
	\overline{(1+|{\mathbb D}{\mathbf u}^{\delta}|^{2})^{\frac{p-2}{2}}{\mathbb D}{\mathbf u}^{\delta}}:{\mathbb D}{\mathbf u}^{\delta})dxdt=0.
\end{equation*}	
So, it is deduced from Lemma \refeq{L2-6} that
\begin{equation*}
	\overline{(1+|{\mathbb D}{\mathbf u}^{\delta}|^{2})^{\frac{p-2}{2}}{\mathbb D}{\mathbf u}^{\delta}}=(1+|{\mathbb D}{\mathbf u}^{\delta}|^{2})^{\frac{p-2}{2}}{\mathbb D}{\mathbf u}^{\delta}.
\end{equation*}
\end{proof}

By the weak lower semi-continuity of the norm, one obtains the uniform estimates on $(\rho^{\delta},{\mathbf u}^{\delta},\phi^{\delta},\mu^{\delta}).$ That is, there exists a positive constant $C$, independent of $\delta$, such that
\begin{eqnarray}\label{E3-19}
\begin{cases}
	\|\rho^\delta\|_{L^{\infty}(0,T;L^{\infty}(\Omega))}\leqslant C,~~\|\partial_{t}\rho^\delta\|_{L^{\infty}(0,T;H^{-1}(\Omega))}\leqslant C,\\
	\|\sqrt{\rho^\delta}{\mathbf u}^\delta\|_{{\bf L}^{\infty}(0,T;{\bf L}^{2}(\Omega))}\leqslant C,
    ~\|{\mathbf u}^\delta\|_{{\bf L}^{p}(0,T;{\bf W}_{0}^{1,p}(\Omega))}\leqslant C,\\
	\|\sqrt[4]{\rho^\delta}\phi^\delta\|_{L^{\infty}(0,T;L^{4}(\Omega))}\leqslant C,~~\|\phi^\delta\|_{L^{\infty}(0,T;H^{1}(\Omega))}\leqslant C,~~
    \|\sqrt{\rho^\delta}\partial_t\phi^\delta\|_{L^{2}(0,T;L^{2}(\Omega))}\leqslant C,\\
    \|\phi^\delta\|_{ L^{2}(0,T;W^{2,\infty}(\Omega))}\leqslant C,~\|\phi^\delta\|_{L^{\infty}(0,T;H^{2}(\Omega))}\leqslant C,
    \|\partial_t\phi^\delta\|_{L^{2}(0,T;H^{1}(\Omega))}\leqslant C, \\
  \|\mu^\delta\|_{L^{2}(0,T;H^{2}(\Omega))}\leqslant C,~\|\mu^\delta\|_{L^{\infty}(0,T;H^{1}(\Omega))}\leqslant C,
  ~\|\sqrt{\rho^\delta}\mu^\delta\|_{L^{\infty}(0,T;L^{2}(\Omega))}\leqslant C.
		\end{cases}
	\end{eqnarray}
This allows us to take similar argument in the process $m\to+\infty$,  pass further to the limit as $\delta\to 0^{+}$, and obtain a limit  $(\rho,{\mathbf u},\phi,\mu)$ satisfying the problem \eqref{E1-1}-\eqref{E1-2} as stated in Theorem \ref{thm2-1} .

\section{Proof of Theorem \ref{thm2-2}}\label{S-4}

\subsection{Approximate scheme and initial date}

\quad Let $\Omega={\mathbb T}^3$ and $\frac{5}{2}<p<3.$ For fixed $\delta\in(0,1)$ and $m\in{\mathbb Z}^+,$ we are looking for $(\rho_{m},{\mathbf u}_{m},\phi_{m},\mu_{m})$ to solve the following system
\begin{eqnarray}\label{E4-3}
\begin{cases}
\partial_{t}\rho_{m}+{\mathbf u}_{m}\cdot\nabla\rho_{m}=0&\mbox{in}\;\Omega\times(0,T),\\
(\rho_{m}\partial_{t}{\mathbf u}_{m},{\mathbf w})+(\rho_{m}{\mathbf u}_{m}\cdot\nabla{\mathbf u}_{m},{\mathbf w})
    +(\nu(\phi_{m})(1+|{\mathbb D}{{\mathbf u}}_{m}|^{2})^{\frac{p-2}{2}}\mathbb{D}{\mathbf u}_{m},\nabla {\mathbf w})\\
~~~~~~~~~~~~~~~~=(\rho_{m}\mu_{m}\nabla\phi_{m},{\mathbf w})-(\rho_{m}\nabla\Psi(\phi_{m}),{\mathbf w})~~~~({\bf w}\in {\bf V}_m)&\mbox{in}\;(0,T),\\
(\rho_{m}\partial_{t}\phi_{m},w)+(\rho_{m} {\mathbf u}_{m}\cdot\nabla\phi_{m},w)+(\nabla\mu_{m},\nabla w)=0~~~~(w\in V_m)&\mbox{in}\;(0,T),\\
(\rho_{m}\mu_{m},w)=(\nabla\phi_{m},\nabla w)+(\rho_m\Psi^\prime(\phi_{m}),w)~~~~(w\in V_m)&\mbox{in}\;(0,T),\\
{\mathbf u}_{m}=0,\partial_{n}\mu_{m}=\partial_{n}\phi_{m}=0 &\mbox{ on }\partial\Omega\times (0,T),\\
\rho_{m}(\cdot,0)=\rho_{0\delta},{\mathbf u}_{m}(\cdot,0)={\mathbf u}_{0m},\phi_{m}(\cdot,0)=\phi_{0m}&\mbox{ in }\Omega.
\end{cases}
\end{eqnarray}
Here, $\{\rho_{0\delta}\},$ $\{{\mathbf u}_{0m}\}$ and $\{\phi_{0m}\}$ are defined as \eqref{E-3-1} and \eqref{E-3-3}, respectively.
According to Proposition A.1. in Appendix, the semi-Galerkin scheme \eqref{E4-3} admits a solution $(\rho_m,{\bf u}_m,\phi_m,\mu_m)$ with
$$\rho_m\in C^1(\overline{\Omega\times (0,T)}),~{\bf u}_m\in {\bf C}^1([0,T];{\bf V}_m),~~\phi_m\in C^1([0,T];V_m)~~\mbox{ and }~~\mu_m\in C([0,T];V_m).$$

To establish the existence of strong solutions to the problem \eqref{E1-1}-\eqref{E1-4}, one needs to get uniform estimates for the approximate solutions $(\rho_m,{\bf u}_m,\phi_m,\mu_m)$ in suitable higher norms.

\begin{lemma}\label{chuzhi4-1}
	Let $\rho_{0}\in L^{\infty}(\Omega),$  $\phi_{0}\in H^{1}(\Omega)$ with $\partial_{n}\phi_{0}=0$ on $\partial\Omega,$ and $\nabla\mu_{0}\in L^{2}(\Omega)$. Then
\begin{equation}\label{chuE4-2}
\phi_{0}\in H^{2}(\Omega),\;\;\mu_{0}\in H^{1}(\Omega).
\end{equation}
\end{lemma}

\begin{proof}
First, evaluating $\eqref{E1-1}_{3}$ at $t=0$, one arrives at
\begin{align*}
\|\Delta\phi_{0}\|_{L^2(\Omega)}&\leqslant\|\rho_{0}\|_{L^\infty(\Omega)}\|\mu_0\|_{L^2(\Omega)}
+\|\rho_{0}\|_{L^\infty(\Omega)}(\|\phi_{0}\|^{3}_{L^{6}(\Omega)}+\|\phi_{0}\|^{2}_{L^2(\Omega)})\\
&\leqslant\|\rho_{0}\|_{L^\infty(\Omega)}\|\mu_0\|_{L^2(\Omega)}
+\|\rho_{0}\|_{L^\infty(\Omega)}(\|\phi_{0}\|^{3}_{H^{1}(\Omega)}+\|\phi_{0}\|^{2}_{L^2(\Omega)}).
	\end{align*}
It is deduced from Lemma \ref{L-2-3} that
\begin{equation*}
	\|\phi_{0}\|_{H^2(\Omega)}\leqslant C.
\end{equation*}

Second, evaluating $\eqref{E1-1}_{4}$ at $t=0$, one arrives that
\begin{align*}
\|\rho_{0}\mu_0\|_{L^2(\Omega)}&\leqslant\|\Delta\phi_{0}\|_{L^2(\Omega)}
+\|\rho_{0}\|_{L^\infty(\Omega)}(\|\phi_{0}\|^{3}_{L^{6}(\Omega)}+\|\phi_{0}\|^{2}_{L^2(\Omega)})\\
		&\leqslant\|\phi_{0}\|_{H^2(\Omega)}+\|\rho_{0}\|_{L^\infty(\Omega)}(\|\phi_{0}\|^{3}_{H^{1}(\Omega)}
+\|\phi_{0}\|^{2}_{L^2(\Omega)}).
\end{align*}
So, $\|\rho_{0}\mu_{0}\|_{L^{2}(\Omega)}\leqslant C$. Furthermore, one obtains
\begin{equation*}
	\|\mu_{0}\|_{H^1(\Omega)}\leqslant C.
\end{equation*}
\end{proof}

\subsection{Uniform estimate}

\quad Following from the subsection \ref{subs-3-2}, one can obtain the following lemma.

\begin{lemma}
Under the condition of Theorem \ref{thm2-2}, there exists a positive constant $C$, independent of $m$ and $\delta$, such that
	\begin{eqnarray}\label{E4-1}
		\begin{cases}
			0<C\delta\leqslant\rho_{m}(x,t)\leqslant\rho^{*}+1\;(\forall(x,t)\in \Omega\times (0,T)),~~\|\partial_{t}\rho_{m}\|_{ L^{\infty}(0,T;H^{-1}(\Omega))}\leqslant C,\\
			\|\sqrt{\rho_{m}}{\mathbf u}_{m}\|_{{\bf L}^{\infty}(0,T;{\bf L}^{2}(\Omega))}\leqslant C,~~\|{\mathbf u}_{m}\|_{{\bf L}^{p}(0,T;{\bf W}_{0}^{1,p}(\Omega))}\leqslant C,\\
			\|\sqrt{\rho_{m}}\phi_{m}\|_{L^{\infty}(0,T;L^{2}(\Omega))}\leqslant C,~~\|\partial_t\phi_{m}\|_{L^{2}(0,T;H^{1}(\Omega))}\leqslant C,\\
			\|\phi_{m}\|_{L^{\infty}(0,T;H^{2}(\Omega))}\leqslant C,~~\|\phi_{m}\|_{ L^{2}(0,T;W^{2,\infty}(\Omega))}\leqslant C,~~\\			\|\mu_{m}\|_{L^{2}(0,T;H^{2}(\Omega))}\leqslant C,
			~~ \|\mu_{m}\|_{L^{\infty}(0,T;H^{1}(\Omega))}\leqslant C.
		\end{cases}
	\end{eqnarray}
\end{lemma}


\begin{lemma}\label{lemmma4-1}
Under the condition of Theorem \ref{thm2-2},  there exists a positive constant $C$, independent of $m$ and $\delta$, such that
\begin{align}\label{E4-9}
	&\frac{d}{dt}\left\{\frac{1}{2}\|\nabla\mu_{m}\|^{2}_{{{\bf L}^{2}(\Omega)}}+\int_{\Omega}\left(\rho_m \mu_{m} {\mathbf u}_{m}\cdot\nabla\phi_{m}
	+\frac{1}{p}\nu(\phi_{m})(1+|{\mathbb D}{\mathbf u}_{m}|^{2})^{\frac{p}{2}}\right)dx\right\}\notag\\
	&+\int_{\Omega}\left(|\nabla\partial_t\phi_{m}|^{2}+3\rho_m \phi_m^2|\partial_t\phi_{m}|^{2}+\rho_m |\partial_t\phi_{m}|^{2}+\rho_m  |\partial_t{\mathbf u}_{m}|^{2}\right)dx\notag\\
	&\leqslant\epsilon\|\nabla^{2}{\mathbf u}_{m}\|^{2}_{{\bf L}^{2}(\Omega)}+C(\epsilon)\|(1+|{\mathbb D}{\mathbf u}_{m}|^{2})^{\frac{p}{4}}\|^{\frac{2(5-p)}{3-p}}_{{\bf L}^{2}(\Omega)}
	+C(\epsilon)\|\nabla {\mathbf u}_{m}\|^{6}_{{\bf L}^{2}(\Omega)}\notag\\
	&+C\left(\|{\mathbf u}_{m}\|^{2}_{\mathbf{W}_{0}^{1,p}(\Omega)}+1)(\|\phi_{m}\|^{8}_{H^{2}(\Omega)}+\|\mu_{m}\|^{4}_{H^{1}(\Omega)}\right)
\end{align}
holds for any fixed $\epsilon\in(0,1).$
\end{lemma}

\begin{proof}
 On one hand, taking ${\mathbf w}=\partial_t{\mathbf u}_{m}$ in $\eqref{E4-3}_2$ and integrating the resultant over $\Omega,$ one obtains that
\begin{align}\label{E4-2}
	&\int_{\Omega}\rho_{m}|\partial_t{\mathbf u}_{m}|^{2}dx+\int_{\Omega}\rho_m({\mathbf u}_{m}\cdot\nabla ){\mathbf u}_{m}\cdot \partial_t{\mathbf u}_{m}dx
     +\frac{1}{p}\frac{d}{dt}\int_{\Omega}(1+|{\mathbb D}{\mathbf u}_{m}|^{2})^{\frac{p}{2}}dx\\
     &-\frac{1}{p}\int_{\Omega}\nu^\prime(\phi_{m})\partial_t\phi_{m}(1+|{\mathbb D}{\mathbf u}_{m}|^{2})^{\frac{p}{2}}dx
	=\int_{\Omega}\rho_{m}\mu_{m}\nabla\phi_{m}\cdot \partial_t{\mathbf u}_{m}dx-\int_{\Omega}\rho_{m}\nabla\Psi(\phi_{m})\cdot \partial_t{\mathbf u}_{m}dx.\notag
\end{align}
On the other hand, taking $w=\partial_t\mu_{m}$ in $\eqref{E4-3}_3$ and integrating the resultant over $\Omega,$ one gets that
\begin{equation}\label{Equation4-6}
   \frac{1}{2}\frac{d}{dt}\|\nabla\mu_{m}\|^{2}_{{\bf L}^{2}(\Omega)}
   +\underset{I_{1}}{\underbrace{\int_{\Omega}\rho_{m}\partial_t\phi_{m}\partial_t\mu_{m}dx}}
   +\underset{I_{2}}{\underbrace{\int_{\Omega}\rho _{m}{\mathbf u}_{m}\cdot\nabla\phi_{m}\partial_t\mu_{m}dx}}=0.
\end{equation}
It is deduced from $\eqref{E4-3}_1$ and $\eqref{E4-3}_4$ that
\begin{align*}
I_{1}&=\int_{\Omega}\rho_{m}\partial_t\phi_{m}\partial_t\mu_{m}dx\\
&     =\int_{\Omega}\left(-\Delta\partial_t\phi_{m}\partial_t\phi_{m}
+[\rho_m\Psi^{\prime}(\phi_{m})]_{t}\partial_t\phi_{m}+div(\rho {\mathbf u}_{m})\mu_{m}\partial_t\phi_{m}\right)dx\\
&=\int_{\Omega}\left(|\nabla\partial_t\phi_{m}|^{2}+\partial_t\rho_m\Psi^{\prime}(\phi_{m})\partial_t\phi_{m}
+\rho_m\Psi^{\prime\prime}(\phi_{m})|\partial_t\phi_{m}|^{2}
       -\rho_m{\mathbf u}_{m}\cdot\nabla(\mu_{m}\partial_t\phi_{m})\right)dx\\
	&=\int_{\Omega}|\nabla\partial_t\phi_{m}|^{2}dx-\int_{\Omega}div(\rho {\mathbf u}_{m})\Psi^{\prime}(\phi_{m})\partial_t\phi_{m}dx\\
    &+\int_{\Omega}\rho_m\Psi^{\prime\prime}(\phi_{m})|\partial_t\phi_{m}|^{2}dx-\int_{\Omega}\rho_m\partial_t\phi_{m} {\mathbf u}_{m}\cdot\nabla\mu_{m}dx-\int_{\Omega}\rho_m\mu_{m} {\mathbf u}_{m}\cdot\nabla\partial_t\phi_{m}dx\\
	&=\int_{\Omega}|\nabla\partial_t\phi_{m}|^{2}dx
+\int_{\Omega}\rho_m{\mathbf u}_{m}\Psi^{\prime\prime}(\phi_{m})\partial_t\phi_{m}{\mathbf u}_{m}\cdot\nabla\phi_{m}dx
     +\int_{\Omega}\rho_m {\mathbf u}_{m}\Psi^{\prime}(\phi_{m}){\mathbf u}_{m}\cdot\nabla\partial_t\phi_{m}dx\\
    &+\int_{\Omega}\rho_m\Psi^{\prime\prime}(\phi_{m})|\partial_t\phi_{m}|^{2}dx
    -\int_{\Omega}\rho\partial_t\phi_{m} {\mathbf u}_{m}\cdot\nabla\mu_{m}dx-\int_{\Omega}\rho\mu_{m} {\mathbf u}_{m}\cdot\nabla\partial_t\phi_{m}dx,\\
I_{2}&=\int_{\Omega}\rho _{m}{\mathbf u}_{m}\cdot\nabla\phi_{m}\partial_t\mu_{m}dx\\
&=\frac{d}{dt}\int_{\Omega}\rho_m\mu_{m} {\mathbf u}_{m}\cdot\nabla\phi_{m} dx
-\int_{\Omega}\partial_t\rho_m\mu_{m} {\mathbf u}_{m}\cdot\nabla\phi_{m}dx\\
&+\int_{\Omega}\rho_m\mu_{m} \partial_t{\mathbf u}_{m}\nabla\phi_{m}dx+\int_{\Omega}\rho_m\mu_{m} {\mathbf u}_{m}\nabla\partial_t\phi_{m}dx\\
&=\frac{d}{dt}\int_{\Omega}\rho_m\mu_{m} {\mathbf u}_{m}\cdot\nabla\phi_{m} dx
+\int_{\Omega}\left(div(\rho_m{\mathbf u}_{m})\mu_{m} {\mathbf u}_{m}\cdot\nabla\phi_{m}
 -\rho_m\mu_{m} \partial_t{\mathbf u}_{m}\nabla\phi_{m}-\rho\mu_{m} {\mathbf u}_{m}\nabla\partial_t\phi_{m}\right)dx\\
&=\frac{d}{dt}\int_{\Omega}\rho_m\mu_{m}{\mathbf u}_{m}\cdot\nabla\phi_{m} dx
-\int_{\Omega}\rho_m({\mathbf u}_{m}\cdot\nabla\mu_{m})({\mathbf u}_{m}\cdot\nabla\phi_{m})dx\\
&-\int_{\Omega}\rho_m\mu_{m} {\mathbf u}_{m}\cdot\nabla({\mathbf u}_{m}\cdot\nabla\phi_{m})dx-\int_{\Omega}\rho_m\mu_{m} \partial_t{\mathbf u}_{m}\nabla\phi_{m}dx-\int_{\Omega}\rho_m\mu_{m} {\mathbf u}_{m}\nabla\partial_t\phi_{m}dx.
\end{align*}
Thus, it follows from \eqref{Equation4-6} that
\begin{align}\label{E4-4}
	&\frac{d}{dt}\left\{\frac{1}{2}\|\nabla\mu_{m}\|^{2}_{{{\bf L}^{2}(\Omega)}}
+\int_{\Omega}\rho_m \mu_{m} {\mathbf u}_{m}\cdot\nabla\phi_{m} dx\right\}
	+\int_{\Omega}(|\nabla\partial_t\phi_{m}|^{2}+\rho_m \Psi^{\prime\prime}(\phi_{m})|\partial_t\phi_{m}|^{2})dx\\
	&=\int_{\Omega}\left(-\rho_m \Psi^{\prime\prime}(\phi_{m})\partial_t\phi_{m}{\mathbf u}_{m}\cdot\nabla\phi_{m}
      -\rho_m \Psi^{\prime}(\phi_{m}){\mathbf u}_{m}\cdot\nabla\partial_t\phi_{m}+\rho_m \partial_t\phi_{m}{\mathbf u}_{m}\cdot\nabla\mu_{m}\right.\notag\\
    &\left.+\rho_m ({\mathbf u}_{m}\cdot\nabla \mu_{m})({\mathbf u}_{m}\cdot\nabla\phi_{m})+\rho_m \mu_{m} {\mathbf u}_{m}\cdot\nabla({\mathbf u}_{m}\cdot\nabla\phi_{m})+\rho_m \mu_{m} \partial_t{\mathbf u}_{m}\cdot\nabla\phi_{m}+2\rho_m \mu_{m} {\mathbf u}_{m}\cdot\nabla\partial_t\phi_{m}\right)dx\notag.
\end{align}
Further, one takes $w=\partial_t\phi_{m}$ in $\eqref{E4-3}_4$ and integrates the resultant over $\Omega$ to arrive at
\begin{equation}\label{E4-5}
	\int_{\Omega}\rho_m |\partial_t\phi_{m}|^{2}dx
=-\int_{\Omega}\left(\rho_m  {\mathbf u}_{m}\cdot\nabla\phi_{m}\partial_t\phi_{m}+\nabla\mu_{m}\cdot\nabla\partial_t\phi_{m}\right)dx.
\end{equation}
Since $\Psi^{''}(s)=3s^{2}-1,$ it is deduced from \eqref{E4-4} and \eqref{E4-5} that
\begin{align}\label{E4-6}
	&\frac{d}{dt}\{\frac{1}{2}\|\nabla\mu_{m}\|^{2}_{{{\bf L}^{2}(\Omega)}}+\int_{\Omega}\rho_m \mu_{m} {\mathbf u}_{m}\cdot\nabla\phi_{m} dx\}
+\int_{\Omega}\left(|\nabla\partial_t\phi_{m}|^{2}+3\rho_m \phi_m^2|\partial_t\phi_{m}|^{2}+\rho_m |\partial_t\phi_{m}|^{2}\right)dx\notag\\
	&=\int_{\Omega}(-3\rho_m \phi^{2}_{m}\partial_t\phi_{m}{\mathbf u}_{m}\cdot\nabla\phi_{m}-\rho_{m}\partial_t\phi_{m}{\mathbf u}_{m}\cdot\nabla\phi_{m}-\rho_m  {\bf u}_m\cdot\nabla\phi_{m}\partial_t\phi_{m}\notag\\
	&-\rho_m \Psi^{\prime}(\phi_{m}){\mathbf u}_{m}\cdot\nabla\partial_t\phi_{m}+\rho_m \partial_t\phi_{m}{\mathbf u}_{m}\cdot\nabla\mu_{m}+\rho_m ({\mathbf u}_{m}\cdot\nabla \mu_{m})({\mathbf u}_{m}\cdot\nabla\phi_{m})\notag\\
	&+\rho_m \mu_{m} {\mathbf u}_{m}\cdot\nabla({\mathbf u}_{m}\cdot\nabla\phi_{m})+\rho_m \mu_{m} \partial_t{\mathbf u}_{m}\cdot\nabla\phi_{m}+2\rho_m \mu_{m} {\bf u}_m\cdot \nabla\partial_t\phi_{m}-2\nabla\mu_{m}\cdot\nabla\partial_t\phi_{m})dx.
\end{align}
Now, one adds \eqref{E4-2} and \eqref{E4-6} to obtain that
\begin{align}\label{E4-7}
  &\frac{d}{dt}\left\{\frac{1}{2}\|\nabla\mu_{m}\|^{2}_{{{\bf L}^{2}(\Omega)}}+\int_{\Omega}\left(\rho_m \mu_{m} {\mathbf u}_{m}\cdot\nabla\phi_{m}
    +\frac{1}{p}\nu(\phi_{m})(1+|{\mathbb D}{\mathbf u}_{m}|^{2})^{\frac{p}{2}}\right)dx\right\}\notag\\
  &~~+\int_{\Omega}\left(|\nabla\partial_t\phi_{m}|^{2}+3\rho_m \phi_m^2|\partial_t\phi_{m}|^{2}+\rho_m |\partial_t\phi_{m}|^{2}
   +\rho_m  |\partial_t{\mathbf u}_{m}|^{2}\right)dx\notag\\
  &~~=\underset{I_{1}}{\underbrace{\int_{\Omega}-3\rho_m \phi^{2}_{m}\partial_t\phi_{m}{\mathbf u}_{m}\cdot\nabla\phi_{m}dx}}
  +\underset{I_{2}}{\underbrace{\int_{\Omega}-\rho_{m}\partial_t\phi_{m}{\mathbf u}_{m}\cdot\nabla\phi_{m}dx}}\notag\\
  &+\underset{I_{3}}{\underbrace{\int_{\Omega}-\rho_m ({\mathbf u}_m\cdot\nabla {\mathbf u}_m)\cdot \partial_t{\mathbf u}_{m}dx}}+\underset{I_{4}}{\underbrace{\frac{1}{p}\int_{\Omega}\nu^{\prime}(\phi_{m})\partial_t\phi_{m}(1+|{\mathbb D}{\mathbf u}_{m}|^{2})^{\frac{p}{2}}dx}}\notag\\
   &+\underset{I_{5}}{\underbrace{\int_{\Omega}\rho_m \mu_{m}\nabla\phi_{m}\cdot \partial_t{\mathbf u}_{m}dx}}
   +\underset{I_{6}}{\underbrace{\int_{\Omega}-\rho_m \nabla\Psi(\phi_{m})\cdot \partial_t{\mathbf u}_{m}dx}}\notag\\
  &+\underset{I_{7}}{\underbrace{\int_{\Omega}-\rho_m{\mathbf u}_m\cdot\nabla\phi_{m}\partial_t\phi_{m}dx}}
    +\underset{I_{8}}{\underbrace{\int_{\Omega}-\rho_m \Psi^\prime(\phi_{m}){\mathbf u}_m\cdot\nabla\partial_t\phi_{m}dx}}\notag\\
    &+\underset{I_{9}}{\underbrace{\int_{\Omega}\rho_m \partial_t\phi_{m}{\mathbf u}_{m}\cdot\nabla\mu_{m} dx}}+\underset{I_{10}}{\underbrace{\int_{\Omega}\rho_m ({\mathbf u}_{m}\cdot\nabla \mu_{m})({\mathbf u}_{m}\cdot\nabla\phi_{m})dx}}\notag\\
    & +\underset{I_{11}}{\underbrace{\int_{\Omega}\rho_m \mu_{m} {\mathbf u}_{m}\cdot\nabla {\mathbf u}_{m}\cdot\nabla\phi_{m}dx}}
     +\underset{I_{12}}{\underbrace{\int_{\Omega}\rho_m \mu_{m} {\mathbf u}_{m}\cdot({\mathbf u}_{m}\cdot\nabla^{2}\phi_{m})dx}}\notag\\
	&+\underset{I_{13}}{\underbrace{\int_{\Omega}\rho_m \mu_{m} \partial_t{\mathbf u}_{m}\cdot\nabla\phi_{m} dx}}
     +\underset{I_{14}}{\underbrace{\int_{\Omega}-2\nabla\mu_{m}\cdot\nabla\partial_t\phi_{m}dx}}
     +\underset{I_{15}}{\underbrace{\int_{\Omega}2\rho_m \mu_{m} {\mathbf u}_{m}\cdot\nabla\partial_t\phi_{m}dx}}.
\end{align}
Each term on the right-hand side of \eqref{E4-7} is estimated as follows
\begin{align*}
&|I_{1}|=|\int_{\Omega}-3\rho_m \phi^{2}_{m}\partial_t\phi_{m}{\mathbf u}_{m}\cdot\nabla\phi_{m}dx|\\
&\;\;\;\;\;\leqslant C\|\sqrt{\rho_{m}}\phi_{m}\partial_{t}\phi_{m}\|_{L^{2}(\Omega)}\|\phi_{m}\|_{L^{\infty}(\Omega)}
            \|{\mathbf u}_{m}\|_{\mathbf{L}^{3}(\Omega)}\|\nabla\phi_{m}\|_{{\bf L}^{6}(\Omega)}\\
&\;\;\;\;\;\leqslant C\|\sqrt{\rho_{m}}\phi_{m}\partial_{t}\phi_{m}\|_{L^{2}(\Omega)}
             \|{\mathbf u}_{m}\|_{\mathbf{W}_{0}^{1,p}(\Omega)}\|\phi_{m}\|^{2}_{H^{2}(\Omega)}\\
&\;\;\;\;\;\leqslant\frac{1}{2}\|\sqrt{\rho_{m}}\phi_{m}\partial_{t}\phi_{m}\|^{2}_{L^{2}(\Omega)}
             +C\|{\mathbf u}_{m}\|^{2}_{\mathbf{W}_{0}^{1,p}(\Omega)}\|\phi_{m}\|^{4}_{H^{2}(\Omega)},\\
&|I_{2}|=|\int_{\Omega}-\rho_{m}\partial_t\phi_{m}{\mathbf u}_{m}\cdot\nabla\phi_{m}dx|\\
&\;\;\;\;\;\leqslant\|\sqrt{\rho_{m}}\partial_{t}\phi_{m}\|_{L^{2}(\Omega)}\|{\mathbf u}_{m}\|_{{\mathbf L}^{3}(\Omega)}
            \|\nabla\phi_{m}\|_{{\bf L}^{6}(\Omega)}\\
&\;\;\;\;\;\leqslant \|\sqrt{\rho_{m}}\partial_{t}\phi_{m}\|_{L^{2}(\Omega)}
             \|{\mathbf u}_{m}\|_{\mathbf{W}_{0}^{1,p}(\Omega)}\|\phi_{m}\|_{H^{2}(\Omega)}\\
&\;\;\;\;\;\leqslant\frac{1}{4}\|\sqrt{\rho_{m}}\partial_{t}\phi_{m}\|^{2}_{L^{2}(\Omega)}
            +C\|{\mathbf u}_{m}\|^{2}_{\mathbf{W}_{0}^{1,p}(\Omega)}\|\phi_{m}\|^{2}_{H^{2}(\Omega)},\\
&|I_{3}|=|\int_{\Omega}-\rho_m ({\mathbf u}_m\cdot\nabla {\mathbf u}_m)\cdot \partial_t{\mathbf u}_{m}dx|\\
&\;\;\;\;\;\leqslant C\|\sqrt{\rho_m }\partial_t{\mathbf u}_m\|_{{\bf L}^{2}(\Omega)}\|{\mathbf u}_m\|_{{\bf L}^{6}(\Omega)}
           \|\nabla {\mathbf u}_{m}\|_{{\bf L}^{3}(\Omega)}\\
&\;\;\;\;\;\leqslant C\|\nabla {\mathbf u}_{m}\|_{{\bf L}^{2}(\Omega)}^{\frac{3}{2}}\|\nabla^{2}{\mathbf u}_{m}
           \|_{{\bf L}^{2}(\Omega)}^{\frac{1}{2}}\|\sqrt{\rho_m }\partial_t{\mathbf u}_m\|_{{\bf L}^{2}(\Omega)}\\
&\;\;\;\;\;\leqslant\frac{\epsilon}{16}\|\nabla^{2}{\mathbf u}_{m}\|^{2}_{L^{2}(\Omega)}+\frac{1}{8}
          \|\sqrt{\rho_m }\partial_t{\mathbf u}_m\|^{2}_{L^{2}(\Omega)}+C(\epsilon)\|\nabla {\mathbf u}_{m}\|^{6}_{L^{2}(\Omega)},\\
&|I_{4}|=|\frac{1}{p}\int_{\Omega}\nu^{\prime}(\phi_{m})\partial_t\phi_{m}(1+|{\mathbb D}{\mathbf u}_{m}|^{2})^{\frac{p}{2}}dx|\\
&\;\;\;\;\;\leqslant C\int_{\Omega}|\partial_t\phi_{m}||{\mathbb D}{\mathbf u}_{m}|^{p}dx\\
&\;\;\;\;\;\leqslant C\|\partial_t\phi_{m}\|_{L^{6}(\Omega)}\|{\mathbb D}{\mathbf u}_{m}\|^{p}_{{\bf L}^{\frac{6p}{5}}(\Omega)}\\
&\;\;\;\;\;\leqslant C\|\partial_t\phi_{m}\|_{L^{6}(\Omega)}\|\nabla^{2}{\mathbf u}_{m}\|^{\frac{p}{6-p}}_{{\bf L}^{2}(\Omega)}
          \|\nabla{\mathbb D}{\mathbf u}_{m}\|^{\frac{p(5-p)}{6-p}}_{{\bf L}^{p}(\Omega)}\\
&\;\;\;\;\;\leqslant\frac{\epsilon}{16}\|\nabla^{2}{\mathbf u}_{m}\|^{2}_{{\bf L}^{2}(\Omega)}
     +C(\epsilon)(\|\partial_t\phi_{m}\|^{2}_{H^{1}(\Omega)}
     +\|(1+|{\mathbb D}{\mathbf u}_{m}|^{2})^{\frac{p}{4}}\|^{\frac{2(5-p)}{3-p}}_{{\bf L}^{2}(\Omega)}),\\
&|I_{5}|=|\int_{\Omega}\rho_m \mu_{m}\nabla\phi_{m}\cdot \partial_t{\mathbf u}_{m}dx|\\
&\;\;\;\;\;\leqslant\|\sqrt{\rho_m }\partial_t{\mathbf u}_m\|_{{\bf L}^{2}(\Omega)}\|\nabla\phi_{m}\|_{{\bf L}^{6}(\Omega)}
     \|\mu_{m}\|_{L^{3}(\Omega)}\\
&\;\;\;\;\;\leqslant\|\sqrt{\rho_m }\partial_t{\mathbf u}_m\|_{L^{2}(\Omega)}\|\phi_{m}\|_{H^{2}(\Omega)}\|\mu_{m}\|_{H^{1}(\Omega)}\\
&\;\;\;\;\;\leqslant\frac{1}{8}\|\sqrt{\rho_m }\partial_t{\mathbf u}_m\|^{2}_{L^{2}(\Omega)}
         +C\|\phi_{m}\|^{2}_{H^{2}(\Omega)}\|\mu_{m}\|^{2}_{H^{1}(\Omega)},\\
&|I_{6}|=|\int_{\Omega}-\rho_m \nabla\Psi(\phi_{m})\cdot \partial_t{\mathbf u}_{m}dx|\\
&\;\;\;\;\;\leqslant\|\sqrt{\rho_m }\partial_t{\mathbf u}_m\|_{L^{2}(\Omega)}
           \|\nabla\phi_{m}\|_{{\bf L}^{3}(\Omega)}\|\phi_{m}^{3}-\phi_{m}\|_{L^{6}(\Omega)}\\
&\;\;\;\;\; \leqslant\|\sqrt{\rho_m }\partial_t{\mathbf u}_m\|_{L^{2}(\Omega)}(\|\phi_{m}\|^{4}_{H^{2}(\Omega)}+\|\phi_{m}\|^{2}_{H^{2}(\Omega)})\\
&\;\;\;\;\;\leqslant\frac{1}{8}\|\sqrt{\rho_m }\partial_t{\mathbf u}_m\|^{2}_{L^{2}(\Omega)}
          +C(\|\phi_{m}\|^{4}_{H^{2}(\Omega)}+\|\phi_{m}\|^{2}_{H^{2}(\Omega)})^{2},\\
&|I_{7}|=|\int_{\Omega}-\rho_m   {\mathbf u}_{m}\cdot\nabla\phi_{m}\partial_t\phi_{m}dx|\\
&\;\;\;\;\;\leqslant\|\sqrt{\rho_m }\partial_t\phi_{m}\|_{L^{2}(\Omega)}\|\nabla\phi_{m}\|_{{\bf L}^{6}(\Omega)}
           \|{\mathbf u}_{m}\|_{{\bf W}_{0}^{1,p}(\Omega)}\\
&\;\;\;\;\; \leqslant\|\sqrt{\rho_m }\partial_t\phi_{m}\|_{L^{2}(\Omega)}\|\phi_{m}\|_{H^{2}(\Omega)}\|{\mathbf u}_{m}\|_{H^{1}(\Omega)}\\
&\;\;\;\;\;\leqslant\frac{1}{4}\|\sqrt{\rho_{m}}\partial_{t}\phi_{m}\|^{2}_{L^{2}(\Omega)}
          +C\|\phi_{m}\|^{2}_{H^{2}(\Omega)}\|{\mathbf u}_{m}\|^{2}_{\mathbf{W}_0{}^{1,p}(\Omega)},\\
&|I_{8}|=|\int_{\Omega}-\rho_m \Psi^\prime(\phi_{m}){\mathbf u}_m\cdot\nabla\partial_t\phi_{m}dx|\\
&\;\;\;\;\;\leqslant C\|\nabla\partial_t\phi_{m}\|_{{\bf L}^{2}(\Omega)}\|\phi_{m}^{3}
           -\phi_{m}\|_{L^{6}(\Omega)}\|{\mathbf u}_{m}\|_{{\bf L}^{3}(\Omega)}\\
&\;\;\;\;\;\leqslant C\|\partial_t\phi_{m}\|_{H^{1}(\Omega)}(\|\phi_{m}\|^{3}_{H^{2}(\Omega)}
            +\|\phi_{m}\|_{H^{2}(\Omega)})\|{\mathbf u}_{m}\|_{H^{1}(\Omega)}\\
&\;\;\;\;\;\leqslant C\|\partial_t\phi_{m}\|^{2}_{H^{1}(\Omega)}+C\|{\mathbf u}_{m}\|^{2}_{{\bf W}_{0}^{1,p}(\Omega)}(\|\phi_{m}\|^{3}_{H^{2}(\Omega)}
         +\|\phi_{m}\|_{H^{2}(\Omega)})^{2},\\
&|I_{9}|=|\int_{\Omega}\rho_m \partial_t\phi_{m}{\mathbf u}_{m}\cdot\nabla\mu_{m} dx|\\
&\;\;\;\;\;\leqslant \|\rho_{m}\|_{L^{\infty}(\Omega)}\|\partial_t\phi_{m}\|_{L^{6}(\Omega)}\|\nabla\mu_{m}\|_{{\bf L}^{2}(\Omega)}
           \|{\mathbf u}_{m}\|_{{\bf L}^{3}(\Omega)}\\
& \;\;\;\;\;\leqslant C\|\partial_t\phi_{m}\|_{H^{1}(\Omega)}\|\mu_{m}\|_{H^{1}(\Omega)}\|{\mathbf u}_{m}\|_{{\bf W}_{0}^{1,p}(\Omega)}\\
& \;\;\;\;\;\leqslant C\|\partial_t\phi_{m}\|^{2}_{H^{1}(\Omega)}
            +C\|{\mathbf u}_{m}\|^{2}_{{\bf W}_{0}^{1,p}(\Omega)}\|\mu_{m}\|^{2}_{H^{1}(\Omega)},\\
&|I_{10}|=|\int_{\Omega}\rho_m ({\mathbf u}_{m}\cdot\nabla \mu_{m})({\mathbf u}_{m}\cdot\nabla\phi_{m})dx|\\
&\;\;\;\;\;\leqslant\|{\mathbf u}_{m}\|^{2}_{L^{6}(\Omega)}\|\nabla\mu_{m}\|_{{\bf L}^{2}(\Omega)}\|\nabla\phi_{m}\|_{{\bf L}^{6}(\Omega)}\\
&\;\;\;\;\;\leqslant C\|{\mathbf u}_{m}\|_{{\bf W}_{0}^{1,p}(\Omega)}\|\mu_{m}\|_{H^{1}(\Omega)}\|\phi_{m}\|_{H^{2}(\Omega)}\\
&\;\;\;\;\;\leqslant C\|{\mathbf u}_{m}\|^{2}_{\mathbf{W}_{0}^{1,p}(\Omega)}(\|\mu_{m}\|^{2}_{H^{1}(\Omega)}+\|\phi_{m}\|^{2}_{H^{2}(\Omega)}),\\
&|I_{11}|=|\int_{\Omega}\rho_m \mu_{m} {\mathbf u}_{m}\cdot\nabla {\mathbf u}_{m}\cdot\nabla\phi_{m}dx|\\
&\;\;\;\;\;\leqslant C\|\mu_{m}\|_{L^{6}(\Omega)}\|{\mathbf u}_{m}\|_{{\bf L}^{6}(\Omega)}
        \|\nabla {\mathbf u}_{m}\|_{{\bf L}^{2}(\Omega)}\|\nabla\phi_{m}\|_{{\bf L}^{6}(\Omega)}\\
&\;\;\;\;\;\leqslant C\|\mu_{m}\|_{H^{1}(\Omega)}\|{\mathbf u}_{m}\|_{{\bf W}_{0}^{1,p}(\Omega)}\|\phi_{m}\|_{H^{2}(\Omega)}\\
&\;\;\;\;\;\leqslant C\|{\mathbf u}_{m}\|^{2}_{\mathbf{W}_{0}^{1,p}(\Omega)}(\|\mu_{m}\|^{2}_{H^{1}(\Omega)}+\|\phi_{m}\|^{2}_{H^{2}(\Omega)}),\\
&|I_{12}|=|\int_{\Omega}\rho_m \mu_{m} {\mathbf u}_{m}\cdot({\mathbf u}_{m}\cdot\nabla^{2}\phi_{m})dx|\\
&\;\;\;\;\;\leqslant C\|\mu_{m}\|_{L^{6}(\Omega)}\|{\mathbf u}_{m}\|^{2}_{L^{6}(\Omega)}\|\nabla^{2}\phi_{m}\|_{{\bf L}^{2}(\Omega)}\\
&\;\;\;\;\;\leqslant C\|\mu_{m}\|_{H^{1}(\Omega)}\|{\mathbf u}_{m}\|_{{\bf W}_{0}^{1,p}(\Omega)}\|\phi_{m}\|_{H^{2}(\Omega)}\\
&\;\;\;\;\;\leqslant C\|{\mathbf u}_{m}\|^{2}_{\mathbf{W}_{0}^{1,p}(\Omega)}(\|\mu_{m}\|^{2}_{H^{1}(\Omega)}+\|\phi_{m}\|^{2}_{H^{2}(\Omega)}),\\
&|I_{13}|=|\int_{\Omega}\rho_m \mu_{m} \partial_t{\mathbf u}_{m}\cdot\nabla\phi_{m}dx|\\
&\;\;\;\;\;\leqslant C\|\sqrt{\rho_m }\partial_t{\mathbf u}_m\|_{{\bf L}^{2}(\Omega)}\|\mu_{m}\|_{L^{6}(\Omega)}
           \|\nabla\phi_{m}\|_{{\bf L}^{3}(\Omega)}\\
&\;\;\;\;\; \leqslant C\|\sqrt{\rho_m }\partial_t{\mathbf u}_m\|_{L^{2}(\Omega)}\|\mu_{m}\|_{H^{1}(\Omega)}\|\phi_{m}\|_{H^{2}(\Omega)}\\
&\;\;\;\;\; \leqslant\frac{1}{8}\|\sqrt{\rho_m }\partial_t{\mathbf u}_m\|^{2}_{L^{2}(\Omega)}+C(\|\mu_{m}\|^{4}_{H^{1}(\Omega)}
           +\|\phi_{m}\|^{4}_{H^{2}(\Omega)}),\\
&|I_{14}|=|\int_{\Omega}-2\nabla\mu_{m}\cdot\nabla\partial_t\phi_{m}dx|\\
&\;\;\;\;\;\leqslant 2\|\nabla\mu_{m}\|_{{\bf L}^{2}(\Omega)}\|\nabla\partial_t\phi_{m}\|_{{\bf L}^{2}(\Omega)}\\
&\;\;\;\;\;\leqslant C\|\partial_t\phi_{m}\|^{2}_{H^{1}(\Omega)}+C\|\mu_{m}\|^{2}_{H^{1}(\Omega)},\\
&|I_{15}|=|\int_{\Omega}2\rho_m \mu_{m} {\mathbf u}_{m}\cdot\nabla\partial_t\phi_{m}dx|\\
&\;\;\;\;\;\leqslant C\|\mu_{m}\|_{L^{6}(\Omega)}\|{\mathbf u}_{m}\|_{{\bf L}^{3}(\Omega)}\|\nabla\partial_t\phi_{m}\|_{{\bf L}^{2}(\Omega)}\\
&\;\;\;\;\;\leqslant C\|\mu_{m}\|_{H^{1}(\Omega)}\|{\mathbf u}_{m}\|_{H^{1}(\Omega)}\|\partial_t\phi_{m}\|_{H^{1}(\Omega)}\\
&\;\;\;\;\;\leqslant C\|\partial_t\phi_{m}\|^{2}_{H^{1}(\Omega)}
+C\|{\mathbf u}_{m}\|^{2}_{\mathbf{W}_{0}^{1,p}(\Omega)}\|\mu_{m}\|^{2}_{H^{1}(\Omega)},
\end{align*}
where $\epsilon\in(0,1)$ is any fixed. Then,
\begin{align*}
    &\frac{d}{dt}\left\{\frac{1}{2}\|\nabla\mu_{m}\|^{2}_{{{\bf L}^{2}(\Omega)}}+\int_{\Omega}\left(\rho_m \mu_{m} {\mathbf u}_{m}\cdot\nabla\phi_{m}
    +\frac{1}{p}\nu(\phi_{m})(1+|{\mathbb D}{\mathbf u}_{m}|^{2})^{\frac{p}{2}}\right)dx\right\}\notag\\
   &+\int_{\Omega}\left(|\nabla\partial_t\phi_{m}|^{2}+3\rho_m \phi_m^2|\partial_t\phi_{m}|^{2}+\rho_m |\partial_t\phi_{m}|^{2}+\rho_m  |\partial_t{\mathbf u}_{m}|^{2}\right)dx\notag\\
  &\leqslant\epsilon\|\nabla^{2}{\mathbf u}_{m}\|^{2}_{{\bf L}^{2}(\Omega)}+C(\epsilon)\|(1+|{\mathbb D}{\mathbf u}_{m}|^{2})^{\frac{p}{4}}\|^{\frac{2(5-p)}{3-p}}_{{\bf L}^{2}(\Omega)}
+C(\epsilon)\|\nabla {\mathbf u}_{m}\|^{6}_{{\bf L}^{2}(\Omega)}\notag\\
&+C(\|{\mathbf u}_{m}\|^{2}_{\mathbf{W}_{0}^{1,p}(\Omega)}+1)(\|\phi_{m}\|^{8}_{H^{2}(\Omega)}+\|\mu_{m}\|^{4}_{H^{1}(\Omega)}),
\end{align*}
where $\epsilon\in(0,1)$ is any fixed.
\end{proof}

\begin{lemma}\label{lemmma4-2}
Under the condition of Theorem \ref{thm2-2},  there exists a position constant $C$, independent of $m$ and $\delta$, such that
\begin{align}\label{E4-10}
&\int_{\Omega}\nu(\phi_{m})(1+|{\mathbb D}{\mathbf u}_{m}|^{2})^{\frac{p-2}{2}}|\nabla {\mathbb D}{\mathbf u}_{m}|^{2}dx\notag\\
&\leqslant C\|(1+|{\mathbb D}{\mathbf u}_{m}|^{2})^{\frac{p-2}{4}}{\mathbb D}{\mathbf u}_{m}\|^{2}_{{\bf L}^{2}(\Omega)}+C\|\sqrt{\rho}\partial_t{\mathbf u}_{m}\|^{2}_{{\bf L}^{2}(\Omega)}
\notag\\
&+C\left(\|(1+|{\mathbb D}{\mathbf u}_{m}|^{2})^{\frac{p-2}{4}}{\mathbb D}{\mathbf u}_{m}\|^{\frac{8}{p}}_{{\bf L}^{2}(\Omega)}
+\|\mu_{m}\|^{4}_{H^{1}(\Omega)}+\|\phi_{m}\|^{6}_{H^{2}(\Omega)}\right).
\end{align}
\end{lemma}

\begin{proof}
First, one takes $w=\Delta {\mathbf u}_{m}$ in $\eqref{E4-3}_2$ and gets that
\begin{align}\label{ee-49}
	&\int_{\Omega}\left(\rho_{m}\partial_{t}{\mathbf u}_{m}\cdot\Delta {\mathbf u}_{m}+\rho_{m} ({\mathbf u}_{m}\cdot\nabla {\mathbf u}_{m})\cdot\Delta {\mathbf u}_{m}-div(\nu(\phi_{m})(1+|\mathbb D{\mathbf u}_{m}|^{2})^{\frac{p-2}{2}}\mathbb D{\mathbf u}_{m})\Delta {\mathbf u}_{m}+\nabla P_{m}\cdot\Delta {\mathbf u}_{m}\right)dx\notag\\
	&=\int_{\Omega}\left(\rho_{m}\mu_{m}\nabla\phi_{m}\cdot\Delta {\mathbf u}_{m}-\rho_{m}\nabla\Psi(\phi_{m})\cdot\Delta{\mathbf u}_{m}\right)dx.
\end{align}
Second, twice integration by parts over a periodic domain leads to
\begin{align}\label{ee-50}
&-\int_{\Omega}div(\nu(\phi_{m})(1+|{\mathbb D}{\mathbf u}_{m}|^{2})^{\frac{p-2}{2}}{\mathbb D}{\mathbf u}_{m})\cdot\Delta {\mathbf u}_mdx\notag\\
&=\int_{\Omega}-\nabla(\nu(\phi_{m})(1+|{\mathbb D}{\mathbf u}_{m}|^{2})^{\frac{p-2}{2}}{\mathbb D}{\mathbf u}_{m}):
        \nabla {\mathbb D}{\mathbf u}_{m}dx\notag\\
&=-\int_{\Omega}\nabla(\nu(\phi_{m}))(1+|{\mathbb D}{\mathbf u}_{m}|^{2})^{\frac{p-2}{2}}{\mathbb D}{\mathbf u}_{m}:
        \nabla {\mathbb D}{\mathbf u}_{m}dx\notag\\
&-\int_{\Omega}\nu(\phi_{m})(p-2)(1+|{\mathbb D}{\mathbf u}_{m}|^{2})^{\frac{p-4}{2}}|{\mathbb D}{\mathbf u}_{m}|^{2}
        |\nabla {\mathbb D}{\mathbf u}_{m}|^{2}dx\notag\\
&-\int_{\Omega}\nu(\phi_{m})(1+|{\mathbb D}{\mathbf u}_{m}|^{2})^{\frac{p-2}{2}}|\nabla {\mathbb D}{\mathbf u}_{m}|^{2}dx.
\end{align}
So, one can deduce from \eqref{ee-49} and \eqref{ee-50} that
\begin{align*}
&\int_{\Omega}\left(\nu(\phi_{m})(p-2)(1+|{\mathbb D}{\mathbf u}_{m}|^{2})^{\frac{p-4}{2}}|{\mathbb D}{\mathbf u}_{m}|^{2}
       |\nabla {\mathbb D}{\mathbf u}_{m}|^{2}+\nu(\phi_{m})(1+|{\mathbb D}{\mathbf u}_{m}|^{2})^{\frac{p-2}{2}}
       |\nabla {\mathbb D}{\mathbf u}_{m}|^{2}\right)dx\notag\\
&=-\int_{\Omega}\nabla\nu(\phi_{m})(1+|{\mathbb D}{\mathbf u}_{m}|^{2})^{\frac{p-2}{2}}{\mathbb D}{\mathbf u}_{m}:
      \nabla {\mathbb D}{\mathbf u}_{m}dx
     +\int_{\Omega}\rho_{m}\partial_t{\mathbf u}_{m}\cdot\Delta{\mathbf u}_{m}dx\notag\\
&+\int_{\Omega}\rho_{m}({\mathbf u}_{m}\cdot\nabla {\mathbf u}_{m})\cdot\Delta{\mathbf u}_{m}dx
     -\int_{\Omega}\rho_{m}\mu_{m}\nabla\phi_{m}\cdot\Delta  {\mathbf u}_{m}dx
     +\int_{\Omega}\rho_{m}\nabla\Psi(\phi_{m})\cdot\Delta  {\mathbf u}_{m}dx\notag\\
&\leqslant C\|(1+|{\mathbb D}{\mathbf u}_{m}|^{2})^{\frac{p-2}{4}}{\mathbb D}{\mathbf u}_{m}\|_{{\bf L}^{2}(\Omega)}
    \|(1+|{\mathbb D}{\mathbf u}_{m}|^{2})^{\frac{p-2}{4}}\nabla {\mathbb D}{\mathbf u}_{m}\|_{{\bf L}^{2}(\Omega)}\notag\\
&+C\|\sqrt{\rho_{m}}\partial_t{\mathbf u}_{m}\|_{{\bf L}^{2}(\Omega)}\|\Delta  {\mathbf u}_{m}\|_{{\bf L}^{2}(\Omega)}\notag
     +\|\Delta  {\mathbf u}_{m}\|_{{\bf L}^{2}(\Omega)}\|{\mathbf u}_{m}\|_{{\bf L}^{\frac{2p}{p-2}}(\Omega)}
     \|\nabla{\mathbf u}_{m}\|_{{\bf L}^{p}(\Omega)}\notag\\
&+C\|\Delta  {\mathbf u}_{m}\|_{{\bf L}^{2}(\Omega)}\|\nabla\phi_{m}\|_{{\bf L}^{4}(\Omega)}\|\mu_{m}\|_{L^{4}(\Omega)}
    +C\|\Delta  {\mathbf u}_{m}\|_{{\bf L}^{2}(\Omega)}\|\phi_{m}^{3}-\phi_{m}\|_{L^{6}(\Omega)}\|\nabla\phi_{m}\|_{{\bf L}^{3}(\Omega)}\notag\\
&\leqslant C(\epsilon)\|(1+|{\mathbb D}{\mathbf u}_{m}|^{2})^{\frac{p-2}{4}}{\mathbb D}{\mathbf u}_{m}\|^{2}_{{\bf L}^{2}(\Omega)}
  +\frac{\epsilon}{5}\|(1+|{\mathbb D}{\mathbf u}_{m}|^{2})^{\frac{p-2}{4}}\nabla {\mathbb D}{\mathbf u}_{m}\|^{2}_{{\bf L}^{2}(\Omega)}\notag\\
&+\frac{1}{\epsilon}\|\sqrt{\rho_{m}}\partial_t{\mathbf u}_{m}\|^{2}_{{\bf L}^{2}(\Omega)}
   +\frac{\epsilon}{5}\|\Delta{\mathbf u}_{m}\|^{2}_{{\bf L}^{2}(\Omega)}
   + C(\epsilon)\|(1+|{\mathbb D}{\mathbf u}_{m}|^{2})^{\frac{p-2}{4}}{\mathbb D}{\mathbf u}_{m}\|^{\frac{8}{p}}_{L^{2}(\Omega)}
   +\frac{\epsilon}{5}\|\Delta{\mathbf u}_{m}\|^{2}_{{\bf L}^{2}(\Omega)}\notag\\
   &+C(\epsilon)\|\mu_{m}\|^{2}_{H^{1}(\Omega)}\|\phi_{m}\|^{2}_{H^{2}(\Omega)}
    +\frac{\epsilon}{5}\|\Delta{\mathbf u}_{m}\|^{2}_{{\bf L}^{2}(\Omega)}
    +C(\epsilon)(\|\phi_{m}\|^{6}_{H^{2}(\Omega)}+\|\phi_{m}\|^{2}_{H^{1}(\Omega)})
    +\frac{\epsilon}{5}\|\Delta{\mathbf u}_{m}\|^{2}_{{\bf L}^{2}(\Omega)}\notag\\
&\leqslant \epsilon\|(1+|{\mathbb D}{\mathbf u}_{m}|^{2})^{\frac{p-2}{4}}\nabla {\mathbb D}{\mathbf u}_{m}\|^{2}_{{\bf L}^{2}(\Omega)}
      +C(\epsilon)\|(1+|{\mathbb D}{\mathbf u}_{m}|^{2})^{\frac{p-2}{4}}{\mathbb D}{\mathbf u}_{m}\|^{2}_{{\bf L}^{2}(\Omega)}\notag\\
&+\frac{1}{\epsilon}\|\sqrt{\rho_{m}}\partial_t{\mathbf u}_{m}\|^{2}_{{\bf L}^{2}(\Omega)}
    +C(\epsilon)\left(\|(1+|{\mathbb D}{\mathbf u}_{m}|^{2})^{\frac{p-2}{4}}{\mathbb D}{\mathbf u}_{m}\|^{\frac{8}{p}}_{{\bf L}^{2}(\Omega)}
    +\|\mu_{m}\|^{4}_{H^{1}(\Omega)} +\|\phi_{m}\|^{6}_{H^{2}(\Omega)}\right).
\end{align*}
Taking suitable small $\epsilon\in(0,1),$ one get \eqref{E4-10} from above inequality.
\end{proof}

\begin{lemma}\label{lemmma4-4}
 Under the condition of Theorem \ref{thm2-2}, there exist a positive constant $C$ and $T_0\in (0,T)$, independent of $m$ and $\delta$, such that
\begin{equation}\label{E4-20}
			\|{\mathbf u}_{m}\|_{{\bf L}^{\infty}(0,T_{0};\mathbf{W}_{0}^{1,p}(\Omega))}\leqslant C,\;\|{\mathbf u}_{m}\|_{{\bf L}^{2}(0,T_{0};{\bf H}^{2}(\Omega))}\leqslant C\;and\;\|\sqrt{\rho}\partial_t{\mathbf u}_{m}\|_{{\bf L}^{2}(0,T_{0};{\bf L}^{2}(\Omega))}\leqslant C.
\end{equation}
\end{lemma}

\begin{proof}
One deduces from Lemma \ref{lemmma4-1} and Lemma \ref{lemmma4-2} that
 \begin{align}\label{E4-15}
 	&\frac{d}{dt}\left\{\frac{1}{2}\|\nabla\mu_{m}\|^{2}_{{{\bf L}^{2}(\Omega)}}+\int_{\Omega}\left(\rho_m \mu_{m} {\mathbf u}_{m}\cdot\nabla\phi_{m}
 	+\frac{1}{p}\nu(\phi_{m})(1+|{\mathbb D}{\mathbf u}_{m}|^{2})^{\frac{p}{2}}\right)dx\right\}\notag\\
 	&+\int_{\Omega}\left(|\nabla\partial_t\phi_{m}|^{2}+3\rho_m \phi_m^2|\partial_t\phi_{m}|^{2}+\rho_m |\partial_t\phi_{m}|^{2}+\rho_m  |\partial_t{\mathbf u}_{m}|^{2}\right)dx\notag\\
 	&\leqslant C\left(\|(1+|{\mathbb D}{\mathbf u}_{m}|^{2})^{\frac{p}{4}}\|^{\frac{2(5-p)}{3-p}}_{{\bf L}^{2}(\Omega)}
 	+\|(1+|{\mathbb D}{\mathbf u}_{m}|^{2})^{\frac{p}{4}}\|^{\frac{8}{p}}_{{\bf L}^{2}(\Omega)}+\|\phi_{m}\|^{8}_{H^{2}(\Omega)}+\|\mu_{m}\|^{4}_{H^{1}(\Omega)}\right),
 \end{align}
by taking suitable small $\epsilon\in(0,1).$ For the sake of brevity, introduce $H_{m}(t)$ and $G_m(t)$ are introduced as follows
 \begin{equation}\label{E4-11}
 	H_{m}(t)=\frac{1}{2}\|\nabla\mu_{m}\|^{2}_{{{\bf L}^{2}(\Omega)}}+\int_{\Omega}\left(\rho_m \mu_{m} {\mathbf u}_{m}\cdot\nabla\phi_{m}
 	+\frac{1}{p}\nu(\phi_{m})(1+|{\mathbb D}{\mathbf u}_{m}|^{2})^{\frac{p}{2}}\right)dx
 \end{equation}
 and
 \begin{equation}\label{E4-12}
 	G_m(t)=\int_{\Omega}\left(|\nabla\partial_t\phi_{m}|^{2}+3\rho_m \phi_m^2|\partial_t\phi_{m}|^{2}+\rho_m |\partial_t\phi_{m}|^{2}
 	+\rho_m  |\partial_t{\mathbf u}_{m}|^{2}\right)dx
 \end{equation}
 respectively. Note that
 \begin{align*}
 |\int_{\Omega}\rho_{m}\mu_{m} {\mathbf u}_{m}\cdot\nabla\phi_{m} dx|
 &\leqslant \|\rho_{m}\|^{\frac{3}{4}}_{L^{\infty}(\Omega)}\|\sqrt{\rho_{m}}{\mathbf u}_{m}\|^{\frac{1}{2}}_{L^{2}(\Omega)}
    \|\mu_{m}\|_{L^{6}(\Omega)}\|{\mathbf u}_{m}\|^{\frac{1}{2}}_{{\bf L}^{6}(\Omega)}\|\nabla\phi_{m}\|_{{\bf L}^{2}(\Omega)}\\
 &\leqslant C\|\mu_{m}\|_{H^{1}(\Omega)}\|\nabla {\mathbf u}_{m}\|^{\frac{1}{2}}_{L^{2}(\Omega)}
    \|\sqrt{\rho_{m}} {\mathbf u}_{m}\|^{\frac{1}{2}}_{L^{2}(\Omega)}\|\nabla\phi_{m}\|_{{\bf L}^{2}(\Omega)}\\
 &\leqslant C(1+\|\nabla\mu_{m}\|_{{\bf L}^{2}(\Omega)})\|\nabla {\mathbf u}_{m}\|^{\frac{1}{2}}_{{\bf L}^{p}(\Omega)}\\
 &\leqslant C(1+\|\nabla\mu_{m}\|_{{\bf L}^{2}(\Omega)})
       \|(1+|{\mathbb D}{\mathbf u}_{m}|^{2})^{\frac{p}{4}}\|^{\frac{1}{p}}_{{ L}^{2}(\Omega)}\\
 &\leqslant C(1+\|\nabla\mu_{m}\|_{{\bf L}^{2}(\Omega)})\frac{1}{(\nu_{*})^{\frac{1}{2p}}}
      \|\sqrt{\nu(\phi_{m})}(1+|{\mathbb D}{\mathbf u}_{m}|^{2})^{\frac{p}{4}}\|^{\frac{1}{p}}_{{ L}^{2}(\Omega)}\\
 &\leqslant C+\frac{1}{4}\|\nabla\mu_{m}\|^{2}_{{\bf L}^{2}(\Omega)}+\frac{1}{2p}
 \|\sqrt{\nu(\phi_{m})}(1+|{\mathbb D}{\mathbf u}_{m}|^{2})^{\frac{p}{4}}\|^{2}_{{ L}^{2}(\Omega)}.
 \end{align*}
 One finds that
 \begin{align*}
 	H_{m}(t)&=\frac{1}{2}\|\nabla\mu_{m}\|^{2}_{{{\bf L}^{2}(\Omega)}}+\int_{\Omega}\left(\rho_m \mu_{m} {\mathbf u}_{m}\cdot\nabla\phi_{m}
 	+\frac{1}{p}\nu(\phi_{m})(1+|{\mathbb D}{\mathbf u}_{m}|^{2})^{\frac{p}{2}}\right)dx\\
 	&\geqslant\frac{1}{2}\|\nabla\mu_{m}\|^{2}_{{{\bf L}^{2}(\Omega)}}+\frac{1}{p}\int_{\Omega}\nu(\phi_{m})(1+|{\mathbb D}{\mathbf u}_{m}|^{2})^{\frac{p}{2}}dx-|\int_{\Omega}\rho_m \mu_{m} {\mathbf u}_{m}\cdot\nabla\phi_{m}dx|\\
 	&\geqslant \frac{1}{4}\|\nabla\mu_{m}\|^{2}_{{{\bf L}^{2}(\Omega)}}+\frac{1}{2p}\int_{\Omega}\nu(\phi_{m})(1+|{\mathbb D}{\mathbf u}_{m}|^{2})^{\frac{p}{2}}dx-C(E_{0},\nu_{*})
 \end{align*}
 and
 \begin{align*}
 	H_{m}(t)
 	\leqslant \frac{3}{4}\|\nabla\mu_{m}\|^{2}_{{{\bf L}^{2}(\Omega)}}+\frac{3}{2p}\int_{\Omega}\nu(\phi_{m})(1+|{\mathbb D}{\mathbf u}_{m}|^{2})^{\frac{p}{2}}dx+C(E_{0},\nu^{*}).
 \end{align*}
 Setting $c=\frac{5-p}{3-p},$ one can find $c>1$ and deduce from \eqref{E4-15} that
 \begin{equation}\label{E4-13}
 	\frac{d}{dt}H_{m}+G_m\leqslant C(E_{0})(C(E_{0})+H_{m})^{c}.
 \end{equation}
 Now, one needs to bound $H_m(0)$ by a positive constant $C,$ which is independent of $m$ and $\delta$. First, one evaluates $\eqref{E3-3}_4$ at $t=0$ and takes $w=\mu_m(0)$ to arrive at
 \begin{align}\label{chuE4-3}
 	&|\int_\Omega\rho_{0\delta}\mu_{m}(0)^2dx|\notag\\
 	&=\int_\Omega\Delta\phi_{m}(0)\mu_{m}(0)dx+\int_\Omega\rho_{0\delta}(\phi_{m}(0)^2-1)\phi_{m}(0)dx\notag\\
 	&\leqslant \|\nabla\phi_{m}(0) \|_{{\bf L}^2(\Omega)}\|\nabla\mu_{m}(0) \|_{{\bf L}^2(\Omega)}
 	+\|\rho_{0\delta}\|_{L^\infty(\Omega)}\left(\|\phi_{m}(0)\|_{L^3(\Omega)}^3+\|\phi_{m}(0)\|_{L^1(\Omega)}\right)\notag\\
 	&\leqslant \frac12\|\nabla\phi_{m}(0) \|_{{\bf L}^2(\Omega)}^2+\frac12\|\nabla\mu_{m}(0) \|_{{\bf L}^2(\Omega)}^2
 	+\|\rho_{0\delta}\|_{L^\infty(\Omega)}\left(\|\phi_{m}(0)\|_{L^6(\Omega)}^3+\|\phi_{m}(0)\|_{L^1(\Omega)}\right).
 \end{align}
 It is deduced from the initial condition in Theorem \ref{thm2-2} that $$\|\sqrt{\rho_{0\delta}}\mu_{m}(0)\|_{L^{2}(\Omega)}\leqslant C.$$ Furthermore, $$\|\mu_{m}(0)\|_{H^{1}(\Omega)}\leqslant C.$$
 Second, one evaluates $\eqref{E3-3}_4$ at $t=0$ and takes $w=\Delta\phi_{0m}$ to arrive at
 \begin{align*}
 	\|\Delta\phi_{m}(0)\|^{2}_{L^2(\Omega)}&\leqslant\|\rho_{0\delta}\|_{L^\infty(\Omega)}^\frac{1}{2}\|\sqrt{\rho_{0\delta}}\mu_m(0)\|_{L^2(\Omega)}\|\Delta\phi_{m}(0)\|_{L^2(\Omega)}\\
 	&+\|\rho_{0\delta}\|_{L^\infty(\Omega)}\|\Delta\phi_{m}(0)\|_{L^2(\Omega)}\left(\|\phi_{m}(0)\|^{3}_{L^6(\Omega)}+\|\phi_{m}(0)\|_{L^2(\Omega)}\right).
 \end{align*}
 Similarly, one finds that $$\|\Delta\phi_{0m}\|_{L^{2}(\Omega)}\leqslant C$$ and so $\|\phi_{0m}\|_{H^{2}(\Omega)}\leqslant C$, where $C$ is independent of $m$ and $\delta$. Thus,
 \begin{align}\label{E4-17}
 	&H_{m}(0)
 	=\frac{1}{2}\|\nabla\mu_{m}(0)\|^{2}_{{{\bf L}^{2}(\Omega)}}+\int_{\Omega}\left(\rho_{0\delta} \mu_{m}(0) {\mathbf u}_{m}(0)\cdot\nabla\phi_{m}(0)
 	+\frac{1}{p}\nu(\phi_{m}(0))(1+|{\mathbb D}{\mathbf u}_{m}(0)|^{2})^{\frac{p}{2}}\right)dx\notag\\
 	&\leqslant \frac{1}{2}\|\nabla\mu_{m}(0)\|^{2}_{{{\bf L}^{2}(\Omega)}}
 	+C\left(\|\sqrt{\rho_{0\delta}}\mu_{m}(0)\|_{{L^{2}(\Omega)}}\|{\mathbf u}_{m}(0)\|_{\mathbf{W}_{0}^{1,p}(\Omega)}\|\nabla\phi_{m}(0)\|_{{\bf L}^6(\Omega)}+\|{\mathbb D}{\mathbf u}_{m}(0)\|_{{\bf L}^{p}(\Omega)}^p+1\right)\notag\\
 	&\leqslant C,
 \end{align}
where $C$ is independent of $m$ and $\delta$.

Obviously, it is deduced from Lemma \ref{gronwall1} to \eqref{E4-13} and find that there exists ${T_0}\in (0,T),$ independent of $m$ and $\delta,$ such that
\begin{equation*}
	[C(E_{0})+H_{m}(0)]^{1-c}+\frac{{T_0}C(E_{0})}{c-1}>0
\end{equation*}
 and
 \begin{equation}\label{E4-18}
 	H_{m}(t)\leqslant C\;\;(\forall t\in[0,{T_0}]),
 \end{equation}
 where $C$ is independent of $m$ and $\delta$. Then
 \begin{equation}\label{E4-19}
 	\int_0^{T_0}\int_{\Omega}\left(|\nabla\partial_t\phi_{m}|^{2}+3\rho_m \phi_m^2|\partial_t\phi_{m}|^{2}+\rho_m |\partial_t\phi_{m}|^{2}+\rho_m  |\partial_t{\mathbf u}_{m}|^{2}\right)dxdt
 	\leqslant C.
 \end{equation}
Therefore, there exist $T_{0}>0$ and $C>0,$ depending only on the initial date, such that
 \begin{equation*}
 	\|\nabla{\mathbf u}_{m}\|_{{\bf L}^{\infty}(0,T_{0};{\bf L}^{p}(\Omega))}\leqslant C\;\mbox{ and }\;\|\sqrt{\rho}\partial_t{\mathbf u}_{m}\|_{{\bf L}^{2}(0,T_{0};{\bf L}^{2}(\Omega))}\leqslant C.
 \end{equation*}
 Note that
 \begin{equation*}
 	\|\nabla^{2}{\mathbf u}_{m}\|_{{\bf L}^{2}(0,T_{0};{\bf L}^{2}(\Omega))}\leqslant\|\sqrt{\nu(\phi_{m})}(1+|{\mathbb D}{\mathbf u}_{m}|^{2})^{\frac{p-2}{4}}
 	|\nabla {\mathbb D}{\mathbf u}_{m}|\|_{{\bf L}^{2}(0,T_{0};{\bf L}^{2}(\Omega))}\leqslant C,
 \end{equation*}
follows from Lemma \ref{lemmma4-2}. It is deduced from Lemma \ref{lem2-3} that
 \begin{equation*}
\|{\mathbf u}_{m}\|_{{\bf L}^{\infty}(0,T_{0};\mathbf{W}_{0}^{1,p}(\Omega))}\leqslant C\;\mbox{ and }\;\|{\mathbf u}_{m}\|_{{\bf L}^{2}(0,T_{0};{\bf H}^{2}(\Omega))}\leqslant C.
\end{equation*}
\end{proof}

Next, one improves the regularity of $\nabla\partial_{t}\phi_{m}.$

\begin{lemma}\label{lemmma4-3}
	Under the condition of Theorem \ref{thm2-2}, there exists a positive constant $C$, independent of $m$ and $\delta$, such that
\begin{equation}\label{E4-25}
	\|\partial_{t}\phi_{m}\|_{L^{\frac{6-p}{3-p}}(0,T_{0};H^{1}(\Omega))}\leqslant C,
\end{equation}
where $T_0$ is obtained in Lemma \ref{lemmma4-4}.
\end{lemma}

\begin{proof}
Recall that the equation \eqref{(E3-35)}, one obtains that
\begin{align}\label{E4-27}
  &\frac{d}{dt}\int_{\Omega}\frac{1}{2}|\nabla\mu_{m}|^{2}dx
  +\int_{\Omega}|\nabla\partial_t\phi_{m}|^{2}dx+\int_{\Omega}3\rho_{m}\phi_{m}^{2}(\partial_t\phi_{m})^{2}dx\notag\\
   &+3\int_{\Omega}\rho_{m}\phi_{m}^{2}|{\mathbf u}_{m}|^{2}|\nabla\phi_{m}|^{2}dx+\int_{\Omega}\rho_{m}(\partial_t\phi_{m})^{2}dx\notag\\
   &=\underset{I_{1}}{\underbrace{\int_{\Omega}\rho_{m} {\mathbf u}_{m}\cdot\nabla\mu_{m}\partial_t\phi_{m}dx}}
  +\underset{I_{2}}{\underbrace{\int_{\Omega}\rho_{m}\mu_{m}{\mathbf u}_{m}\cdot\nabla\partial_t\phi_{m}dx}}
  +\underset{I_{3}}{\underbrace{\int_{\Omega}\rho_{m}\phi_{m} {\mathbf u}_{m}\cdot\nabla\partial_t\phi_{m}dx}}\notag\\
   &-\underset{I_{4}}{\underbrace{\int_{\Omega}\rho_{m}\phi^{3}_{m}{\mathbf u}_{m}\cdot\nabla\partial_t\phi_{m}dx}}
   -\underset{I_{5}}{\underbrace{6\int_{\Omega}\rho_{m}\phi^{2}_{m}{\mathbf u}_{m}\cdot\nabla\phi_{m}\partial_t\phi_{m}dx}}
   -\underset{I_{6}}{\underbrace{\int_{\Omega}\nabla\partial_t\phi_{m}\cdot(\nabla {\mathbf u}_{m}\cdot\nabla\phi_{m})dx}}\notag\\
	&-\underset{I_{7}}{\underbrace{\int_{\Omega}{\mathbf u}_{m}\cdot(\nabla^{2} \phi_{m}\cdot\nabla\partial_t\phi_{m})dx}}
   +\underset{I_{8}}{\underbrace{\int_{\Omega}\rho_{m} {\mathbf u}_{m}\cdot\nabla\mu_{m}({\mathbf u}_{m}\cdot\nabla\phi_{m})dx}}\notag\\
	&+\underset{I_{9}}{\underbrace{\int_{\Omega}\rho_{m}\mu_{m} {\mathbf u}_{m}\cdot(\nabla {\mathbf u}_{m}\cdot\nabla\phi_{m})dx}}
    +\underset{I_{10}}{\underbrace{\int_{\Omega}\rho_{m}\mu_{m} {\mathbf u}_{m}\cdot({\mathbf u}_{m}\cdot\nabla^{2}\phi_{m})dx}}
   -\underset{I_{11}}{\underbrace{\int_{\Omega}\rho_{m}\phi^{3}_{m}{\mathbf u}_{m}\cdot({\mathbf u}_{m}\cdot\nabla^{2}\phi_{m})dx}}\notag\\
   &-\underset{I_{12}}{\underbrace{\int_{\Omega}\rho_{m}\phi^{3}_{m}{\mathbf u}_{m}\cdot(\nabla {\mathbf u}_{m}\cdot\nabla\phi_{m})dx}}+\underset{I_{13}}{\underbrace{\int_{\Omega}\rho_{m}|{\mathbf u}_{m}|^{2}|\nabla\phi_{m}|^{2}dx}}
  +\underset{I_{14}}{\underbrace{\int_{\Omega}\rho_{m}\phi_{m} {\mathbf u}_{m}\cdot(\nabla {\mathbf u}_{m}\cdot\nabla\phi_{m})dx}}\notag\\
	&  +\underset{I_{15}}{\underbrace{\int_{\Omega}\rho_{m}\phi_{m} {\mathbf u}_{m}\cdot( {\mathbf u}_{m}\cdot\nabla^{2}\phi_{m})dx}}.
\end{align}
The estimate of $I_{1}-I_{15}$ is given from a new perspectives as follows
\begin{align*}
	&|I_{1}|=|\int_{\Omega}\rho_{m} {\mathbf u}_{m}\cdot\nabla\mu_{m}\partial_t\phi_{m}dx|\\
	&\;\;\;\;\leqslant C\|\partial_t\phi_{m}\|_{L^{6}(\Omega)}\|\nabla\mu_{m}\|_{{\bf L}^{2}(\Omega))}\|{\mathbf u}_{m}\|_{{\bf L}^{3}(\Omega))}\\
	&\;\;\;\;\leqslant C\|\nabla\partial_t\phi_{m}\|_{{\bf L}^{2}(\Omega)}\|\mu_{m}\|_{H^{1}(\Omega)}\|{\mathbf u}_{m}
        \|_{{\bf W}_{0}^{1,p}(\Omega)},\\
    &\;\;\;\;\leqslant \frac{1}{14}\|\nabla\partial_t\phi_{m}\|_{{\bf L}^{2}(\Omega)}^2+C\|\mu_{m}\|_{H^{1}(\Omega)}^2\|{\mathbf u}_{m}
    \|_{{\bf W}_{0}^{1,p}(\Omega)}^2,\\
	&|I_{2}|=|\int_{\Omega}\rho_{m}\mu_{m}{\mathbf u}_{m}\cdot\nabla\partial_t\phi_{m}dx|\\
	&\;\;\;\;\leqslant C\|\nabla\partial_t\phi_{m}\|_{{\bf L}^{2}(\Omega)}\|{\mathbf u}_{m}\|_{{\bf L}^{3}(\Omega))}
        \|\mu_{m}\|_{L^{6}(\Omega)}\|\rho_{m}\|_{L^{\infty}(\Omega)}\\
	&\;\;\;\;\leqslant C\|\nabla\partial_t\phi_{m}\|_{{\bf L}^{2}(\Omega)}
       \|{\mathbf u}_{m}\|_{{\bf W}_{0}^{1,p}(\Omega)}\|\mu_{m}\|_{H^{1}(\Omega)}\\
   &\;\;\;\;\leqslant \frac{1}{14}\|\nabla\partial_t\phi_{m}\|_{{\bf L}^{2}(\Omega)}^2+C\|\mu_{m}\|_{H^{1}(\Omega)}^2\|{\mathbf u}_{m}
    \|_{{\bf W}_{0}^{1,p}(\Omega)}^2,\\
	&|I_{3}|=|\int_{\Omega}\rho_{m}\phi_{m} {\mathbf u}_{m}\cdot\nabla\partial_t\phi_{m}dx|\\
	&\;\;\;\;\leqslant C\|\nabla\partial_t\phi_{m}\|_{{\bf L}^{2}(\Omega)}\|{\mathbf u}_{m}\|_{{\bf L}^{3}(\Omega)}\|\phi_{m}\|_{L^{6}(\Omega)}\\
    &\;\;\;\;\leqslant C\|\nabla\partial_t\phi_{m}\|_{{\bf L}^{2}(\Omega)}\|{\mathbf u}_{m}\|_{{\bf W}_{0}^{1,p}(\Omega)}
    \|\phi_{m}\|_{H^{1}(\Omega)}\\
   &\;\;\;\;\leqslant \frac{1}{14}\|\nabla\partial_t\phi_{m}\|_{{\bf L}^{2}(\Omega)}^2+C\|\phi_{m}\|_{H^{1}(\Omega)}^2\|{\mathbf u}_{m}
    \|_{{\bf W}_{0}^{1,p}(\Omega)}^2,\\
	&|I_{4}|=|\int_{\Omega}\rho_{m}\phi^{3}_{m}{\mathbf u}_{m}\cdot\nabla\partial_t\phi_{m}dx|\\
	&\;\;\;\;\leqslant C\|\phi^{3}_{m}\|_{L^{3}(\Omega)}\|{\mathbf u}_{m}\|_{{\bf L}^{6}(\Omega)}\|\nabla\partial_t\phi_{m}\|_{{\bf L}^{2}(\Omega)}\\
	&\;\;\;\;\leqslant C\|\nabla\partial_t\phi_{m}\|_{{\bf L}^{2}(\Omega)}\|{\mathbf u}_{m}\|_{{\bf W}_{0}^{1,p}(\Omega)}\|\phi_{m}\|^{3}_{H^{2}(\Omega)}\\
   &\;\;\;\;\leqslant \frac{1}{14}\|\nabla\partial_t\phi_{m}\|_{{\bf L}^{2}(\Omega)}^2
   +C\|{\mathbf u}_{m}\|_{{\bf W}_{0}^{1,p}(\Omega)}^2\|\phi_{m}\|^{6}_{H^{2}(\Omega)},\\
	&|I_{5}|=|6\int_{\Omega}\rho_{m}\phi^{2}_{m}{\mathbf u}_{m}\cdot\nabla\phi_{m}\partial_t\phi_{m}dx|\\
	&\;\;\;\;\leqslant C\|\nabla\phi_{m}\|_{ {\bf L}^{6}(\Omega)}\|\phi^{2}_{m}\|_{L^{6}(\Omega)}\|\partial_t\phi_{m}\|_{L^{6}(\Omega)}\|{\mathbf u}_{m}\|_{\bf{L}^{6}(\Omega)}\\
	&\;\;\;\;\leqslant C\|\nabla\partial_t\phi_{m}\|_{{\bf L}^{2}(\Omega)}\|{\mathbf u}_{m}\|_{{\bf W}_{0}^{1,p}(\Omega)}\|\phi_{m}\|^{3}_{H^{2}(\Omega)}\\
   &\;\;\;\;\leqslant \frac{1}{14}\|\nabla\partial_t\phi_{m}\|_{{\bf L}^{2}(\Omega)}^2
   +C\|{\mathbf u}_{m}\|_{{\bf W}_{0}^{1,p}(\Omega)}^2\|\phi_{m}\|^{6}_{H^{2}(\Omega)},\\
	&|I_{6}|=|\int_{\Omega}\nabla\partial_t\phi_{m}\cdot(\nabla {\mathbf u}_{m}\cdot\nabla\phi_{m})dx|\\
	&\;\;\;\;\leqslant C\|\nabla\partial_t\phi_{m}\|_{{\bf L}^{2}(\Omega)}\|\nabla {\mathbf u}_{m}\|_{{\bf L}^{3}(\Omega)}\|\nabla\phi_{m}\|_{{\bf L}^{6}(\Omega)}\\
	&\;\;\;\;\leqslant C\|\nabla\partial_t\phi_{m}\|_{{\bf L}^{2}(\Omega)}\|\nabla {\mathbf u}_{m}\|_{{\bf L}^{3}(\Omega)}\|\phi_{m}\|_{H^{2}(\Omega)}\\
	&\;\;\;\;\leqslant C\|\nabla\partial_t\phi_{m}\|_{{\bf L}^{2}(\Omega)}\|{\mathbf u}_{m}\|^{\frac{p}{6-p}}_{{\bf W}_{0}^{1,p}(\Omega)}\|{\mathbf u}_{m}\|^{\frac{6-2p}{6-p}}_{{\bf W}_{0}^{2,2}(\Omega)}\|\phi_{m}\|_{H^{2}(\Omega)}\\
   &\;\;\;\;\leqslant \frac{1}{14}\|\nabla\partial_t\phi_{m}\|_{{\bf L}^{2}(\Omega)}^2
   +C\|{\mathbf u}_{m}\|^{\frac{2p}{6-p}}_{{\bf W}_{0}^{1,p}(\Omega)}
    \|{\mathbf u}_{m}\|^{\frac{4(3-p)}{6-p}}_{{\bf W}_{0}^{2,2}(\Omega)}\|\phi_{m}\|_{H^{2}(\Omega)}^2,\\
	&|I_{7}|=|\int_{\Omega}{\mathbf u}_{m}\cdot(\nabla^{2} \phi_{m}\cdot\nabla\partial_t\phi_{m})dx|\\
	&\;\;\;\;\leqslant C\|\nabla\partial_t\phi_{m}\|_{{\bf L}^{2}(\Omega)}\| {\mathbf u}_{m}\|_{{\bf L}^{\frac{3p}{3-p}}(\Omega)}\|\nabla^{2}\phi_{m}\|_{{\bf L}^{\frac{6p}{5p-6}}(\Omega)}\\
	&\;\;\;\;\leqslant C\|\nabla\partial_t\phi_{m}\|_{{\bf L}^{2}(\Omega)}\| {\mathbf u}_{m}\|_{{\bf W}_{0}^{1,p}(\Omega)}\|\phi_{m}\|^{\frac{5p-6}{3p}}_{H^{2}(\Omega)}\|\phi_{m}\|^{\frac{6-2p}{3p}}_{W^{2,\infty}(\Omega)}\\
   &\;\;\;\;\leqslant \frac{1}{14}\|\nabla\partial_t\phi_{m}\|_{{\bf L}^{2}(\Omega)}^2
   +\| {\mathbf u}_{m}\|_{{\bf W}_{0}^{1,p}(\Omega)}^2\|\phi_{m}\|^{\frac{2(5p-6)}{3p}}_{H^{2}(\Omega)}\|\phi_{m}\|^{\frac{4(3-p)}{3p}}_{W^{2,\infty}(\Omega)},\\
	&|I_{8}|=|\int_{\Omega}\rho_{m} {\mathbf u}_{m}\cdot\nabla\mu_{m}({\mathbf u}_{m}\cdot\nabla\phi_{m})dx|\\
	&\;\;\;\;\leqslant C\|{\mathbf u}_{m}\|^{2}_{{\bf L}^{6}(\Omega)}\|\nabla\mu_{m}\|_{{\bf L}^{2}(\Omega)}\|\nabla\phi_{m}\|_{{\bf L}^{6}(\Omega)}\\
	&\;\;\;\;\leqslant C\|{\mathbf u}_{m}\|^{2}_{{\bf W}_{0}^{1,p}(\Omega)}\|\mu_{m}\|_{H^{1}(\Omega)}\|\phi_{m}\|_{H^{2}(\Omega)},\\
	&|I_{9}|=|\int_{\Omega}\rho_{m}\mu_{m} {\mathbf u}_{m}\cdot(\nabla {\mathbf u}_{m}\cdot\nabla\phi_{m})dx|\\
	&\;\;\;\;\leqslant C\|\nabla {\mathbf u}_{m}\|_{{\bf L}^{2}(\Omega)}\|{\mathbf u}_{m}\|_{{\bf L}^{6}(\Omega)}\|\mu_{m}\|_{L^{6}(\Omega)}\|\nabla\phi_{m}\|_{{\bf L}^{6}(\Omega)}\\
	&\;\;\;\;\leqslant C\|{\mathbf u}_{m}\|^{2}_{{\bf W}_{0}^{1,p}(\Omega)}\|\mu_{m}\|_{H^{1}(\Omega)}\|\phi_{m}\|_{H^{2}(\Omega)},\\
	&|I_{10}|=|\int_{\Omega}\rho_{m}\mu_{m} {\mathbf u}_{m}\cdot({\mathbf u}_{m}\cdot\nabla^{2}\phi_{m})dx|\\
	&\;\;\;\;\leqslant\|{\mathbf u}_{m}\|^{2}_{{\bf L}^{6}(\Omega)}\|\mu_{m}\|_{L^{6}(\Omega)}\|\nabla^{2}\phi_{m}\|_{{\bf L}^{2}(\Omega)}\\
	&\;\;\;\;\leqslant C\|{\mathbf u}_{m}\|^{2}_{{\bf W}_{0}^{1,p}(\Omega)}\|\mu_{m}\|_{H^{1}(\Omega)}\|\phi_{m}\|_{H^{2}(\Omega)},\\
	&|I_{11}|=|\int_{\Omega}\rho_{m}\phi^{3}_{m}{\mathbf u}_{m}\cdot({\mathbf u}_{m}\cdot\nabla^{2}\phi_{m})dx|\\
	&\;\;\;\;\leqslant C\|{\mathbf u}_{m}\|^{2}_{{\bf L}^{6}(\Omega)}\|\nabla^{2}\phi_{m}\|_{{\bf L}^{2}(\Omega)}\|\phi_m\|^{3}_{L^{3}(\Omega)}\\
	&\;\;\;\;\leqslant C\|{\mathbf u}_{m}\|^{2}_{{\bf W}_{0}^{1,p}(\Omega)}\|\phi_{m}\|^{4}_{H^{2}(\Omega)},\\
	&|I_{12}|=|\int_{\Omega}\rho_{m}\phi^{3}_{m}{\mathbf u}_{m}\cdot(\nabla {\mathbf u}_{m}\cdot\nabla\phi_{m})dx|\\
	&\;\;\;\;\leqslant C\|{\mathbf u}_{m}\|_{{\bf L}^{6}(\Omega)}\|\nabla{\mathbf u}_{m}\|_{{\bf L}^{2}(\Omega)}\|\nabla\phi_{m}\|_{{\bf L}^{6}(\Omega)}\|\phi_{m}^{3}\|_{L^{6}(\Omega)}\\
	&\;\;\;\;\leqslant C\|{\mathbf u}_{m}\|^{2}_{{\bf W}_{0}^{1,p}(\Omega)}\|\phi_{m}\|^{4}_{H^{2}(\Omega)},\\
	&|I_{13}|=|\int_{\Omega}\rho_{m}{\mathbf u}_{m}^{2}\nabla\phi_{m}^{2}dx|\\
	&\;\;\;\;\leqslant C\|{\mathbf u}_{m}\|^{2}_{{\bf L}^{4}(\Omega)}\|\nabla\phi_{m}\|^{2}_{{\bf L}^{4}(\Omega)}\\
	&\;\;\;\;\leqslant C\|{\mathbf u}_{m}\|^{2}_{{\bf W}_{0}^{1,p}(\Omega)}\|\phi_{m}\|^{2}_{H^{2}(\Omega)},\\
	&|I_{14}|=|\int_{\Omega}\rho_{m}\phi_{m} {\mathbf u}_{m}\cdot(\nabla {\mathbf u}_{m}\cdot\nabla\phi_{m})dx|\\
	&\;\;\;\;\leqslant C\|\phi_m\|_{L^{6}(\Omega)}\|{\mathbf u}_{m}\|_{{\bf L}^{6}(\Omega)}\|\nabla{\mathbf u}_{m}\|_{{\bf L}^{2}(\Omega)}\|\nabla\phi_{m}\|_{{\bf L}^{6}(\Omega)}\\
	&\;\;\;\;\leqslant C\|{\mathbf u}_{m}\|^{2}_{{\bf W}_{0}^{1,p}(\Omega)}\|\phi_{m}\|^{2}_{H^{2}(\Omega)},\\
	&|I_{15}|=|\int_{\Omega}\rho_{m}\phi_{m} {\mathbf u}_{m}\cdot( {\mathbf u}_{m}\cdot\nabla^{2}\phi_{m})dx|\\
	&\;\;\;\;\leqslant C\|{\mathbf u}_{m}\|^{2}_{{\bf L}^{6}(\Omega)}\|\nabla^{2}\phi_{m}\|_{L^{2}(\Omega)}\|\phi_{m}\|_{L^{6}(\Omega)}\\
&\;\;\;\;\leqslant C\|{\mathbf u}_{m}\|^{2}_{{\bf W}_{0}^{1,p}(\Omega)}\|\phi_{m}\|^{2}_{H^{2}(\Omega)}.
\end{align*}
So, it follows from  \eqref{E4-27} that
\begin{align}\label{E4-27}
&\frac{d}{dt}\int_{\Omega}\frac{1}{2}|\nabla\mu_{m}|^{2}dx
  +\int_{\Omega}|\nabla\partial_t\phi_{m}|^{2}dx+\int_{\Omega}\rho_{m}\phi_{m}^{2}(\partial_t\phi_{m})^{2}dx\notag\\
&+\int_{\Omega}\rho_{m}\phi_{m}^{2}|{\mathbf u}_{m}|^{2}|\nabla\phi_{m}|^{2}dx+\int_{\Omega}\rho_{m}(\partial_t\phi_{m})^{2}dx\notag\\
&\leqslant C\|{\mathbf u}_{m}\|_{{\bf W}_{0}^{1,p}(\Omega)}^2\left(\|\mu_{m}\|_{H^{1}(\Omega)}^2+\|\phi_{m}\|_{H^{2}(\Omega)}^2
   +\|\phi_{m}\|^{6}_{H^{2}(\Omega)}+\|\phi_{m}\|^{\frac{2(5p-6)}{3p}}_{H^{2}(\Omega)}\|\phi_{m}\|^{\frac{4(3-p)}{3p}}_{W^{2,\infty}(\Omega)}\right)
      \notag\\
 &~~~~+C\|{\mathbf u}_{m}\|^{\frac{2p}{6-p}}_{{\bf W}_{0}^{1,p}(\Omega)}
 \|{\mathbf u}_{m}\|^{\frac{4(3-p)}{6-p}}_{{\bf W}_{0}^{2,2}(\Omega)}\|\phi_{m}\|_{H^{2}(\Omega)}^2\notag\\
&\leqslant C\|{\mathbf u}_{m}\|_{{\bf W}_{0}^{1,p}(\Omega)}^2\left(1+\|\phi_{m}\|^{\frac{4(3-p)}{3p}}_{W^{2,\infty}(\Omega)}\right)
      +C\|{\mathbf u}_{m}\|^{\frac{4(3-p)}{6-p}}_{{\bf W}_{0}^{1,p}(\Omega)}
   \|{\mathbf u}_{m}\|^{\frac{2p}{6-p}}_{{\bf W}_{0}^{2,2}(\Omega)}\notag.
\end{align}
With the help of Gronwall's inequality and Lemma \ref{L-2-3}, one gets that
\begin{equation*}
\|\nabla\mu_{m}\|_{{\bf L}^\infty(0,T_0;{\bf L}^2(\Omega))}\leqslant C,\;\;\|\partial_t\phi_{m}\|_{{ L}^2(0,T_0;{H}^1(\Omega))}\leqslant C.
\end{equation*}
Furthermore, one obtains that
\begin{align*}
	\|\partial_{t}\phi_{m}\|^{\frac{6-p}{3-p}}_{H^{1}(\Omega)}&\leqslant  C\|{\mathbf u}_{m}\|^{\frac{6-p}{3-p}}_{{\bf W}_{0}^{1,p}(\Omega)}\left(1+\|\phi_{m}\|^{\frac{2(6-p)}{3p}}_{W^{2,\infty}(\Omega)}\right)
	+C\|{\mathbf u}_{m}\|^{\frac{p}{3-p}}_{{\bf W}_{0}^{1,p}(\Omega)}
	\|{\mathbf u}_{m}\|^{2}_{{\bf W}_{0}^{2,2}(\Omega)}\\
	&\leqslant C\|{\mathbf u}_{m}\|^{\frac{6-p}{3-p}}_{{\bf W}_{0}^{1,p}(\Omega)}\left(1+\|\phi_{m}\|^{2}_{W^{2,\infty}(\Omega)}\right)
	+C\|{\mathbf u}_{m}\|^{\frac{p}{3-p}}_{{\bf W}_{0}^{1,p}(\Omega)}
	\|{\mathbf u}_{m}\|^{2}_{{\bf W}_{0}^{2,2}(\Omega)}.
\end{align*}
Thus, it is deduced from Lemma \ref{chuzhi4-1}, \eqref{E4-1} and \eqref{E4-20} that
\begin{equation*}
		\|\partial_{t}\phi_{m}\|_{L^{\frac{6-p}{3-p}}(0,T_{0};H^{1}(\Omega))}\leqslant  C.
\end{equation*}
\end{proof}

\begin{lemma}\label{lemma4-6}
Under the condition of Theorem \ref{thm2-2}, there exists a positive constant $C,$ independent of $m$ and $\delta$, such that
\begin{equation}\label{(E4-28)}
\|\sqrt{\rho}\partial_t{\mathbf u}_{m}\|_{{\bf L}^{\infty}(0,T_{0};{\bf L}^{2}(\Omega))}\leqslant C\;\;\; \|\partial_t{\mathbf u}_{m}\|_{\mathbf{L}^{2}(0,T_{0};\mathbf{W}_{0}^{1,2}(\Omega))}\leqslant C,
\end{equation}
where $T_0$ is obtained in Lemma \ref{lemmma4-4}.
\end{lemma}

\begin{proof}
One takes ${\bf w}=\partial_t{\mathbf u}_{m}$ as the test function for $\eqref{E4-3}_2$ and gets that
\begin{align}\label{E4-22}
  &\frac{1}{2}\frac{d}{dt}\int_{\Omega}\rho_{m}|\partial_t{\mathbf u}_{m}|^{2}dx\notag\\
  &+\int_{\Omega}\left(\nu(\phi_{m})(p-2)(1+|{\mathbb D}{\mathbf u}_{m}|^{2})^{\frac{p-4}{2}}|{\mathbb D}{\mathbf u}_{m}|^{2}|\nabla \partial_t{\mathbf u}_{m}|^{2}
  +\nu(\phi_{m})(1+|{\mathbb D}{\mathbf u}_{m}|^{2})^{\frac{p-2}{2}}|\nabla \partial_t{\mathbf u}_{m}|^{2}\right)dx\notag\\
  &\leqslant\underset{K_{1}}{\underbrace{\int_{\Omega}\rho_{m}|{\mathbf u}_{m}||\partial_t{\mathbf u}_{m}||\nabla {\mathbf u}_{m}|^{2}dx}}
  +\underset{K_{2}}{\underbrace{\int_{\Omega}\rho_{m}|{\mathbf u}_m|^{2}|\nabla^{2}{\mathbf u}_{m}||\partial_t{\mathbf u}_{m}|dx}}\notag\\
  &+\underset{K_{3}}{\underbrace{\int_{\Omega}\rho_{m}|{\mathbf u}_{m}|^{2}|\nabla {\mathbf u}_{m}||\nabla \partial_t{\mathbf u}_{m}|dx}}
  +\underset{K_{4}}{\underbrace{\int_{\Omega}\rho_{m}|\partial_t{\mathbf u}_{m}|^{2}|\nabla {\mathbf u}_{m}|dx}}
  +\underset{K_{5}}{\underbrace{2\int_{\Omega}\rho_{m}|{\mathbf u}_m||\partial_t{\mathbf u}_{m}||\nabla \partial_t{\mathbf u}_{m}|dx}}\notag\\
 &+\underset{K_{6}}{\underbrace{\int_{\Omega}\nu^\prime(\phi_{m})|\partial_t\phi_{m}|(1+|{\mathbb D}{\mathbf u}_{m}|^{2})^{\frac{p-2}{2}}
  |{\mathbb D}{\mathbf u}_{m}||\nabla \partial_t{\mathbf u}_{m}|dx}}
 +\underset{K_{7}}{\underbrace{2\int_{\Omega}|\nabla\phi_{m}||\nabla\partial_t\phi_{m}||\nabla \partial_t{\mathbf u}_{m}|dx}}.
\end{align}
Here $K_{1}-K_{7}$ are estimated as follows
	\begin{align*}
		&K_{1}=\int_{\Omega}\rho_{m}|{\mathbf u}_{m}||\partial_t{\mathbf u}_{m}||\nabla {\mathbf u}_{m}|^{2}dx\\
		&\;\;\;\;\leqslant C\|\nabla {\mathbf u}_{m}\|_{{\bf L}^{3}(\Omega)}^2\|{\mathbf u}_{m}\|_{{\bf L}^{6}(\Omega)}\|\partial_t{\mathbf u}_{m}\|_{{\bf L}^{6}(\Omega)}\\
		&\;\;\;\;\leqslant C\|\nabla^{2}{\mathbf u}_{m}\|_{{\bf L}^{2}(\Omega)}\|\nabla {\mathbf u}_{m}\|^{2}_{{\bf L}^{2}(\Omega)}\|\nabla \partial_t{\mathbf u}_{m}\|_{{\bf L}^{2}(\Omega)}\\
		&\;\;\;\;\leqslant C\|{\mathbf u}_{m}\|_{{\bf H}^{2}(\Omega)}\| {\mathbf u}_{m}\|^{2}_{{\bf W}_{0}^{1,p}(\Omega)}\|\nabla \partial_t{\mathbf u}_{m}\|_{{\bf L}^{2}(\Omega)}\\
		&\;\;\;\;\leqslant \frac{\eta}{7}\|\nabla \partial_t{\mathbf u}_{m}\|^{2}_{{\bf L}^{2}(\Omega)}+C\|{\mathbf u}_{m}\|^{2}_{{\bf H}^{2}(\Omega)}\|{\mathbf u}_{m}\|^{4}_{{\bf W}_{0}^{1,p}(\Omega)},\\
		&K_{2}=\int_{\Omega}\rho_{m}|{\mathbf u}_m|^{2}|\nabla^{2}{\mathbf u}_{m}||\partial_t{\mathbf u}_{m}|dx\\
		&\;\;\;\; \leqslant C\|\partial_t{\mathbf u}_{m}\|_{{\bf L}^{6}(\Omega)}\|\nabla^{2} {\mathbf u}_{m}\|_{{\bf L}^{2}(\Omega)}\|{\mathbf u}_{m}\|_{{\bf L}^{6}(\Omega)}^{2}\\
		&\;\;\;\;\leqslant C\|\nabla \partial_t{\mathbf u}_{m}\|_{{\bf L}^{2}(\Omega)}\| {\mathbf u}_{m}\|_{{\bf H}^{2}(\Omega)}\|{\mathbf u}_{m}\|^{2}_{{\bf W}_{0}^{1,p}(\Omega)}\\
		&\;\;\;\;\leqslant\frac{\eta}{7}\|\nabla \partial_t{\mathbf u}_{m}\|^{2}_{{\bf L}^{2}(\Omega)}+C\|{\mathbf u}_{m}\|^{2}_{{\bf H}^{2}(\Omega)}\| {\mathbf u}_{m}\|^{4}_{{\bf W}_{0}^{1,p}(\Omega)},\\
		&K_{3}=\int_{\Omega}\rho_{m}|{\mathbf u}_{m}|^{2}|\nabla {\mathbf u}_{m}||\nabla \partial_t{\mathbf u}_{m}|dx\\
		&\;\;\;\;\leqslant C\|\nabla \partial_t{\mathbf u}_{m}\|_{{\bf L}^{2}(\Omega)}\|\nabla {\mathbf u}_{m}\|_{{\bf L}^{6}(\Omega)}\|{\mathbf u}_{m}\|^{2}_{{\bf L}^{6}(\Omega)}\\
	&\;\;\;\;\leqslant C\|\nabla \partial_t{\mathbf u}_{m}\|_{{\bf L}^{2}(\Omega)}\| {\mathbf u}_{m}\|_{{\bf H}^{2}(\Omega)}\|{\mathbf u}_{m}\|^{2}_{{\bf W}_{0}^{1,p}(\Omega)}\\
	&\;\;\;\;\leqslant\frac{\eta}{7}\|\nabla \partial_t{\mathbf u}_{m}\|^{2}_{{\bf L}^{2}(\Omega)}+C\|{\mathbf u}_{m}\|^{2}_{{\bf H}^{2}(\Omega)}\| {\mathbf u}_{m}\|^{4}_{{\bf W}_{0}^{1,p}(\Omega)},\\
		&K_{4}=\int_{\Omega}\rho_{m}|\partial_t{\mathbf u}_{m}|^{2}|\nabla {\mathbf u}_{m}|dx\\
		&\;\;\;\;\leqslant C\|\partial_t{\mathbf u}_{m}\|^{2}_{{\bf L}^{\frac{2p}{p-1}}(\Omega)}\|\nabla {\mathbf u}_m\|_{{\bf L}^{p}(\Omega)}\\
		&\;\;\;\;\leqslant\frac{\eta}{7}\|\nabla \partial_t{\mathbf u}_{m}\|^{2}_{{\bf L}^{2}(\Omega)}+C\| {\mathbf u}_{m}\|^{2}_{{\bf W}_{0}^{1,p}(\Omega)},\\
		&K_{5}=2\int_{\Omega}\rho_{m}|{\mathbf u}_m||\partial_t{\mathbf u}_{m}||\nabla \partial_t{\mathbf u}_{m}|dx\\
		&\;\;\;\;\leqslant C\|\sqrt{\rho_{m}}\partial_t{\mathbf u}_{m}\|^{\frac{1}{2}}_{{\bf L}^{2}(\Omega)}\|\partial_t{\mathbf u}_{m}\|^{\frac{1}{2}}_{{\bf L}^{6}(\Omega)}\|\nabla\partial_t{\mathbf u}_{m}\|_{{\bf L}^{2}(\Omega)}\| {\mathbf u}_{m}\|_{{\bf L}^{6}(\Omega)}\\
		&\;\;\;\;\leqslant C\|\sqrt{\rho_{m}}\partial_t{\mathbf u}_{m}\|^{\frac{1}{2}}_{{\bf L}^{2}(\Omega)}\|\nabla\partial_t{\mathbf u}_{m}\|^{\frac{3}{2}}_{{\bf L}^{2}(\Omega)}\| {\mathbf u}_{m}\|_{{\bf W}_{0}^{1,p}(\Omega)}\\
		&\;\;\;\;\leqslant\frac{\eta}{7}\|\nabla \partial_t{\mathbf u}_{m}\|^{2}_{{\bf L}^{2}(\Omega)}+C\|\sqrt{\rho}\partial_t{\mathbf u}_{m}\|^{2}_{{\bf L}^{2}(\Omega)}\| {\mathbf u}_{m}\|^{4}_{{\bf W}_{0}^{1,p}(\Omega)},\\
		&K_{6}=\int_{\Omega}\nu^{\prime}(\phi_{m})|\partial_t\phi_{m}|(1+|{\mathbb D}{\mathbf u}_{m}|^{2})^{\frac{p-2}{2}}|{\mathbb D}{\mathbf u}_{m}||\nabla \partial_t{\mathbf u}_{m}|dx\\
		&\;\;\;\;\leqslant C\|\partial_t\phi_{m}\|_{{\bf L}^{6}(\Omega)}\|(1+|{\mathbb D}{\mathbf u}_{m}|^{2})^{\frac{p-2}{4}}\nabla \partial_t{\mathbf u}_{m}\|_{{\bf L}^{2}(\Omega)}\|(1+|{\mathbb D}{\mathbf u}_{m}|^{2})^{\frac{p-2}{4}}{\mathbb D}{\mathbf u}_{m}\|_{{\bf L}^{3}(\Omega)}\\
		&\;\;\;\;\leqslant C\|(1+|{\mathbb D}{\mathbf u}_{m}|^{2})^{\frac{p-2}{4}}\nabla \partial_t{\mathbf u}_{m}\|_{{\bf L}^{2}(\Omega)}\|\nabla {\mathbf u}_{m}\|^{\frac{p}{2}}_{{\bf L}^{\frac{3p}{2}}(\Omega)}\|\partial_t\phi_{m}\|_{H^{1}(\Omega)}\\
		&\;\;\;\;\leqslant C\|(1+|{\mathbb D}{\mathbf u}_{m}|^{2})^{\frac{p-2}{4}}\nabla \partial_t{\mathbf u}_{m}\|_{{\bf L}^{2}(\Omega)}\|\nabla {\mathbf u}_{m}\|^{\frac{p(4-p)}{2(6-p)}}_{{\bf L}^{p}(\Omega)}\|\nabla {\mathbf u}_{m}\|^{\frac{p}{6-p}}_{{\bf L}^{6}(\Omega)}\|\partial_t\phi_{m}\|_{H^{1}(\Omega)}\\
		&\;\;\;\;\leqslant\frac{\eta}{7}\|(1+|{\mathbb D}{\mathbf u}_{m}|^{2})^{\frac{p-2}{4}}\nabla \partial_t{\mathbf u}_{m}\|^{2}_{{\bf L}^{2}(\Omega)}+C(\|{\mathbf u}_{m}\|^{4-p}_{{\bf W}_{0}^{1,p}(\Omega)}\| {\mathbf u}_{m}\|^{2}_{{\bf H}^{2}(\Omega)}+\|\partial_t\phi_{m}\|^{\frac{6-p}{3-p}}_{H^{1}(\Omega)})\\
			&\;\;\;\;\leqslant\frac{\eta}{7}\|(1+|{\mathbb D}{\mathbf u}_{m}|^{2})^{\frac{p-2}{4}}\nabla \partial_t{\mathbf u}_{m}\|^{2}_{{\bf L}^{2}(\Omega)}+C(\|{\mathbf u}_{m}\|^{4}_{{\bf W}_{0}^{1,p}(\Omega)}\| {\mathbf u}_{m}\|^{2}_{{\bf H}^{2}(\Omega)}+\|\partial_t\phi_{m}\|^{\frac{6-p}{3-p}}_{H^{1}(\Omega)}),\\
		&K_{7}=2\int_{\Omega}|\nabla\phi_{m}||\nabla\partial_t\phi_{m}||\nabla \partial_t{\mathbf u}_{m}|dx\\
		&\;\;\;\;\leqslant C\|\nabla \partial_t{\mathbf u}_{m}\|_{{\bf L}^{2}(\Omega)}\|\nabla\phi_{m}\|_{{\bf L}^{\infty}(\Omega)}\|\nabla\partial_t\phi_{m}\|_{{\bf L}^{2}(\Omega)}\\
		&\;\;\;\;\leqslant C\|\nabla \partial_t{\mathbf u}_{m}\|_{{\bf L}^{2}(\Omega)}\|\nabla^{2}\phi_{m}\|_{{\bf L}^{\frac{6-p}{3-p}}(\Omega)}\|\nabla\partial_t\phi_{m}\|_{{\bf L}^{2}(\Omega)}\\
		&\;\;\;\;\leqslant C\|\nabla \partial_t{\mathbf u}_{m}\|_{{\bf L}^{2}(\Omega)}\|\nabla^{2}\phi_{m}\|^{\frac{6-2p}{6-p}}_{{\bf L}^{2}(\Omega)}\|\nabla^{2}\phi_{m}\|^{\frac{p}{6-p}}_{{\bf L}^{\infty}(\Omega)}\|\nabla\partial_t\phi_{m}\|_{H^{1}(\Omega)}\\
		&\;\;\;\;\leqslant C\|\nabla \partial_t{\mathbf u}_{m}\|_{{\bf L}^{2}(\Omega)}\|\phi_{m}\|^{\frac{6-2p}{6-p}}_{ H^{2}(\Omega)}\|\phi_{m}\|^{\frac{p}{6-p}}_{W^{2,\infty}(\Omega)}\|\nabla\partial_t\phi_{m}\|_{H^{1}(\Omega)}\\
		&\;\;\;\;\leqslant\frac{\eta}{7}\|\nabla \partial_t{\mathbf u}_{m}\|^{2}_{{\bf L}^{2}(\Omega)}+C(\|\partial_t\phi_{m}\|^{\frac{6-p}{3-p}}_{H^{1}(\Omega)}+\|\phi_{m}\|^{2}_{W^{2,\infty}(\Omega)}),
	\end{align*}
where $\eta\in (0,1)$ is fixed any.	Obviously, it is deduced from \eqref{E4-22} that
\begin{align*}
&\frac{d}{dt}\int_{\Omega}\rho|\partial_t{\mathbf u}_{m}|^{2}dx+\int_{\Omega}\nu(\phi_{m})\left((1+|{\mathbb D}{\mathbf u}_{m}|^{2})^{\frac{p-4}{2}}|{\mathbb D}{\mathbf u}_{m}|^{2}|\nabla \partial_t{\mathbf u}_{m}|^{2}+(1+|{\mathbb D}{\mathbf u}_{m}|^{2})^{\frac{p-2}{2}}|\nabla \partial_t{\mathbf u}_{m}|^{2}\right)dx\notag\\
&\leqslant C \left[\|{\mathbf u}_{m}\|^{4}_{{\bf W}_{0}^{1,p}(\Omega)}\left(\|{\mathbf u}_{m}\|^{2}_{{\bf H}^{2}(\Omega)}+\|\sqrt{\rho}\partial_t{\mathbf u}_{m}\|^{2}_{{\bf L}^{2}(\Omega)}+1\right)
+\|\partial_t\phi_{m}\|^{\frac{6-p}{3-p}}_{H^{1}(\Omega)}
+\|\phi_{m}\|^{2}_{W^{2,\infty}(\Omega)}\right].
\end{align*}
Based on \eqref{E4-20} and \eqref{E4-25}, one uses Gronwall's inequality to arrive at
\begin{equation*}
	\|\sqrt{\rho_{m}}\partial_t{\mathbf u}_{m}\|_{{\bf L}^{\infty}(0,T_{0};{\bf L}^{2}(\Omega))}\leqslant C
\end{equation*}
and
\begin{equation*}
	\|\nabla \partial_t{\mathbf u}_{m}\|_{{\bf L}^{2}(0,T_{0};{\bf L}^{2}(\Omega))}\leqslant C\|\sqrt{\nu(\phi_{m})}(1+|{\mathbb D}{\mathbf u}_{m}|^{2})^{\frac{p-2}{4}}\nabla \partial_t{\mathbf u}_{m}\|_{{\bf L}^{2}(0,T_{0};{\bf L}^{2}(\Omega))}\leqslant C,
\end{equation*}
where the positive constant $C$ is independent of $m$ and $\delta.$ Moreover, it is deduced from Lemma \ref{L-2-3} that
\begin{equation*}
\|\partial_t{\mathbf u}_{m}\|_{\mathbf{L}^{2}(0,T_{0};\mathbf{W}_{0}^{1,2}(\Omega))}\leqslant C.
\end{equation*}
\end{proof}

In summary, there exists $T_0\in(0,T),$ independent of $m$ and $\delta,$ such that
\begin{eqnarray}\label{E4-30}
	\begin{cases}
		0<C\delta\leqslant\rho_{m}(x,t)\leqslant\rho^{*}+1\;(\forall(x,t)\in Q_{T}),~~\|\partial_{t}\rho_{m}\|_{ L^{\infty}(0,T_{0};H^{-1}(\Omega))}\leqslant C,\\
		\|\sqrt{\rho_{m}}{\mathbf u}_{m}\|_{{\bf L}^{\infty}(0,T_{0};{\bf L}^{2}(\Omega))}\leqslant C,~~\|\partial_{t}{\mathbf u}_{m}\|_{{\bf L}^{2}(0,T_{0};{\bf W}_{0}^{1,2}(\Omega))}\leqslant C,\\
		\|{\mathbf u}_{m}\|_{{\bf L}^{\infty}(0,T_{0};{\bf W}_{0}^{1,p}(\Omega))}\leqslant C,~~\|{\mathbf u}_{m}\|_{{\bf L}^{2}(0,T_{0};{\bf H}^{2}(\Omega))}\leqslant C,\\
		\|\sqrt[4]{\rho_{m}}\phi_{m}\|_{L^{\infty}(0,T_{0};L^{4}(\Omega))}\leqslant C,~~\|\sqrt{\rho_{m}}\phi_{m}\|_{L^{\infty}(0,T_{0};L^{2}(\Omega))}\leqslant C,\\
		\|\phi_{m}\|_{L^{\infty}(0,T_{0};H^{2}(\Omega))}\leqslant C,~~\|\phi_{m}\|_{ L^{2}(0,T_{0};W^{2,\infty}(\Omega))}\leqslant C,~~\\
		\|\partial_t\phi_{m}\|_{L^{\frac{6-p}{3-p}}(0,T_{0};H^{1}(\Omega))}\leqslant C,\\
		\|\mu_{m}\|_{L^{2}(0,T_{0};H^{2}(\Omega))}\leqslant C,
		~~ \|\mu_{m}\|_{L^{\infty}(0,T_{0};H^{1}(\Omega))}\leqslant C,
	\end{cases}
\end{eqnarray}
where $C$ is independent of $m$ and $\delta.$

\subsection{Passage to the limit}

\quad Now, we are in a position to take limits as $m\rightarrow\infty$ and then as $\delta\rightarrow0^{+}.$ Clearly, one deduces from \eqref{E4-30} and Lions-Aubin Lemma that
\begin{eqnarray}\label{E-4-28}
\begin{cases}
	{\mathbf u}_{m}\rightarrow{\mathbf u}^{\delta} \;strongly\; in\; {\bf C}([0,T_{0}];{\bf L}^{q}(\Omega))\;(\forall q\in[1,\frac{3p}{3-p})),\\
    {\mathbf u}_{m}\rightarrow{\mathbf u}^{\delta} \;strongly\; in\; {\bf L}^2(0,T_{0};{\bf W}_0^{1,p}(\Omega)),\\
    \rho_{m}\rightarrow\rho^{\delta} \;strongly\; in\; C([0,T_{0}];L^{r}(\Omega))\;(\forall r\in(1,\infty)),\\
    \phi_{m}\rightarrow\phi^\delta\;strongly\;in\; C([0,T_0];L^{r}(\Omega))\;\;(\forall\;r\in(1,\infty)),\\
    \phi_{m}\rightarrow\phi^\delta\;strongly\;in\; L^2(0,T_0;W^{1,r}(\Omega))\;\;(\forall\;r\in(1,\infty)),\\
    \mu_{m}\rightarrow\mu^\delta\;weakly\;in\; L^{\infty}(0,T_{0};H^{1}(\Omega)).
\end{cases}
\end{eqnarray}
So, one gets that
\begin{eqnarray}\label{E-4-29}
\begin{cases}
	\rho_{m}{\mathbf u}_{m}\rightarrow\rho^{\delta}{\mathbf u}^{\delta} \;strongly\; in\; {\bf C}([0,T_{0}];{\bf L}^{2}(\Omega)),\\
    \rho_{m}{\mathbf u}_{m}\otimes{\mathbf u}_{m}\rightarrow\rho^{\delta}{\mathbf u}^{\delta}\otimes{\mathbf u}^{\delta} \;strongly\;in\; {\bf C}([0,T_{0}];{\bf L}^{2}(\Omega)),\\
    \nabla\phi_{m}\otimes\nabla\phi_{m}\rightarrow\nabla\phi^\delta\otimes\nabla\phi^\delta\;strongly\;in\;{\bf C}([0,T_{0}];{\bf L}^{2}(\Omega)),\\
	\rho_{m}\phi_{m}\rightarrow\rho^\delta\phi^\delta\;strongly\;in\;C([0,T_{0}];L^{2}(\Omega)),\\
	\rho_{m}\Psi^\prime(\phi_{m})\rightarrow\rho^\delta\Psi^\prime(\phi^\delta)\;strongly\;in\;C([0,T_{0}];L^{2}(\Omega)),\\
	\rho_{m}\Psi^\prime(\phi_{m})\nabla\phi_{m}\rightarrow\rho^\delta\Psi^\prime(\phi^\delta)\nabla\phi^\delta\;strongly\;in\;{\bf C}([0,T_{0}];{\bf L}^{2}(\Omega)),\\
    \rho_{m}{\mathbf u}_{m}\phi_{m}\rightarrow \rho^\delta  {\mathbf u}^\delta\phi^\delta\;strongly\;in\;{\bf L}^{2}(0,T_{0};{\bf L}^{2}(\Omega)),\\
    \rho_{m}\mu_{m}\rightarrow\rho^{\delta}\mu^{\delta}\;weakly\;in\;L^{2}(0,T;L^{2}(\Omega)),\\
    \rho_{m}\mu_{m}\nabla\phi_{m}\rightarrow\rho^{\delta}\mu^{\delta}\nabla\phi^{\delta}\;weakly\;in\;{\bf L}^{2}(0,T;{\bf L}^{\frac{3}{2}}(\Omega)).
\end{cases}
\end{eqnarray}
Recall that
$${\cal T}({\mathbf u}_{m})\triangleq(1+|{\mathbb D}{\mathbf u}_m|^2)^\frac{p-2}{2}{\mathbb D}{\mathbf u}_m.$$
It follows from the convergence of $\phi_{m}$ and $\mathbf{u}_{m}$ obtained in \eqref{E-4-28} that
\begin{eqnarray}\label{E-4-30}
\nu(\phi_{m}){\cal T}({\mathbf u}_{m})\rightarrow\nu(\phi^{\delta}){\cal T}({\mathbf u}^\delta)\;strongly\;in\;{\bf L}^{\frac{2p}{3p-2}}(0,T;{\bf L}^{\frac{p}{p-1}}(\Omega)).
\end{eqnarray}

By the weak lower semi-continuity of the norm, one can get uniform estimates on the approximate solution $(\rho^{\delta},{\mathbf u}^{\delta},\phi^{\delta},\mu^{\delta}).$ That is, there exists a positive constant $C$, independent of $\delta$, such that
\begin{eqnarray*}
	\begin{cases}
		\|\rho^\delta\|_{L^{\infty}(0,T_{0};L^{\infty}(\Omega))}\leqslant C,~~\|\partial_{t}\rho^\delta\|_{L^{\infty}(0,T_{0};H^{-1}(\Omega))}\leqslant C,\\
		\|\sqrt{\rho^\delta}{\mathbf u}^\delta\|_{{\bf L}^{\infty}(0,T_{0};{\bf L}^{2}(\Omega))}\leqslant C,
		~\|{\mathbf u}^\delta\|_{{\bf L}^{\infty}(0,T_{0};{\bf W}_{0}^{1,p}(\Omega))}\leqslant C,\\
		\|{\mathbf u}^\delta\|_{{\bf L}^{2}(0,T_{0};{\bf H}^{2}(\Omega))}\leqslant C,
		~\|\partial_{t}{\mathbf u}^\delta\|_{{\bf L}^{2}(0,T_{0};{\bf W}_{0}^{1,2}(\Omega))}\leqslant C,\\
		\|\phi^\delta\|_{ L^{2}(0,T_{0};W^{2,\infty}(\Omega))}\leqslant C,~\|\phi^\delta\|_{L^{\infty}(0,T_{0};H^{2}(\Omega))}\leqslant C,
		\|\partial_t\phi^\delta\|_{L^{2}(0,T_{0};H^{1}(\Omega))}\leqslant C, \\
		\|\mu^\delta\|_{L^{2}(0,T_{0};H^{2}(\Omega))}\leqslant C,~\|\mu^\delta\|_{L^{\infty}(0,T_{0};H^{1}(\Omega))}\leqslant C,
		~\|\sqrt{\rho^\delta}\mu^\delta\|_{L^{\infty}(0,T_{0};L^{2}(\Omega))}\leqslant C.
	\end{cases}
\end{eqnarray*}
This allows us to take similar argument in the process $m\to+\infty$,  pass to the limit as $\delta\to 0^{+}$, and to obtain a limit  $(\rho,{\mathbf u},\phi,\mu)$ satisfying the problem \eqref{E1-1}-\eqref{E1-2} as stated in Theorem \ref{thm2-2}.

\begin{corollary}\label{cor4-8}
Let $(\rho,{\mathbf u},\phi,\mu)$  be the local strong solution to the problem \eqref{E1-1}-\eqref{E1-2} obtained in Theorem \ref{thm2-2}. Then 	\begin{equation*}
	P\in L^{2}(0,T_{0};W^{1,\frac{2p}{3p-4}}(\Omega)).
\end{equation*}
\end{corollary}

\begin{proof}
One can rewrite $\eqref{E1-1}_{2}$ and obtain that
\begin{align*}
\nabla P
   &=-div(\nabla\phi\otimes\nabla\phi)-\rho\partial_{t}{\mathbf u}+\rho ({\mathbf u}\cdot\nabla ){\mathbf u}
       +div(\nu(\phi)(1+|\mathbb{D}{\mathbf u}|^{2})^{\frac{p-2}{2}}\mathbb{D}{\mathbf u})\\
  &=\nabla\phi\nabla^{2}\phi+\Delta\phi\nabla\phi-\rho\partial_{t}{\mathbf u}+\rho({\mathbf u}\cdot\nabla) {\mathbf u}
   +\nabla\nu(\phi)(1+|\mathbb{D}{\mathbf u}|^{2})^{\frac{p-2}{2}}\mathbb{D}{\mathbf u}\\
  &+(p-2)\nu(\phi)(1+|\mathbb{D}{\mathbf u}|^{2})^{\frac{p-4}{2}}\nabla \mathbb{D}{\mathbf u}:\mathbb{D}{\mathbf u}
  +\nu(\phi)(1+|\mathbb{D}{\mathbf u}|^{2})^{\frac{p-2}{2}}\Delta{\mathbf u}.
\end{align*}
It is deduced from $(1+|\mathbb{D}{\mathbf u}|^{2})^{\frac{p-2}{2}}\in\mathbf{L}^{\infty}(0,T_{0};\mathbf{L}^{\frac{p}{p-2}}(\Omega))$ and $\Delta{\mathbf u}\in\mathbf{L}^{2}(0,T_{0};\mathbf{L}^{2}(\Omega))$ that
\begin{eqnarray*}
 &&\nu(\phi)(1+|\mathbb{D}{\mathbf u}|^{2})^{\frac{p-2}{2}}\Delta{\mathbf u}\in\mathbf{L}^{2}(0,T_{0};\mathbf{L}^{\frac{2p}{3p-4}}(\Omega)),\\
 &&\nu(\phi)(1+|\mathbb{D}{\mathbf u}|^{2})^{\frac{p-4}{2}}\nabla \mathbb{D}{\mathbf u}:\mathbb{D}{\mathbf u}\in\mathbf{L}^{2}(0,T_{0};\mathbf{L}^{\frac{2p}{3p-4}}(\Omega)).
\end{eqnarray*}
Note that $\nabla\nu(\phi)=\nu^{\prime}(\phi)\nabla\phi\in \mathbf{L}^{2}(0,T_{0};\mathbf{L}^{\infty}(\Omega)).$ One finds that
\begin{equation*}
 \nabla\nu(\phi)(1+|\mathbb{D}{\mathbf u}|^{2})^{\frac{p-2}{2}}\mathbb{D}{\mathbf u}\in\mathbf{L}^{2}(0,T_{0};\mathbf{L}^{\frac{p}{p-1}}(\Omega)).
\end{equation*}
Moreover, $$\rho\partial_{t}{\mathbf u},\;\rho{\mathbf u}\cdot\nabla {\mathbf u},\;\nabla\phi\cdot\nabla^{2}\phi+\Delta\phi\nabla\phi\in{\mathbf L}^{2}(0,T_{0};\mathbf{L}^{p}(\Omega))\subset {\mathbf L}^{2}(0,T_{0};\mathbf{L}^{\frac{2p}{3p-4}}(\Omega)).$$
Thus,
\begin{equation*}
\nabla P\in\mathbf{L}^{2}(0,T_{0};\mathbf{L}^{\frac{2p}{3p-4}}(\Omega)).
\end{equation*}
Further, one deduces from Lemma \ref{lem2-3}
\begin{equation*}
	P \in L^{2}(0,T_{0};W^{1,\frac{2p}{3p-4}}(\Omega)).
\end{equation*}
	\end{proof}

\section{Appendix }

{\bf Lemma A.1.} {\sl Let $\Omega\subset {\mathbb R}^3$ be a bounded domain with $\partial\Omega\in {\mathcal C}^2,$ $1<q<+\infty\;and \;{\bf v}\in {\bf W}^{1,q}(\Omega)$. If $\rho$ is a non-negative function satisfying
\begin{eqnarray}\label{L-2-3}
0<M<\int_{\Omega}\rho dx~~\mbox{ and }~~\int_{\Omega}\rho^{\gamma}dx\leqslant E_{0}
\end{eqnarray}
for some $\gamma>1$, then there exists a constant $C(M,E_{0})$ such that
	\begin{equation}
		\|{\bf v}\|^{q}_{{\bf L}^{q}(\Omega)}\leqslant C(M,E_{0})\left(\|\nabla {\bf v}\|^{q}_{{\bf L}^{q}(\Omega)}+(\int_{\Omega}\rho|{\bf v}|dx)^{q}\right).
	\end{equation}
}

\begin{proof}
Assume the conclusion is incorrect, then there  exist a  $\{\rho_{n}\}^{\infty}_{n=1}$ satisfy \eqref{L-2-3} and $\{{\bf v}_{n}\}^{\infty}_{n=1}\subset W^{1,p}(\Omega)$ such that \begin{equation*}
		\|{\bf v}_{n}\|^{p}_{L^{p}(\Omega)}\geqslant c_{n}(\|\nabla {\bf v}_{n}\|^{p}_{L^{p}(\Omega)}+(\int_{\Omega}\rho_{n}|v_{n}|dx)^{p})\mbox{ with }\;c_{n}\rightarrow+\infty
	\end{equation*}
	Setting ${\bf w}_{n}=\frac{{\bf v}_{n}}{\|{\bf v}_{n}\|_{L^{p}(\Omega)}},$ one has that $\|{\bf w}_{n}\|_{L^{p}(\Omega)}=1$ and
	$${\bf w}_{n}\rightarrow {\bf w}~~ in~~ L^{q}(B) \mbox{  for every compact set } B\subset \Omega\mbox{ and any } 2\leqslant q\leqslant p^*,$$
	where $p^*$ is the critical Sobolev exponent and
	$$p^*=\frac{3p}{3-p}~\mbox{for }1<p<3\mbox{ and }p^*~\mbox{arbitrary finite}\mbox{for }p\geqslant 3.$$
	Since
	\begin{equation*}
		\frac{1}{c_{n}}\geqslant
		\int_{\Omega}\frac{|\nabla {\bf v}_{n}|^{p}}{\|{\bf v}_{n}\|^{p}_{L^{p}(\Omega)}}dx+(\int_{\Omega}\frac{\rho_{n}|{\bf v}_{n}|}{\|{\bf v}_{n}\|_{L^{p}(\Omega)}}dx)^{p}
		=\int_{\Omega}\nabla|{\bf w}_{n}|^{p}dx+(\int_{\Omega}\rho_{n}|{\bf w}_{n}|dx)^{p},
	\end{equation*}
	one takes $n\rightarrow\infty$ to arrive at $\int_{\Omega}\nabla|{\bf w}|^{p}dx+(\int_{\Omega}\rho|{\bf w}|dx)^{p}\rightarrow0.$ So,
	${\bf w}$ is a constant $c_0.$
	
	On the other hand, we introduce a cut-off function $T_K$ defined as
	$$T_k(z)\equiv kT(\frac{z}{k}),$$
	where $T\in C^\infty({\mathbb R})$ is chosen such that
	\begin{eqnarray*}
		T(z)\equiv\begin{cases}T(z)=z&\mbox{for}~~z\in[0,1],\\
			T(z)\ \mbox{ concave}&\mbox{on}~~[0,+\infty),\\
			T(z)=2&\mbox{for}~~z\geqslant 2,\\
			T(z)=-T(-z)&\mbox{for}~~z\in(-\infty,0].
		\end{cases}
	\end{eqnarray*}
	So, there exists a $k=k(M,E_{0})$ large enough such that
	$$T_{k}(\rho_{n})\rightarrow\overline{T_{k}(\rho)}\mbox{ weakly in }L^{\beta}(\Omega)\mbox{ for any finite }\beta\geqslant 1,$$
	where
	$$\int_{\Omega}\overline{T_{k}(\rho)}>\frac{M}{2}.$$
	Obviously,
	\begin{align*}
		0=\lim\limits_{n\rightarrow\infty}\int_{B}\rho_{n}|{\bf w}_{n}|dx
		&\geqslant\lim\limits_{n\rightarrow\infty}\int_{B}T_{k}(\rho_{n})|{\bf w}_{n}|dx\\
		&=\int_{\Omega}\overline{T_{k}(\rho)}|{\bf w}|dx=c_0\int_{B}\overline{T_{k}(\rho)}dx>0.
	\end{align*}
	holds for every compact set $B\subset \Omega.$ It is a contradiction and the Proof of Lemma 2.4 is completed.
\end{proof}

\quad Given $\rho_{0}\in L^{\infty}(\Omega)$ with $0\leqslant \rho_0\leqslant\rho^*$ for almost every $x\in\Omega,$ there exist a sequence $\{\rho_{0\delta}\}$  such that
\begin{equation}\label{E-3-1}
	\rho_{0\delta}\in C^{\infty}(\overline{\Omega}),\;\;0<\delta\leqslant\rho_{0\delta}(x)\leqslant\rho^{*}+1\;\;(\forall x\in\overline{\Omega})
\end{equation}
and
\begin{equation}\label{E-3-2}
\rho_{0\delta}\rightarrow\rho_{0}\;strongly\;in\;L^{ r }(\Omega)~(\forall r \in[1,\infty)),\;\; \;
\rho_{0\delta}\rightarrow\rho_{0}\;weak-star\;in\;L^{\infty}(\Omega)
\end{equation}
as $\delta\rightarrow 0^+.$ Consider the family of eigenfunctions $\{w_{j}\}_{j=1}^{\infty}$ and eigenvalues $\{\lambda_{j}\}_{j=1}^{\infty}$ of the Laplace operator $A=-\Delta+I$ with homogeneous Neumann boundary condition and family of eigenfunctions $\{{\mathbf w}_{j}\}_{j=1}^{\infty}$ and eigenvalues $\{\lambda^{S}_{j}\}_{j=1}^{\infty}$ of the Stokes operator $A$. For any integer $m\geqslant 1$, the $m$-dimensional orthogonal projections $\Pi_{m}$ and $P_{m}$ are defined on $V_{m}$ and ${\mathbf V}_{m}$ with respect to the inner product in $L^{2}(\Omega)$ and in ${\mathbf H}_{\sigma},$ respectively. Let
\begin{equation}\label{(chuzhi2)}
{\mathbf u}_{0m}=P_{m}{\mathbf u}_{0}~~~~~~~~\mbox{        and         }~~~~~~~~~~~~~~\phi_{0m}=\Pi_{m}\phi_{0}.
\end{equation}
then
\begin{equation*}
	{\mathbf u}_{0m}\rightarrow{\mathbf u}_{0}\;strongly\;in\;{\bf L}^{2}(\Omega)~~
\mbox{        and         }~~\phi_{0m}\rightarrow\phi_{0}\;strongly\;in\;H^{1}(\Omega)
\end{equation*}
respectively as $m\to+\infty.$
Then, we are looking for $(\rho_{m},{\mathbf u}_{m},\phi_{m},\mu_{m})$ to solve the following system
\begin{eqnarray}\label{E3-3}
\begin{cases}
\partial_{t}\rho_{m}+{\mathbf u}_{m}\cdot\nabla\rho_{m}=0&\mbox{in}\;\Omega\times(0,T),\\
(\rho_{m}\partial_{t}{\mathbf u}_{m},{\mathbf w})+(\rho_{m}{\mathbf u}_{m}\cdot\nabla{\mathbf u}_{m},{\mathbf w})
    +(\nu(\phi_{m})(1+|{\mathbb D}{{\mathbf u}}_{m}|^{2})^{\frac{p-2}{2}}\mathbb{D}{\mathbf u}_{m},\nabla {\mathbf w})\\
~~~~~~~~~~~~~~~~=(\rho_{m}\mu_{m}\nabla\phi_{m},{\mathbf w})-(\rho_{m}\nabla\Psi(\phi_{m}),{\mathbf w})~~~~({\bf w}\in {\bf V}_m)&\mbox{in}\;(0,T),\\
(\rho_{m}\partial_{t}\phi_{m},w)+(\rho_{m} {\mathbf u}_{m}\cdot\nabla\phi_{m},w)+(\nabla\mu_{m},\nabla w)=0~~~~(w\in V_m)&\mbox{in}\;(0,T)\\
(\rho_{m}\mu_{m},w)=(\nabla\phi_{m},\nabla w)+(\rho_m\Psi^\prime(\phi_{m}),w)~~~~(w\in V_m)&\mbox{in}\;(0,T),\\
{\mathbf u}_{m}=0,\partial_{n}\mu_{m}=\partial_{n}\phi_{m}=0 &\mbox{ on }\partial\Omega\times (0,T),\\
\rho_{m}(\cdot,0)=\rho_{0\delta},{\mathbf u}_{m}(\cdot,0)={\mathbf u}_{0m},\phi_{m}(\cdot,0)=\phi_{0m}&\mbox{ in }\Omega.
\end{cases}
\end{eqnarray}
The following proposition states that the semi-Galerkin scheme \eqref{E3-3} admits a unique solution.

{\bf Proposition A.1.}{\sl(Solvability of the semi-Galerkin scheme). Let $1<p<+\infty.$ Suppose that $T\in(0,+\infty),$
\begin{eqnarray*}\label{E1-3}
0<\nu_*\leqslant \nu(s)\leqslant \nu^* (\forall s\in {\mathbb R}),\;\;\;
\Psi(s)=\frac14(s^2-1)^2 (\forall s\in {\mathbb R}).
\end{eqnarray*}
Given $m\in{\mathbb Z}^+,$ the semi-Galerkin scheme \eqref{E3-3} with \eqref{(chuzhi2)} admits a unique solution $(\rho_m,{\mathbf u}_m,\phi_m,\mu_m)$ on $[0,T]$ satisfying
\begin{eqnarray*}
\rho_m\in C([0,T]\times \overline{\Omega}),~~{\bf u}_m\in{\bf H}^1(0,T;{\bf V}_m),~~
\phi_m\in H^1(0,T;V_m),~~\mu_m\in C([0,T];V_m).
\end{eqnarray*}
}

\begin{proof}
For any $m\in {\mathbb Z}^+,$ fixed
\begin{eqnarray}\label{A-1}\widetilde{{\mathbf u}}_{m}\in {\bf C}([0,T];{\bf V}_{m}),~~~\widetilde{\phi}_{m}\in C([0,T];V_{m}),\end{eqnarray}
with
\begin{equation*}\|\widetilde{{\mathbf u}}_{m}\|_{{\bf L}^\infty(\Omega)}\leqslant \widetilde{M}~~~~~~\mbox{      and       }~~~~~~
\|\widetilde{\phi}_{m}\|_{L^\infty(\Omega)}\leqslant \widetilde{M}.\end{equation*}
where $\widetilde{M}$ satisfies
\begin{equation*}
	\widetilde{M}
	>2\left(\int_{\Omega}\left(\frac{1}{2}\rho_{0\delta}|{\mathbf u}_{0m}|^{2}+\frac{1}{2}|\nabla\phi_{0m}|^{2}+\rho_{m}\Psi(\phi_{0m})\right)dx+1\right).
\end{equation*}and the exact value of $\widetilde{M}$ will be determined later. For
\begin{eqnarray}\label{E3-4-1}
\begin{cases}
\partial_{t}\rho_{m}+\widetilde{{\mathbf u}}_{m}\cdot\nabla\rho_{m}=0,\\
\rho_{m}(\cdot,0)=\rho_{0\delta}~~~~~\mbox{ in }\Omega,
\end{cases}
\end{eqnarray}
it is deduced from \cite{Ladyzenskaja-1975}, Lemma 1.3 that
\begin{equation}\label{E3-4}
	\rho_{m}(x,t)=\rho_{0\delta}(\widetilde{{\mathbf X}}_{m}(0,t,x))
\end{equation}
where
\begin{equation*}
	\widetilde{{\mathbf X}}_{m}(s,t,x)=x+\int_{t}^{s}\widetilde{{\mathbf u}}_{m}(\widetilde{{\mathbf X}}_{m}(\tau,t,x),\tau)d\tau\;\;\;(\forall s,t\in[0,T]).
\end{equation*}
The solution $\rho_{m}$ satisfies the following estimates
\begin{equation}\label{E-2}
	0<C_1\delta\leqslant\rho_{m}(x,t)\leqslant\rho^{*}+1\;\;\;(\forall(x,t)\in Q_{T}).
\end{equation}

Next, for given $\rho^m$ as \eqref{E3-4}, we look for the triple $({\mathbf u}_{m},\phi_{m},\mu_{m})$ defined as
\begin{eqnarray*}
{\mathbf u}_{m}(x,t)=\sum\limits_{j=1}^{m}a_{j}^{m}(t){\mathbf w}_{j}(x),\;
\phi_{m}(x,t)=\sum\limits_{j=1}^{m}b_{j}^{m}(t)w_{j}(x),\;
\mu_{m}(x,t)=\sum\limits_{j=1}^{m}c_{j}^{m}(t)w_{j}(x),
\end{eqnarray*}
to solution the following system
\begin{eqnarray}\label{E3-5}
\begin{cases}
(\rho_{m}\partial_{t}{\mathbf u}_{m},{\mathbf w}_l)+(\rho_{m}\left(\widetilde{{\mathbf u}}_{m}\cdot\nabla\right){\mathbf u}_{m},{\mathbf w}_l)
    +(\nu(\widetilde{\phi}_{m})(1+|{\mathbb D}{{\mathbf u}}_{m}|^{2})^{\frac{p-2}{2}}\mathbb{D}{\mathbf u}_{m},\nabla {\mathbf w}_l)\\
~~~~~~~~~~~~~~~~=(\rho_{m}\mu_{m}\nabla\widetilde{\phi}_{m},{\mathbf w}_l)-(\rho_{m}\nabla\Psi(\widetilde{\phi}_{m}),{\mathbf w}_l),\\
(\rho_{m}\partial_{t}\phi_{m},w_l)+(\rho_{m} {\mathbf u}_{m}\cdot\nabla\widetilde{\phi}_{m},w_l)+(\nabla\mu_{m},\nabla w_l)=0,\\
(\rho_{m}\mu_{m},w_l)=(\nabla\phi_{m},\nabla w_l)+(\rho_m\Psi^\prime(\phi_{m}),w_l),\\
{\mathbf u}_{m}=0,\partial_{n}\mu_{m}=\partial_{n}\phi_{m}=0, ~~~\mbox{ on }\partial\Omega\times (0,T),\\
{\mathbf u}_{m}(\cdot,0)={\mathbf u}_{0m},\phi_{m}(\cdot,0)=\phi_{0m}, ~~~\mbox{ in }\Omega
\end{cases}
\end{eqnarray}
for $l=1,2,\cdots,m.$ Set
\begin{equation*}
	{\mathbf a}^{m}(t)=(a^{m}_{1}(t),...,a^{m}_{m}(t)),\;{\mathbf b}^{m}(t)=(b^{m}_{1}(t),...,b^{m}_{m}(t)),\;{\mathbf c}^{m}(t)=(c^{m}_{1}(t),...,c^{m}_{m}(t)),
\end{equation*}
the system \eqref{E3-5}  can be written as
\begin{equation}\label{E3-6}
	\begin{cases}
		\mathbf{M}^{m}_{1}(t)\frac{d}{dt}\mathbf{a}^{m}={\mathbf L}^{m}_{1}(t,\mathbf{a}^{m})+{\mathbf L}^{m}_{2}(t)\mathbf{c}^{m}+\mathbf{F}^{m}_{1}(t),\\
		\mathbf{M}^{m}_{2}(t)\frac{d}{dt}\mathbf{b}^{m}=\mathbf{L}^{m}_{3}(t)\mathbf{a}^{m}+\mathbf{L}^{m}_{4}(t)\mathbf{c}^{m},\\
		\mathbf{M}^{m}_{2}(t)\mathbf{c}^{m}=-\mathbf{L}^{m}_{4}(t)\mathbf{b}^{m}+\mathbf{F}^{m}_{2}(t,\mathbf{b}^{m}),\\
        {\mathbf a}^{m}(0)={\mathbf a}^{m}_0,~~{\mathbf b}^{m}(0)={\mathbf b}^{m}_0,
	\end{cases}
\end{equation}
where
\begin{align*}
	&(\mathbf{M}^{m}_{1})_{l,j}=\int_{\Omega}\rho_{m}(x,t)\mathbf{w}_{j}(x)\cdot \mathbf{w}_{l}(x)dx,\;\;\;
    (\mathbf{M}^{m}_{2})_{l,j}=\int_{\Omega}\rho_{m}(x,t)w_{j}(x)w_{l}(x)dx,\\
	&({\mathbf L}^{m}_{1}(t,\mathbf{a}^{m}))_{l}\\
     &=\int_{\Omega}\left(\rho_{m}(x,t)(\widetilde{{\mathbf u}}_{m}(x,t)\cdot\nabla) {\mathbf u}_{m}(x,t)\cdot \mathbf{w}_{l}(x)-\nu(\widetilde{\phi}_{m}(x,t))(1+|{\mathbb D}{\mathbf u}_m(x,t)|^{2})^{\frac{p-2}{2}}{\mathbb D}{\mathbf u}_{m}(x,t):{\mathbb D} \mathbf{w}_{l}(x)\right)dx\\
	&(\mathbf{L}^{m}_{2})_{l,j}=\int_{\Omega}\rho_{m}(x,t)w_{j}(x)\nabla\widetilde{\phi}_{m}(x,t)\cdot {\mathbf w}_{l}(x)dx,~~~~~~~~~
    ({\mathbf L}^{m}_{3})_{l,j}=-(\mathbf{L}^{m}_{2})_{j,l},\\
	&(\mathbf{L}^{m}_{4})_{l,j}=\int_{\Omega}\nabla w_{j}(x)\cdot\nabla w_{l}(x)dx=-\int_{\Omega}\Delta w_{j}(x)w_{l}(x)=\int_{\Omega}w_{l}(x)\lambda_{j}w_{j}(x)dx=\lambda_{j},\\
   &(\mathbf{F}^{m}_{1}(t))_{l}=\int_{\Omega}\rho_{m}(x,t)\nabla\Psi(\widetilde{\phi}_{m}(x,t))\cdot {\mathbf w}_{l}(x,t)dx,\\
	&(\mathbf{F}^{m}_{2}(t,\mathbf{b}^{m}))_{l}=\int_{\Omega}\rho_{m}(x,t)\Psi^\prime(\phi_{m}(x,t))w_{l}(x,t)dx
\end{align*}
for $l,j=1,...,m$, where $\lambda_{1},...,\lambda_{m}$ are the eigenvalues of the stokes operator.

Note that ${\mathbf L}^{m}_{1}(t,\mathbf{a}^{m})$ is continuous on $[0,T]\times\mathbb{R}^{m}$ and locally Lipschitz in $\mathbb{R}^{m}$ uniformly in $t,$ $\mathbf{F}^{m}_{1},$ $\mathbf{L}_{2}^{m}$ and $\mathbf{L}_{4}^{m}$ are continuous on $[0,T]$, and $\mathbf{F}_{2}^{m}(t,{\mathbf b}^{m})$ is continuous on $[0,T]\times\mathbb{R}^{m}$ and locally Lipschitz in $\mathbb{R}^{m}$ uniformly in $t$. Moreover, ${\mathbf M}^{m}_{1}$ and ${\mathbf M}^{m}_{2}$ are invertible due to the fact that $\rho_m$ is bounded below by $C\delta,$ and the regularity of $\rho_m$ ensures that ${\mathbf M}^{m}_{1}$ and ${\mathbf M}^{m}_{2}$  are in $C[0,T].$ So, the standard ODE theory implies the solvability of \eqref{E3-5} in $(0,t_m)\subset (0,T).$

Next, we find that the following mapping is well defined
\begin{eqnarray*}
&&{\cal T}^1_m:~~C([0,T];{\bf V}_m)~~\mapsto ~~ C^1([0,T]\times \overline{\Omega}),\\
&&~~~~~~~~~~~~~~~~~~~~~~~~~~~~\widetilde{{\bf u}}_m~~\mapsto \rho_m=\rho_m[\widetilde{{\bf u}}_m].
\end{eqnarray*}
Let ${\mathbf v}^1_m,{\mathbf v}^2_m\in {\bf C}([0,T];{\bf V}_m)$ be fixed such that
\begin{equation*}
	\begin{cases}
	  \partial_{t}\rho_{m}^1+{\mathbf v}^1_m\cdot\nabla\rho_{m}^1=0,\\
      \rho_{m}^1(\cdot,0)=\rho_{0\delta}~~~~~\mbox{ in }\Omega
     \end{cases}
\end{equation*}
and
\begin{equation*}
	\begin{cases}
	  \partial_{t}\rho_{m}^2+{\mathbf v}^2_m\cdot\nabla\rho_{m}^2=0,\\
      \rho_{m}^2(\cdot,0)=\rho_{0\delta}~~~~~\mbox{ in }\Omega.
     \end{cases}
\end{equation*}
Set $$\delta\rho_{m}=\rho^{1}_{m}-\rho_{m}^2.$$
It is easy to see that
\begin{equation}\label{(31a)}
\begin{cases}
	\partial_{t}\delta\rho_{m}+{\mathbf v}^{1}\cdot\nabla\delta\rho_{m}+({\mathbf v}^{1}_m-{\mathbf v}^{2}_m)\cdot\nabla\rho_{m}^2=0,\\
     \delta\rho_{m}(\cdot,0)=0~~~~~\mbox{ in }\Omega.
\end{cases}
\end{equation}
Multiplying \eqref{(31a)} by $|\delta\rho_{m}|^{ r -2}\delta\rho_{m}( r \geqslant 2),$ one gets that
\begin{align*}
	\frac{d}{dt}\int_{\Omega}|\delta\rho_{m}|^{ r }dx
    &=\int_{\Omega}({\mathbf v}^{1}_m-{\mathbf v}^{2}_m)\cdot\nabla\rho_{m}^2(|\delta\rho_{m}|^{ r -2}\delta\rho_{m})dx\\
	&\leqslant\|{\mathbf v}^{1}_m-{\mathbf v}^{2}_m\|_{{\bf L}^{ r }(\Omega)}\|\nabla\rho_{m}^2\|_{L^{\infty}(\Omega)}
    \||\delta\rho_{m}|^{ r -2}\delta\rho_{m}\|_{L^{\frac{ r }{ r -1}}(\Omega)}\\
	&\leqslant\|{\mathbf v}^{1}_m-{\mathbf v}^{2}_m\|_{{\bf L}^{ r }(\Omega)}\|\nabla\rho_{m}^2\|_{L^{\infty}(\Omega)}
    \|\delta\rho_{m}\|^{ r -1}_{L^{ r }(\Omega)}\\
	&\leqslant C\|{\mathbf v}^{1}_m-{\mathbf v}^{2}_m\|_{{\bf L}^{ r }(\Omega)}\|\delta\rho_{m}\|^{ r -1}_{L^{ r }(\Omega)}
\end{align*}
and
\begin{equation}\label{E3-12}
	\|\delta\rho_{m}(t)\|_{L^{ r }(\Omega)}\leqslant C\int_{0}^{t}\|{\mathbf v}^{1}_m-{\mathbf v}^{2}_m\|_{{\bf L}^{ r }(\Omega)}d\tau~~(\forall t\in [0,t_m]).
\end{equation}
Obviously, the mapping
\begin{eqnarray}\label{map1}
{\cal T}^1_m:~~\widetilde{{\bf u}}_m\in {\bf C}([0,T];{\bf V}_m)\mapsto \rho_m[\widetilde{{\bf u}}_m]\in {\bf LC}^1([0,T]\times \overline{\Omega})
\end{eqnarray}
is continuous. Then, the following map is also well defined
\begin{eqnarray*}
&&{\cal T}^2_m:~~{\bf C}([0,T];{\bf V}_m)\times C([0,T];V_m)~~\mapsto ~~ {\bf C}([0,T];{\bf V}_m)\times C([0,T];V_m),\\
&&~~~~~~~~~~~~~~~~~~~~~~~~~~~~~~~~~~~~~~~~~~~(\widetilde{{\bf u}}_m,\widetilde{\phi}_m)~~\mapsto ({\bf u}_m,\phi_m).
\end{eqnarray*}
In the following, we describe properties of the map ${\cal T}^2_m.$ First, multiplying $\eqref{E3-5}_1$ by $a_{l}^{m}(t)$ and summing over $l$, one finds that
\begin{align*}
	&\frac{1}{2}\int_{\Omega}\rho_{m}\partial_{t}|{\mathbf u}_{m}|^{2}dx+\frac{1}{2}\int_{\Omega}\rho_{m}\widetilde{{\mathbf u}}_{m}\nabla|{\mathbf u}_{m}|^{2}dx
+\int_{\Omega}\nu(\widetilde{\phi}_{m})(1+|{\mathbb D}{\mathbf u}_m|^{2})^{\frac{p-2}{2}}|\mathbb{D}{\mathbf u}_{m}|^{2}dx\\
	&=\int_{\Omega}\rho_{m}\mu_{m}\nabla\widetilde{\phi}_{m}\cdot {\mathbf u}_{m}dx-\int_{\Omega}\rho_{m}{\mathbf u}_{m}\cdot\nabla\Psi(\widetilde{\phi}_{m})dx
\end{align*}
Multiplying \eqref{E3-4-1} by $\frac{1}{2}|{\mathbf u}_{m}|^{2}$ and integrating it over $\Omega,$ one gets that
\begin{equation*}
	\frac{1}{2}\int_{\Omega}\left(\partial_{t}\rho_{m}|{\mathbf u}_{m}|^{2}+|{\mathbf u}_{m}|^{2}\widetilde{{\mathbf u}}_{m}\cdot\nabla\rho_{m}\right)dx=0.
\end{equation*}
So,
\begin{align}
	&\frac{1}{2}\frac{d}{dt}\int_{\Omega}\rho_{m}|{\mathbf u}_{m}|^{2}dx+\int_{\Omega}\nu(\widetilde{\phi}_{m})(1+|{\mathbb D}{\mathbf u}_m|^{2})^{\frac{p-2}{2}}|\mathbb{D}{\mathbf u}_{m}|^{2}dx\notag\\
	&=\int_{\Omega}\rho_{m}\mu_{m}\nabla\widetilde{\phi}_{m}\cdot {\mathbf u}_{m}dx-\int_{\Omega}\rho_{m}{\mathbf u}_{m}\cdot\nabla\Psi(\widetilde{\phi}_{m})dx\label{E3-9}
\end{align}
One multiplies $\eqref{E3-5}_2$ by $c_{l}^{m}(t)$ and sums over $l$ to arrive at
\begin{align*}
	\int_{\Omega}\rho_{m}\partial_t\phi_{m}\mu_{m}dx+\int_{\Omega}\rho_{m}\mu_{m}{\mathbf u}_{m}\cdot\nabla\widetilde{\phi}_{m}dx
+\int_{\Omega}|\nabla\mu_{m}|^{2}dx=0.
\end{align*}
And, one multiplies $\eqref{E3-5}_3$ by $\frac{d}{dt}b_{l}^{m}(t)$ and sums over $l$ to get that
\begin{align*}
	\int_{\Omega}\rho_{m}\mu_{m}\partial_t\phi_{m}dx&=\int_{\Omega}\nabla\phi_{m}(\nabla\phi_{m})_{t}dx
+\int_{\Omega}\rho_{m}\Psi^\prime(\phi_{m})\partial_t\phi_{m}dx\\
&=\frac{1}{2}\frac{d}{dt}\int_{\Omega}|\nabla\phi_{m}|^{2}dx+\int_{\Omega}\rho_{m}(\Psi(\phi_{m}))_{t}dx
\end{align*}
It is easy to obtain that
\begin{align*}
	\frac{1}{2}\frac{d}{dt}\int_{\Omega}|\nabla\phi_{m}|^{2}dx
+\int_{\Omega}\rho_{m}(\Psi(\phi_{m}))_{t}dx+\int_{\Omega}\rho_{m}\mu_{m}\nabla\widetilde{\phi}_{m}\cdot {\mathbf u}_{m}dx+\int_{\Omega}|\nabla\mu_{m}|^{2}dx=0
\end{align*}
and so
\begin{align}
	&\frac{d}{dt}\int_{\Omega}\left(\frac{1}{2}|\nabla\phi_{m}|^{2}+\rho_{m}\Psi(\phi_{m})\right)dx\notag\\
	&=-\int_{\Omega}\rho_{m}\mu_{m}\nabla\widetilde{\phi}_{m}\cdot {\mathbf u}_{m}dx-\int_{\Omega}|\nabla\mu_{m}|^{2}dx
    +\int_{\Omega}\partial_{t}\rho_{m}\Psi(\phi_{m})dx.\label{{E3-10}}
\end{align}
Combining \eqref{E3-9} and \eqref{{E3-10}}, one gets that
\begin{eqnarray*}
	&&\frac{d}{dt}\int_{\Omega}\left(\frac{1}{2}\rho_{m}|{\mathbf u}_{m}|^{2}+\frac{1}{2}|\nabla\phi_{m}|^{2}+\rho_{m}\Psi(\phi_{m})\right)dx\\
   &&~~~~ +\int_{\Omega}\left(\nu(\widetilde{\phi}_{m})(1+|{\mathbb D}{\mathbf u}_m|^{2})^{\frac{p-2}{2}}|\mathbb{D}{\mathbf u}_{m}|^{2}+|\nabla\mu_{m}|^{2}\right)dx\nonumber\\
    &&=\int_{\Omega}\left(\rho_{m}{\mathbf u}_{m}\cdot\nabla\Psi(\widetilde{\phi}_{m})-\rho_{m}\widetilde{{\mathbf u}}_{m}\cdot\nabla\Psi(\phi_{m})\right)dx.
\end{eqnarray*}
Note that
\begin{eqnarray*}
&&|\int_{\Omega}\rho_{m}{\mathbf u}_{m}\cdot\nabla\Psi(\widetilde{\phi}_{m})dx|
=|\int_{\Omega}\rho_{m}(\widetilde{\phi}_{m}^2-1)\widetilde{\phi}_{m}{\mathbf u}_{m}\cdot\nabla\widetilde{\phi}_{m}dx|\\
&&\leqslant \frac12\int_{\Omega}\rho_{m}|{\mathbf u}_{m}|^2dx+C(\|\widetilde{\phi}_{m}\|_{L^\infty(\Omega)}^6+1)\|\rho_m\|_{L^\infty(\Omega)}\|\widetilde{\phi}_{m}\|^2_{L^2(\Omega)}
\end{eqnarray*}
and
\begin{eqnarray*}
&&|\int_{\Omega}\rho_{m}\widetilde{{\mathbf u}}_{m}\cdot\nabla\Psi(\phi_{m})dx|
=|\int_{\Omega}(\widetilde{{\mathbf u}}_{m}\cdot\nabla(\rho_{m}\Psi(\phi_{m}))-\widetilde{{\mathbf u}}_{m}\cdot\nabla\rho_{m}\Psi(\phi_{m}))dx|\\
&&\leqslant \|\widetilde{{\mathbf u}}_{m}\|_{{\bf L}^\infty(\Omega)}\|\frac{1}{\rho_{m}}\|_{L^\infty(\Omega)}\int_{\Omega}\rho_{m}\Psi(\phi_{m})dx
\end{eqnarray*}
Thus,
\begin{eqnarray}
	&&\frac{d}{dt}\int_{\Omega}\left(\frac{1}{2}\rho_{m}|{\mathbf u}_{m}|^{2}+\frac{1}{2}|\nabla\phi_{m}|^{2}+\rho_{m}\Psi(\phi_{m})\right)dx
    +\int_{\Omega}\left(\nu(\widetilde{\phi}_{m})(1+|{\mathbb D}{\mathbf u}_m|^{2})^{\frac{p-2}{2}}|\mathbb{D}{\mathbf u}_{m}|^{2}+|\nabla\mu_{m}|^{2}\right)dx\nonumber\\
 &&\leqslant
  \left(\|\widetilde{{\mathbf u}}_{m}\|_{L^\infty(\Omega)}\|\frac{1}{\rho_{m}}\|_{L^\infty(\Omega)}+1 \right)
   \left(\frac12\int_{\Omega}\rho_{m}|{\mathbf u}_{m}|^2dx+\int_{\Omega}\rho_{m}\Psi(\phi_{m})dx\right)\nonumber\\
    &&+C(\|\widetilde{\phi}_{m}\|_{L^\infty(\Omega)}^6+1)\|\rho_m\|_{L^\infty(\Omega)}\|\widetilde{\phi}_{m}\|^2_{L^2(\Omega)}.
\end{eqnarray}
Thanks to Gronwall's inequality, one finds that
\begin{eqnarray}
&&\left(\int_{\Omega}\left(\frac{1}{2}\rho_{m}|{\mathbf u}_{m}|^{2}+\frac{1}{2}|\nabla\phi_{m}|^{2}+\rho_{m}\Psi(\phi_{m})\right)dx\right)(s)\\
&&\leqslant e^{\int_0^t\left(m\|\widetilde{{\mathbf u}}_{m}\|_{{\bf L}^\infty(\Omega)}+1 \right)d\tau}\cdot\nonumber\\
&&\left(\int_{\Omega}\left(\frac{1}{2}\rho_{0\delta}|{\mathbf u}_{0m}|^{2}+\frac{1}{2}|\nabla\phi_{0m}|^{2}+\rho_{m}\Psi(\phi_{0m})\right)dx
    +C\int_0^t\left( (\|\widetilde{\phi}_{m}\|_{L^\infty(\Omega)}^6+1)(\rho^*+1)\|\widetilde{\phi}_{m}\|^2_{L^2(\Omega)}\right)d\tau\right)\nonumber\\
&&\leqslant e^{(m\widetilde{M}+1)t}\cdot \left(\int_{\Omega}\left(\frac{1}{2}\rho_{0\delta}|{\mathbf u}_{0m}|^{2}+\frac{1}{2}|\nabla\phi_{0m}|^{2}+\rho_{0\delta}\Psi(\phi_{0m})\right)dx
    +C( (\widetilde{M}^6+1)(\rho^*+1)\widetilde{M}^2)t\right)\nonumber\\
&&\leqslant e^{(m\widetilde{M}+1)t}\cdot \left(\int_{\Omega}\left(\frac{1}{2}\rho_{0}|{\mathbf u}_{0}|^{2}+\frac{1}{2}|\nabla\phi_{0}|^{2}+\rho_{0}\Psi(\phi_{0})\right)dx+1\right)
    +C( (\widetilde{M}^6+1)(\rho^*+1)\widetilde{M}^2)t\nonumber
\end{eqnarray}
for any $t\in (0,t_m].$ Select $T_0\in (0,\min\{1,t_m\})$ such that
$$e^{(m\widetilde{M}+1)t}\cdot \left(\int_{\Omega}\left(\frac{1}{2}\rho_{0}|{\mathbf u}_{0}|^{2}+\frac{1}{2}|\nabla\phi_{0}|^{2}+\rho_{0}\Psi(\phi_{0})\right)dx+1\right)
    +C( (\widetilde{M}^6+1)(\rho^*+1)\widetilde{M}^2)t\leqslant \widetilde{M}.$$
Note that $$\rho_{m}(t)\Psi(\phi_{m}(t))\geqslant \rho_{m}(\frac{1}{8}|\phi_{m}(t)|^{4}-\frac{1}{4})\geqslant \frac{1}{8}\rho_{m}|\phi_{m}(t)|^{4}-\frac{1}{4}\rho_{m}.$$
One gets the following uniform estimates
\begin{eqnarray}
\begin{cases}\label{E-3-16}
    \|\sqrt{\rho_{m}}{\mathbf u}_{m}\|_{{\bf L}^{\infty}(0,T_{0};{\bf L}^{2}(\Omega))}\leqslant C(\widetilde{M}),~~
    \|\nabla\phi_{m}\|_{L^{\infty}(0,T_{0};L^{2}(\Omega))}\leqslant C(\widetilde{M}),\\
    \|\sqrt[4]{\rho_{m}}\phi_{m}\|_{L^{\infty}(0,T_{0};L^{4}(\Omega))}\leqslant C(\widetilde{M}),\\
	\|\nabla\mu_{m}\|_{L^{2}(0,T_{0};L^{2}(\Omega))}\leqslant C(\widetilde{M}),~~
    \|(1+|\mathbb{D}{\mathbf u}_{m}|^{2})^{\frac{p-2}{4}}\mathbb{D}{\mathbf u}_{m}\|_{{\bf L}^{\infty}(0,T_{0};{\bf L}^{2}(\Omega))}\leqslant C(\widetilde{M}).
\end{cases}
\end{eqnarray}
Further, it follows from
\begin{align*}
  \int_{\Omega}\rho_{m}|\phi_{m}|^{2}dx&
  =\int_{\Omega}\sqrt{\rho_{m}}\sqrt{\rho_{m}}|\phi_{m}|^{2}dx
  \leqslant(\int_{\Omega}\rho_{m}dx)^{\frac{1}{2}}(\int_{\Omega}\rho_{m}|\phi_{m}|^{4}dx)^{\frac{1}{2}}
   \leqslant C(\widetilde{M})
\end{align*}
that
\begin{equation*}
	\|\sqrt{\rho_{m}}\phi_{m}\|_{L^{\infty}(0,T_{0};L^{2}(\Omega))}\leqslant C(\widetilde{M}).
\end{equation*}
Recalling that $\rho_{m}$ is bounded below and the Korn inequality, one gets that
$$\|{\mathbf u}_{m}\|_{{\bf L}^{\infty}(0,T_{0};{\bf L}^{2}(\Omega))}\leqslant C(\widetilde{M})~~\mbox{ and }~~
\|\phi_{m}\|_{L^{\infty}(0,T_{0};L^{2}(\Omega))}\leqslant C(\widetilde{M}).$$
So, we consider the following set
\begin{equation*}
	{\cal S}_{m}=\left\{({\bf v}_m,\psi_m)\in C([0,T_0];\mathbf{V}_{m})\times C([0,T_0];V_{m}):
    \underset{0\leqslant t\leqslant T_0}{max}\|{\bf v}_{m}(t)\|_{L^{2}(\Omega)}
          +\underset{0\leqslant t\leqslant T_0}{max}\|\psi_{m}(t)\|_{H^{1}(\Omega)}\leqslant C_{0}\right\}.
\end{equation*}
Due to \eqref{E3-5}, one obtains that
\begin{eqnarray*}
	&&\int_{\Omega}\rho_{m}|\partial_{t}{\mathbf u}_{m}|^{2}dx
     +\int_{\Omega}\rho_{m}(\widetilde{{\mathbf u}}_{m}\cdot\nabla) {\mathbf u}_{m}\cdot\partial_{t}{\mathbf u}_{m}dx
     +\int_{\Omega}\nu(\widetilde{\phi}_{m})(1+|{\mathbb D}{\mathbf u}_{m}|^{2})^{\frac{p-2}{2}}\mathbb{D}{\mathbf u}_{m}:\mathbb{D}\partial_{t} {\mathbf u}_{m}dx\\
	&&=\int_{\Omega}\rho_{m}\mu_{m}\nabla\widetilde{\phi}_{m}\cdot\partial_{t} {\mathbf u}_{m}dx-\int_{\Omega}\rho_{m}\nabla\Psi(\phi_{m})\cdot\partial_{t} {\mathbf u}_{m}dx,
\end{eqnarray*}
and so
\begin{eqnarray}\label{E-3-18}
	&&\int_{\Omega}\rho_{m}|\partial_{t}{\mathbf u}_{m}|^{2}dx\\
	&&\leqslant\|\sqrt{\rho_{m}}\|_{L^{\infty}(\Omega)}\|\sqrt{\rho_{m}}\partial_{t}{\mathbf u}_{m}\|_{{\bf L}^{2}(\Omega)}
        \|\widetilde{{\mathbf u}}_{m}\|_{{\bf L}^{6}(\Omega)}\|\nabla {\mathbf u}_{m}\|_{{\bf L}^{3}(\Omega)}\nonumber\\
    &&~~+\|\sqrt{\rho_{m}}\mu_{m}\|_{L^{2}(\Omega)}\|\sqrt{\rho_{m}}\|_{L^{\infty}(\Omega)}
     \|\partial_{t}{\mathbf u}_{m}\|_{{\bf L}^{3}(\Omega)}\|\nabla\widetilde{\phi}_{m}\|_{L^{6}(\Omega)}\nonumber\\
	&&~~+\|\sqrt{\rho_{m}}\|_{L^{\infty}(\Omega)}
     \|\sqrt{\rho_{m}}\partial_{t}{\mathbf u}_{m}\|_{{\bf L}^{2}(\Omega)}\|\nabla\phi_{m}\|_{L^{6}(\Omega)}(\|\phi_{m}^{3}\|_{L^{3}(\Omega)}+\|\phi_{m}\|_{L^{3}(\Omega)})\nonumber\\
   &&~~+\nu^{*}\|\nabla\partial_{t}{\mathbf u}_{m}\|_{{\bf L}^{2}(\Omega)}
       \|(1+|\mathbb{D}\widetilde{{\mathbf u}}_{m}|^{2})^{\frac{p-2}{2}}\mathbb{D}{\mathbf u}_{m}\|_{{\bf L}^{2}(\Omega)}\nonumber\\
	&&\leqslant C(\rho^{*})\|\sqrt{\rho_{m}}\partial_{t}{\mathbf u}_{m}\|_{{\bf L}^{2}(\Omega)}
         \|\widetilde{{\mathbf u}}_{m}\|_{{\bf H}^{2}(\Omega)}
          \|{\mathbf u}_{m}\|_{\mathbf{V}_{\sigma}}
           +C(\rho^{*})\|\sqrt{\rho_{m}}\mu_{m}\|_{L^{2}(\Omega)}
            \|\nabla\partial_{t}{\mathbf u}_{m}\|_{{\bf L}^{2}(\Omega)}\|\widetilde{\phi}_{m}\|_{H^{2}(\Omega)}\nonumber\\
	&&~~+C(\rho^{*})\|\sqrt{\rho_{m}}\partial_{t}{\mathbf u}_{m}
    \|_{{\bf L}^{2}(\Omega)}\|\phi_{m}\|_{H^{2}(\Omega)}(\|\phi_{m}\|^{3}_{H^{2}(\Omega)}
    +\|\phi_{m}\|_{H^{1}(\Omega)})\nonumber\\
	&&~~+\nu^{*}\|\partial_{t}{\mathbf u}_{m}\|_{\mathbf{V}_{\sigma}}
       \|(1+|\mathbb{D}{\mathbf u}_{m}|^{2})^{\frac{p-2}{2}}\mathbb{D}{\mathbf u}_{m}\|_{{\bf L}^{2}(\Omega)}.\nonumber
\end{eqnarray}
Note that
\begin{equation*}
\|\widetilde{{\mathbf u}}_{m}\|_{{\bf H}^{2}(\Omega)}\leqslant\lambda^{S}_{m}\|\widetilde{{\mathbf u}}_{m}\|_{{\bf L}^{2}(\Omega)},~~
\|{\mathbf u}_{m}\|_{\mathbf{V}_{\sigma}}\leqslant\sqrt{\lambda^{S}_{m}}\|{\mathbf u}_{m}\|_{{\bf L}^{2}(\Omega)},~~
\|\partial_{t}{\mathbf u}_{m}\|_{\mathbf{V}_{\sigma}}\leqslant\sqrt{\lambda^{S}_{m}}\|\partial_{t}{\mathbf u}_{m}\|_{{\bf L}^{2}(\Omega)},
\end{equation*}
where the $\lambda^{S}_{m}$ are the eigenvalues of the Stokes operator $\mathbf{A}$, and
\begin{equation*}
\|\widetilde{\phi}_{m}\|_{H^{2}(\Omega)}\leqslant\sqrt{\lambda_{m}}\|\widetilde{\phi}_{m}\|_{H^{1}(\Omega)},~~
\|\phi_{m}\|_{H^{2}(\Omega)}\leqslant\sqrt{\lambda_{m}}\|\phi_{m}\|_{H^{1}(\Omega)}.
\end{equation*}
Moreover,
\begin{align*}
	\int_{\Omega}\rho_{m}\mu_{m}^{2}dx&
    =-\int_{\Omega}\Delta\phi_{m}\mu_{m}dx+\int_{\Omega}\rho_{m}\mu_{m}(\phi_{m}^{3}-\phi_{m})dx\\
	&\leqslant\frac{1}{C\delta}\|\Delta\phi_{m}\|_{L^{2}(\Omega)}\|\sqrt{\rho_{m}}\mu_{m}\|_{L^{2}(\Omega)}
    +\|\rho_{m}\|^{\frac{1}{2}}_{L^{\infty}(\Omega)}\|\sqrt{\rho_{m}}\mu_{m}\|_{L^{2}(\Omega)}(\|\phi_{m}^{3}\|_{L^{2}(\Omega)}
    +\|\phi_{m}\|_{L^{2}(\Omega)})\\
	&\leqslant C(\sqrt{\lambda_{m}},\delta)\|\phi_{m}\|_{H^{1}(\Omega)}\|\sqrt{\rho_{m}}\mu_{m}\|_{L^{2}(\Omega)}
     +C\|\sqrt{\rho_{m}}\mu_{m}\|_{L^{2}(\Omega)}(\|\phi_{m}\|^{3}_{H^{1}(\Omega)}+\|\phi_{m}\|_{L^{2}(\Omega)})\\
	&\leqslant C(\sqrt{\lambda_{m}},\delta)\|\phi_{m}\|^{2}_{H^{1}(\Omega)}
     +\frac12\|\sqrt{\rho_{m}}\mu_{m}\|^{2}_{L^{2}(\Omega)}+C(\|\phi_{m}\|^{6}_{H^{1}(\Omega)}+\|\phi_{m}\|^{2}_{L^{2}(\Omega)}),
\end{align*}
and then
\begin{equation*}
	\|\sqrt{\rho_{m}}\mu_{m}\|_{L^{\infty}(0,T;L^{2}(\Omega))}\leqslant C(\sqrt{\lambda_{m}},\rho_{*}).
\end{equation*}
It follows from \eqref{E-3-18} that there exists $\tilde{C_{1}}=\tilde{C_{1}}(\delta,\rho^{*},\lambda^{S}_{m},\lambda_{m})$ such that
\begin{equation*}
	\underset{0\leqslant t\leqslant T}{max}\|\sqrt{\rho_{m}}\partial_{t}{\mathbf u}_{m}\|_{{\bf L}^{2}(\Omega)}\leqslant \tilde{C_{1}}
\end{equation*}
At the same time, one gets that
\begin{equation*}
	\int_{\Omega}\rho_{m}|\partial_{t}\phi_{m}|^{2}dx
     +\int_{\Omega}\rho_{m}({\mathbf u}_{m}\cdot\nabla)\widetilde{\phi}_{m}\partial_{t}\phi_{m}dx
     +\int_{\Omega}\nabla\mu_{m}\cdot\nabla\partial_{t}\phi_{m}dx=0,
\end{equation*}
and
\begin{align*}
	\int_{\Omega}\rho_{m}|\partial_{t}\phi_{m}|^{2}dx&
    \leqslant\|\sqrt{\rho_{m}}\partial_{t}\phi_{m}\|_{L^{2}(\Omega)}
    \|\sqrt{\rho_{m}}{\mathbf u}_{m}\|_{{\bf L}^{2}(\Omega)}\|\nabla\widetilde{\phi}_{m}\|_{L^{\infty}(\Omega)}
    +\|\nabla\mu_{m}\|_{L^{2}(\Omega)}\|\nabla\partial_{t}\phi_{m}\|_{L^{2}(\Omega)}\\
	&\leqslant\|\sqrt{\rho_{m}}\partial_{t}\phi_{m}\|_{L^{2}(\Omega)}
    \|\sqrt{\rho_{m}}{\mathbf u}_{m}\|_{{\bf L}^{2}(\Omega)}\|\widetilde{\phi}_{m}\|_{H^{2}(\Omega)}
     +\|\nabla\mu_{m}\|_{L^{2}(\Omega)}\|\nabla\partial_{t}\phi_{m}\|_{L^{2}(\Omega)}\\
	&\leqslant C\|\sqrt{\rho_{m}}\partial_{t}\phi_{m}\|_{L^{2}(\Omega)}
     \|\sqrt{\rho_{m}}{\mathbf u}_{m}\|_{{\bf L}^{2}(\Omega)}\sqrt{\lambda_{m}}\|\widetilde{\phi}_{m}\|_{H^{1}(\Omega)}
      +\lambda_{m}\|\mu_{m}\|_{L^{2}(\Omega)}\|\partial_{t}\phi_{m}\|_{L^{2}(\Omega)}\\
	&\leqslant C(\lambda_{m})\|\sqrt{\rho_{m}}\partial_{t}\phi_{m}\|_{L^{2}(\Omega)}\|\widetilde{\phi}_{m}\|_{H^{1}(\Omega)}
    +C(\lambda_{m},\rho_{*})\|\sqrt{\rho_{m}}\mu_{m}\|_{L^{2}(\Omega)}\|\sqrt{\rho_{m}}\partial_{t}\phi_{m}\|_{L^{2}(\Omega)}\\
	&\leqslant\frac12\|\sqrt{\rho_{m}}\partial_{t}\phi_{m}\|^{2}_{L^{2}(\Omega)}
    +C(\lambda_{m})\|\widetilde{\phi}_{m}\|^{2}_{H^{1}(\Omega)}+C(\lambda_{m})\|\sqrt{\rho_{m}}\mu_{m}\|_{L^{2}(\Omega)}.
\end{align*}
So, there exists $\tilde{C_{2}}=\tilde{C_{2}}(\lambda_{m},\lambda^{S}_{m})$ such that
\begin{equation*}
	\underset{0\leqslant t\leqslant T}{max}\|\partial_{t}\phi_{m}\|_{L^{2}(\Omega)}\leqslant \tilde{C_{2}},
\end{equation*}
Above all,
\begin{equation}\label{E-5-16}
	\underset{0\leqslant t\leqslant T}{max}\|\sqrt{\rho_{m}}\partial_{t}{\mathbf u}_{m}\|_{{\bf L}^{2}(\Omega)}\leqslant C_{1},~~~\underset{0\leqslant t\leqslant T}{max}\|\partial_{t}\phi_{m}\|_{L^{2}(\Omega)}\leqslant C_{2}.
\end{equation}
Set $\tilde{C_{0}}=C_{1}+C_{2}$ and consider the subset $\tilde{\cal S}_{m}$ of ${\cal S}_{m}$ defined as
\begin{align*}
	&\tilde{\cal S}_{m}=\left\{({\bf v}_m,\psi_m)\in C^{1}([0,T_0];{\bf V}_{m})\times C^{1}([0,T_0];V_{m}):\right.\\
      &\left.~~\underset{0\leqslant t\leqslant T_0}{max}\|{\bf v}_m(t)\|_{{\bf L}^{2}(\Omega)}
      +\underset{0\leqslant t\leqslant T_0}{max}\|\psi_m(t)\|_{H^{1}(\Omega)}\leqslant C_{0},
      \underset{0\leqslant t\leqslant T_0}{max}\|\partial_{t}{\bf v}_m(t)\|_{{\bf L}^{2}(\Omega)}
   +\underset{0\leqslant t\leqslant T_0}{max}\|\partial_{t}\psi_{m}(t)\|_{H^{1}(\Omega)}\leqslant \tilde{C_{0}}\right\}.
\end{align*}
Thus, the map ${\cal T}^2_m$ satisfies that
\begin{equation*}
	{\cal T}^2_m:~~S_{m}\mapsto \tilde{S}_{m}.
\end{equation*}
and $\tilde{\cal S}_{m}$ is compact in ${\cal S}_{m}$ due to the Ascoli-Arzel\`{a} theorem. Now, we focus on proving the continuity of the map ${\cal T}^2_m.$ Let us consider a sequence of $\{({\bf v}_m^n,\psi_m^n)\}_{n=1}^\infty\subset{\cal S}_m$ such that $({\bf v}_m^n,\psi_m^n)\rightarrow ({\bf v}_m,\psi_m)$ in $C([0,T_0],{\bf V}_m)\times C([0,T_0],V_m)$ as $n\rightarrow\infty.$
For fixed $({\bf v}_m^n,\psi_m^n)\in{\cal S}_m,$ the system
\begin{eqnarray*}\label{E-3-21}
\begin{cases}
\partial_{t}\rho_{m}+{\bf v}_m^n\cdot\nabla\rho_{m}=0,\\
(\rho_{m}\partial_{t}{\mathbf u}_{m},{\mathbf w}_l)+(\rho_{m}\left({\bf v}_m^n\cdot\nabla\right){\mathbf u}_{m},{\mathbf w}_l)
    +(\nu(\psi^n_{m})(1+|{\mathbb D}{{\mathbf u}}_{m}|^{2})^{\frac{p-2}{2}}\mathbb{D}{\mathbf u}_{m},\nabla {\mathbf w}_l)\\
~~~~~~~~~~~~~~~~=(\rho_{m}\mu_{m}\nabla\psi^n_{m},{\mathbf w}_l)-(\rho_{m}\nabla\Psi(\psi^n_{m}),{\mathbf w}_l),\\
(\rho_{m}\partial_{t}\phi_{m},w_l)+(\rho_{m} {\mathbf u}_{m}\cdot\nabla\psi^n_{m},w_l)+(\nabla\mu_{m},\nabla w_l)=0,\\
(\rho_{m}\mu_{m},w_l)=(\nabla\phi_{m},\nabla w_l)+(\rho_m\Psi^\prime(\phi_{m}),w_l),\\
div{\mathbf u}_{m}=0\\
{\mathbf u}_{m}=0,\partial_{n}\mu_{m}=\partial_{n}\phi_{m}=0, ~~~\mbox{ on }\partial\Omega\times (0,T),\\
\rho_{m}(\cdot,0)=\rho_{0\delta},~{\mathbf u}_{m}(\cdot,0)={\mathbf u}_{0m},~\phi_{m}(\cdot,0)=\phi_{0m}, ~~~\mbox{ in }\Omega
\end{cases}~(l=1,2,\cdots,m)
\end{eqnarray*}
admits a solution $(\rho_m^n, {\bf u}_m^n,\phi_m^n,\mu_m^n).$ One takes similar argument for \eqref{E-2}-\eqref{E-5-16} to deduce that
there exists $\tilde{K_{0}}$ and $\tilde{K_{1}},$ dependent on $m,{\rho}_{0m},{\mathbf u}_{0m}$ and $\phi_{0m}$ but independent on $n,$ such that $\tilde{K_{0}}\leqslant C_0$ and
\begin{align*}
    &\underset{0\leqslant t\leqslant T_0}{max}\|{\mathbf u}^{n}_{m}(t)\|_{{\bf L}^{2}(\Omega)}
    +\underset{0\leqslant t\leqslant T_0}{max}\|\phi^{n}_{m}(t)\|_{H^{1}(\Omega)}\leqslant \tilde{K_{0}},\\
	&\underset{0\leqslant t\leqslant T_0}{max}\|\partial_{t}{\mathbf u}^{n}_{m}(t)\|_{{\bf L}^{2}(\Omega)}
    +\underset{0\leqslant t\leqslant T_0}{max}\|\partial_{t}\phi^{n}_{m}(t)\|_{H^{1}(\Omega)}\leqslant \tilde{K_{1}}.
\end{align*}
So, the Ascoli-Arzel\`{a} theorem implies that, there exist a subsequence $({\mathbf u}^{n_{j}}_{m},\phi^{n_{j}}_{m})$ of $({\mathbf u}^{n}_{m},\phi^{n}_{m})$ and
$(\widetilde{{\bf u}},\widetilde{\phi})$ such that
\begin{align*}
	&{\mathbf u}^{n_{j}}_{m}\rightarrow \widetilde{{\bf u}} \;strongly\;in\;{\bf C}([0,T];{\mathbf V}_{m}),~~
   \phi^{n_{j}}_{m}\rightarrow \widetilde{\phi} \;strongly\;in\;C([0,T];V_{m}),\\
	&\partial_{t}{\mathbf u}^{n_{j}}_{m}\rightarrow \partial_{t}\widetilde{{\bf u}} \;weak-star\;in\;{\bf L}^{\infty}(0,T;{\mathbf V}_{m}),~~
    \partial_{t}\phi^{n_{j}}_{m}\rightarrow \partial_{t}\widetilde{\phi} \;weak-star\;in\;L^{\infty}(0,T;V_{m}).
\end{align*}
as $n_j\rightarrow \infty.$ Further, one follows the similar argument of \eqref{E3-12} to deduce  that there exist a subsequence $\rho_m^{n_j}$ of $\rho_m^n$ and $\widetilde{\rho}$ such that for any fixed $r\geqslant 2$
\begin{eqnarray*}
\rho_m^{n_j}\rightarrow \widetilde{\rho}   \;strongly\;in\;C([0,T];L^r(\Omega))~~\mbox{as}~~n_j\rightarrow \infty.
\end{eqnarray*}
Moreover, it is deduced from the equation $\rho^{n}_{m}\mu^{n}_{m}=-\Delta\phi^{n}_{m}+\rho^{n}_{m}\Psi^\prime(\phi^{n}_{m})$ that
\begin{equation*}
	\mu_{m}^{n_{j}}\rightarrow\widetilde{\mu}\;weakly\;in\;L^{2}(0,T;V_{m})~~\mbox{as}~~n_j\rightarrow \infty.
\end{equation*}
Thus, one obtains that
\begin{eqnarray*}\label{E-3-21}
\begin{cases}
\partial_{t}{\widetilde \rho}+{\bf v}_m\cdot\nabla{\widetilde \rho}=0,\\
({\widetilde \rho}\partial_{t}{\widetilde {\bf u}},{\mathbf w}_l)+({\widetilde \rho}\left({\bf v}_m\cdot\nabla\right){\widetilde {\bf u}},{\mathbf w}_l)
    +(\nu(\psi_m)(1+|{\mathbb D}{{\mathbf u}}_{m}|^{2})^{\frac{p-2}{2}}\mathbb{D}{\widetilde {\bf u}},\nabla {\mathbf w}_l)\\
~~~~~~~~~~~~~~~~=({\widetilde \rho}{\widetilde \mu}\nabla\psi_m,{\mathbf w}_l)-({\widetilde \rho}\nabla\Psi(\psi_m),{\mathbf w}_l),\\
({\widetilde \rho}\partial_{t}{\widetilde \phi},w_l)+({\widetilde \rho} {\widetilde {\bf u}}\cdot\nabla\psi_m,w_l)+(\nabla{\widetilde \mu},\nabla w_l)=0,\\
({\widetilde \rho}{\widetilde \mu},w_l)=(\nabla{\widetilde \phi},\nabla w_l)+(\rho_m\Psi^\prime({\widetilde \phi}),w_l),\\
{\widetilde {\bf u}}=0,\partial_{n}{\widetilde \mu}=\partial_{n}{\widetilde \phi}=0, ~~~\mbox{ on }\partial\Omega\times (0,T),\\
{\widetilde \rho}(\cdot,0)=\rho_{0\delta},~{\widetilde {\bf u}}(\cdot,0)={\mathbf u}_{0m},~{\widetilde \phi}(\cdot,0)=\phi_{0m}, ~~~\mbox{ in }\Omega
\end{cases}
\end{eqnarray*}
for $l=1,2,\cdots,m$ and for all $t\in [0,T_0].$ On the other hand, for $({\bf v}_m,\psi_m),$ the system
\begin{eqnarray*}
\begin{cases}
\partial_{t}\rho_{m}+{\bf v}_m\cdot\nabla\rho_{m}=0,\\
(\rho_{m}\partial_{t}{\mathbf u}_{m},{\mathbf w}_l)+(\rho_{m}\left({\bf v}_m\cdot\nabla\right){\mathbf u}_{m},{\mathbf w}_l)
    +(\nu(\psi_{m})(1+|{\mathbb D}{{\mathbf u}}_{m}|^{2})^{\frac{p-2}{2}}\mathbb{D}{\mathbf u}_{m},\nabla {\mathbf w}_l)\\
~~~~~~~~~~~~~~~~=(\rho_{m}\mu_{m}\nabla\psi_{m},{\mathbf w}_l)-(\rho_{m}\nabla\Psi(\psi_{m}),{\mathbf w}_l),\\
(\rho_{m}\partial_{t}\phi_{m},w_l)+(\rho_{m} {\mathbf u}_{m}\cdot\nabla\psi_{m},w_l)+(\nabla\mu_{m},\nabla w_l)=0,\\
(\rho_{m}\mu_{m},w_l)=(\nabla\phi_{m},\nabla w_l)+(\rho_m\Psi^\prime(\phi_{m}),w_l),\\
{\mathbf u}_{m}=0,\partial_{n}\mu_{m}=\partial_{n}\phi_{m}=0, ~~~\mbox{ on }\partial\Omega\times (0,T),\\
\rho_{m}(\cdot,0)=\rho_{0\delta},~{\mathbf u}_{m}(\cdot,0)={\mathbf u}_{0m},~\phi_{m}(\cdot,0)=\phi_{0m}, ~~~\mbox{ in }\Omega
\end{cases}~~(l=1,2,\cdots,m)
\end{eqnarray*}
admits a unique solution $(\rho^0,{\bf u}^0,\phi^0,\mu^0)$ on $[0,T_0].$ It follows from uniqueness that
$$(\widetilde{\rho},\widetilde{{\bf u}},\widetilde{\phi},\widetilde{\mu})=(\rho^0,{\bf u}^0,\phi^0,\mu^0).$$
Hence, every subsequence of $(\rho_m^n, {\bf u}_m^n,\phi_m^n,\mu_m^n)$ possesses a subsequence whose limit coincide with
$(\widetilde{\rho},\widetilde{{\bf u}},\widetilde{\phi},\widetilde{\mu}).$
Therefore, the map ${\cal T}_m^2: {\cal S}_{m}\mapsto{\cal S}_{m}$ is continuous satisfying that $\overline{{\cal T}_m^2({\cal S}_{m})}=\widetilde{{\cal S}}_{m}$ is compact in $C([0,T];{\bf V}_m)\times C([0,T];V_m).$

According to the Schauder fixed point theorem, there exists $({\bf u}_m,\phi_m)\in C([0,T];{\bf V}_m)\times C([0,T];V_m)$ such that ${\cal T}_m^2({\bf u}_m,\phi_m)=({\bf u}_m,\phi_m).$ Furthermore, $\rho_m$ and $\mu_m$ can be determined by $({\bf u}_m,\phi_m).$ This yields a local solution $(\rho_m,{\bf u}_m,\phi_m,\mu_m)$ to the following semi-Galerkin scheme
\begin{eqnarray}\label{E3-3}
\begin{cases}
\partial_{t}\rho_{m}+{\mathbf u}_{m}\cdot\nabla\rho_{m}=0&\mbox{in}\;\Omega\times(0,T),\\
(\rho_{m}\partial_{t}{\mathbf u}_{m},{\mathbf w})+(\rho_{m}{\mathbf u}_{m}\cdot\nabla{\mathbf u}_{m},{\mathbf w})
    +(\nu(\phi_{m})(1+|{\mathbb D}{{\mathbf u}}_{m}|^{2})^{\frac{p-2}{2}}\mathbb{D}{\mathbf u}_{m},\nabla {\mathbf w})\\
~~~~~~~~~~~~~~~~=(\rho_{m}\mu_{m}\nabla\phi_{m},{\mathbf w})-(\rho_{m}\nabla\Psi(\phi_{m}),{\mathbf w})~~~~({\bf w}\in {\bf V}_m)&\mbox{in}\;(0,T),\\
(\rho_{m}\partial_{t}\phi_{m},w)+(\rho_{m} {\mathbf u}_{m}\cdot\nabla\phi_{m},w)+(\nabla\mu_{m},\nabla w)=0~~~~(w\in V_m)&\mbox{in}\;(0,T)\\
(\rho_{m}\mu_{m},w)=(\nabla\phi_{m},\nabla w)+(\rho_m\Psi^\prime(\phi_{m}),w)~~~~(w\in V_m)&\mbox{in}\;(0,T),\\
{\mathbf u}_{m}=0,\partial_{n}\mu_{m}=\partial_{n}\phi_{m}=0 &\mbox{ on }\partial\Omega\times (0,T),\\
\rho_{m}(\cdot,0)=\rho_{0\delta},{\mathbf u}_{m}(\cdot,0)={\mathbf u}_{0m},\phi_{m}(\cdot,0)=\phi_{0m}&\mbox{ in }\Omega.
\end{cases}
\end{eqnarray}
on the interval $[0,T_0].$ Uniqueness of the approximate solution $(\rho_m,{\bf u}_m,\phi_m,\mu_m)$ can be proved by the standard energy method and we omit the details here.

Since $(\rho_m,{\bf u}_m,\phi_m,\mu_m)$ is the unique solution to the semi-Galerkin scheme \eqref{E3-3} on the interval $[0,T_0],$ one finds that
\begin{eqnarray*}
   &&0<C\delta\leqslant\rho_{m}(x,t)\leqslant\rho^{*}+\frac{1}{m}\;\;\;(\forall(x,t)\in Q_{T}),\\
   &&\frac{d}{dt}\int_{\Omega}\left(\frac{1}{2}\rho_{m}|{\mathbf u}_{m}|^{2}+\frac{1}{2}|\nabla\phi_{m}|^{2}+\rho_{m}\Psi(\phi_{m})\right)dx\\
   &&+\int_{\Omega}\left(\nu({\phi}_{m})(1+|{\mathbb D}{\mathbf u}_m|^{2})^{\frac{p-2}{2}}|\mathbb{D}{\mathbf u}_{m}|^{2}+|\nabla\mu_{m}|^{2}\right)dx=0,
\end{eqnarray*}
and so
\begin{eqnarray}\label{E3-22}
\begin{cases}
     0<C\delta\leqslant\rho_{m}(x,t)\leqslant\rho^{*}+1\;\;\;(\forall(x,t)\in Q_{T}),\\
  \|\sqrt{\rho_{m}}{\partial_{t}\mathbf u}_{m}\|_{{\bf L}^{2}(0,T_{0};{\bf L}^{2}(\Omega))}\leqslant C,~~
 \|\sqrt{\rho_{m}}\partial_{t}\phi_{m}\|_{L^{2}(0,T_{0};L^{2}(\Omega))}\leqslant C,\\
    \|\sqrt{\rho_{m}}{\mathbf u}_{m}\|_{{\bf L}^{\infty}(0,T_{0};{\bf L}^{2}(\Omega))}\leqslant C,~~
    \|\nabla\phi_{m}\|_{L^{\infty}(0,T_{0};L^{2}(\Omega))}\leqslant C,\\
    \|\sqrt[4]{\rho_{m}}\phi_{m}\|_{L^{\infty}(0,T_{0};L^{4}(\Omega))}\leqslant C,\\
	\|\nabla\mu_{m}\|_{L^{2}(0,T_{0};L^{2}(\Omega))}\leqslant C,~~
    \|(1+|\mathbb{D}{\mathbf u}_{m}|^{2})^{\frac{p-2}{4}}\mathbb{D}{\mathbf u}_{m}\|_{{\bf L}^{\infty}(0,T_{0};{\bf L}^{2}(\Omega))}\leqslant C.
\end{cases}
\end{eqnarray}
where $C$ is independent of $m$ and $\delta.$ Further, it follows from
\begin{align*}
  \int_{\Omega}\rho_{m}|\phi_{m}|^{2}dx&
  =\int_{\Omega}\sqrt{\rho_{m}}\sqrt{\rho_{m}}|\phi_{m}|^{2}dx
  \leqslant(\int_{\Omega}\rho_{m}dx)^{\frac{1}{2}}(\int_{\Omega}\rho_{m}|\phi_{m}|^{4}dx)^{\frac{1}{2}}
   \leqslant C,
\end{align*}
that
\begin{equation}\label{E3-23}
	\|\sqrt{\rho_{m}}\phi_{m}\|_{L^{\infty}(0,T_{0};L^{2}(\Omega))}\leqslant C.
\end{equation}
Since $\rho_m$ is bounded from below, the estimates \eqref{E3-22} and \eqref{E3-23} permit to extend the local solution $(\rho_m,{\bf u}_m,\phi_m,\mu_m)$ from $[0,T_0]$ to the whole interval $[0,T].$ This yields a unique global solution to the semi-Galerkin scheme  \eqref{E3-3}.
\end{proof}

\end{document}